\input amstex
\documentstyle{amsppt}
\magnification=\magstep1 \NoRunningHeads
\topmatter
\NoBlackBoxes

\title
Rank-one actions, their $(C,F)$-models and constructions with bounded parameters
 \endtitle

\author
Alexandre I. Danilenko 
\endauthor

\email
alexandre.danilenko@gmail.com
\endemail

\address
 Institute for Low Temperature Physics
\& Engineering of National Academy of Sciences of Ukraine, 47 Nauki Ave.,
 Kharkov, 61164, UKRAINE
\endaddress
\email alexandre.danilenko\@gmail.com
\endemail

\abstract
Let $G$ be a discrete countable infinite group. 
We show that each topological $(C,F)$-action $T$ of $G$ 
on a locally compact non-compact Cantor set
is a free minimal amenable action admitting  a unique up to scaling non-zero invariant Radon measure (answer to a question by Kellerhals, Monod and R{\o}rdam).
We find necessary and sufficient conditions  under which two such actions are topologically conjugate in terms of the underlying $(C,F)$-parameters.
If $G$ is linearly ordered Abelian then the topological centralizer of $T$ is trivial.
If $G$ is monotileable and amenable, denote by $\Cal A_G$ the set of all probability preserving actions of $G$ on the unit interval with Lebesgue measure and endow it with the natural topology.
We show that the set of  $(C,F)$-parameters of all $(C,F)$-actions of $G$ furnished with a suitable topology is a model for $\Cal A_G$ in the sense of Foreman, Rudolph and Weiss.
If $T$ is a rank-one transformation with  bounded sequences of cuts and spacer maps then we found simple necessary and sufficient conditions on the related $(C,F)$-parameters under which  (i) $T$ is rigid, (ii) $T$ is totally ergodic.
It is found an alternative proof of Ryzhikov's theorem that if $T$ is totally ergodic and non-rigid rank-one map with bounded parameters then $T$ has MSJ.
We also give a  more general version of the criterion (by Gao and Hill) for isomorphism and disjointness  of two commensurate  non-rigid totally ergodic rank-one maps with bounded parameters.
It is shown that the  rank-one transformations with bounded parameters and no spacers over the last subtowers is a proper subclass of the rank-one transformations with bounded parameters.
\endabstract

\endtopmatter

\document

 \loadbold

  \head 0. Introduction
 \endhead
 
 The original goal  of  this work was to find new, shorter and more explicit proofs of the main results from  recent papers by Gao and Hill \cite{GaHi1}--\cite{GaHi4}, \cite{Hi1} and \cite{Hi2} on topological and measure theoretical properties of  rank-one transformations.
 As appeared,  our approach  is  well suited  to extend most of those results   to more general classes of rank-one transformations and (some of the results)  to  rank-one actions of more general groups, including non-amenable ones\footnote{The exact relationship  between the original Gao-Hill results and our generalization of them is clarified below (in Introduction) in detail.}.
 This is achieved by applying---in our opinion---more natural, intrinsic techniques to the problems under consideration. 
 While Gao and Hill consider the rank-one transformations as shift-maps on invariant subsets of $\{0,1\}^\Bbb Z$ and study them via tools of symbolic dynamics, our approach is based on analysis of the standard  {\it cutting-and-stacking} constructing algorithm.
 Since this  classical geometric 
 algorithm looses its  clarity beyond the framework of $\Bbb Z^d$-actions, we utilize  instead of it 
 the $(C,F)$-{\it construction} (see  \cite{dJ}, \cite{Da5}, \cite{Da6}) which we consider  as an {\it arithmetic version}  of the cutting-and-stacking. 
 The construction  is  convenient to  produce and investigate  rank-one actions of {\it arbitrary} locally compact groups.

In  \S 1 we briefly review the $(C,F)$-construction and (in the case of $\Bbb Z$-actions) discuss  relation between the $(C,F)$-notions and the classical concepts of  the cutting-and-stacking.
The class of probability preserving $(C,F)$-actions of $\Bbb Z$  is up to isomorphism the class of $\Bbb Z$-actions of funny rank-one (see \cite{Fe1} for Thouvenot's definition of funny rank one).
Every such action is associated with two sequences of finite subsets in $\Bbb Z$:  a sequence of {\it tiling shapes} and a sequence of {\it tiling centers}. 
In a similar way, given an arbitrary  discrete countable  infinite group $G$ and two sequences of finite subsets $(C_n)_{n\ge 1}$ and $(F_n)_{n\ge 0}$ in $G$ satisfying some conditions (see \S 1.2), we associate  a minimal continuous action of $G$ on a perfect totally disconnected  Polish  space  equipped  with a canonical $\sigma$-finite invariant measure (see \S 1.3).
Under an additional condition (see Lemma~1.4(i)) the space is locally compact and the measure is Radon.
We use this  general definition  to answer 
affirmatively the following  non-trivial  questions in the theory of $C^*$-algebras and topological group actions (see \S 1.5):
\roster
\item"---" does $G$ admit a free minimal amenable (in the topological sense according to Definition~1.6 below) action on
a locally compact non-compact Cantor space $X$  \cite{KeMoR\o, Question~7.2(i)}\footnote{Compare with a recent paper \cite{MaR\o} devoted to a partial answer  to this question  in the case where $G$ is {\it exact}, i.e. $G$ admits an amenable action on a compact set.}?
\item"---"  does $G$ admit a free minimal amenable action on
 $X$ which leaves invariant a non-zero Radon measure on $X$ \cite{KeMoR\o, Question~7.2(ii)}?
\endroster

In \S 2 we  study  {\it topological properties} of continuous $(C,F)$-actions (or, equivalently, rank-one actions) of $G$ on locally compact Cantor spaces.
Our main concern is the topological classification of these actions.
Namely, we want to determine  when two $(C,F)$-actions of an {\it arbitrary} discrete countable  infinite group are topologically isomorphic in terms
of the underlying sequences  $(C_n)_{n\ge 1}$ and $(F_n)_{n\ge 0}$ viewed as the {\it parameters} of the actions.
We  solve this general problem in Theorem~2.3 as follows:

\proclaim{Theorem A}
Let $T=(T_g)_{g\in G}$ and $T'=(T'_g)_{g\in G}$ be two $(C,F)$-actions of $G$ on locally compact Cantor spaces associated with some sequences $(C_n,F_{n-1})_{n\ge 1}$ and $(C_n',F_{n-1}')_{n\ge 1}$ respectively.
Then $T$ and $T'$  are topologically isomorphic if and only if 
 there is an
 increasing sequence of integers
 $0=l_0<l_1'<l_1<l_2'<l_2<\cdots$ 
 and subsets
$A_n \subset  F_{l_n'}'$,  $B_n \subset  F_{l_n}$, 
such that $A_{n}B_{n}=C_{l_{n-1}+1}\cdots C_{l_{n}}$, $B_{n} A_{n+1}=C'_{l_{n}'+1}\cdots C'_{l_{n+1}'}$, $F_{l_n}A_{n+1}\subset F'_{l_{n+1}}$, $F_{l_n}A_{n+1}\subset F'_{l_{n+1}}$
and $(F'_{l_n'})^{-1}F'_{l_n'}\cap B_nB_n^{-1}=F_{l_n}^{-1}F_{l_n}\cap A_{n+1}A_{n+1}^{-1}=\{1\}$ for each $n>0$.
\endproclaim

Then we show that in some important special cases the general criterion from Theorem A takes simpler (or more ``compact'') forms.
For instance, if $G$ is 
{\it linearly ordered} Abelian  (Definition~2.4) and the actions are {\it commensurate},
which means that $F_n=F_n'$ eventually, then we obtain the following (see Theorem 2.8):

\proclaim{Theorem B}
Let $T$ and $T'$ be as above and let 
$(G,G_+)$ be a linearly  ordered discrete countable Abelian group. 
Suppose that $C_n\cup C_n'\subset G_+$  for all $n$ and $F_n=F_n'$ for all $n>N$ (for some $N>0$).
Then $T$ and $T'$ are
topologically isomorphic if and only if
 $C_n=C_n'$ for all $n>M$
(for some $M>0$).
\endproclaim

In this connection we mention that the class of linearly ordered Abelian groups is  huge.
Indeed, by  \cite{Le},  an Abelian group admits a linear order if and only if it is torsion free.

Another important case of Theorem~A is  when $G=\Bbb Z$.
In this case we show the following (see~Theorem~2.7):

\proclaim{Theorem C}  Let $T=(T_g)_{g\in G}$ and $T'=(T'_g)_{g\in G}$ be as in Theorem~A 
and let $G=\Bbb Z$.
Suppose that $C_n\cup C_n'\subset \Bbb Z_+$.
Then $T$ and $T'$ are topologically isomorphic if and only if 
there is  $r\ge 0$ and a subset $R\subset F_r'\cap\Bbb Z_+$ such that
$
\sum_{i>0} C_i = R+\sum_{i>r}C'_{i}.
$
\endproclaim

Given a continuous action $(T_g)_{g\in G}$ of $G$ on a topological space $X$, we call the group of homeomorphisms of $X$ commuting with each $T_g$, $g\in G$, the {\it topological centralizer} of $T$.
Denote it by $C_{\text{top}}(T)$.
In Corollary~2.9  we characterize  the topological centralizer of the $(C,F)$-actions of linearly ordered groups:

\proclaim{Theorem D} Let $(G,G_+)$ and $T$ be as in Theorem~B.
Then $C_{\text{top}}(T)=\{T_g\mid g\in G\}$.
\endproclaim

In the case  of topological rank-one $\Bbb Z$-actions we provide a satisfactory solution of the {\it inverse} problem (see Corollary~2.11 for a more general result): when an action is isomorphic to its inverse?

\proclaim{Theorem E} Let $(G,G_+)=(\Bbb Z,\Bbb Z_+)$ and let $T$ be as in Theorem~C.
If $F_n=\{-a_n,\dots,0,\dots,h_n-1\}$ for some $a_n,h_n>0$ and each $n\in\Bbb N$  then
 $T_1$ and $T_{-1}$ are topologically isomorphic if and only if $C_n=\{\max C_n-c\mid c\in C_n\}$ for all $n>M$
(for some $M>0$).
\endproclaim

In the particular case where $G=\Bbb Z$ and the actions are of topological rank-one (which means that $C_n\subset \Bbb Z_+$ and $F_n=\{-a_n,\dots,b_n\}$ for some $a_n,b_n>0$ and each $n\in\Bbb N$ with $\lim_{n\to\infty}a_n=\lim_{n\to\infty}b_n=+\infty$), Theorems~B, C, D and  E are close to main results from \cite{GaHi2} (and partly from \cite{Hi1}).
However there are two points of difference:
\roster
\item
The topological systems considered by Gao and Hill are defined on compact Cantor spaces while  $(C,F)$-actions from \S 2 are defined on locally compact non-compact Cantor spaces.
This seeming difference is eliminated easily by passing to the one-point compactification.
Then the $(C,F)$-actions extend to these compactifications as {\it almost minimal} continuous actions on compact Cantor spaces.
We recall that a continuous action is called {\it almost minimal \cite{Da2}} if it has a single fixed point and  the other orbits are dense.
\item
The second  difference is that Gao and Hill consider  symbolic  models of  the so-called {\it adapted} rank-one transformations.
This  means that no spacers are added on the top the last (highest) subtower on any step of the inductive cutting-and-stacking-construction.
When  constructing {\it topological}   $(C,F)$-models of  rank-one transformations we 
introduce different ``boundary conditions''  on spacers: we add  spacers  on infinitely many steps as on the top of the last subtower as under the bottom of the first  subtower.
\footnote{No specific ``boundary conditions''  on spacers are imposed when we consider {\it measure theoretical} $(C,F)$-transformations.}
 
\endroster

Now, apart from these differences, we compare Theorems~A--E with  main results of \cite{GaHi2}.
Theorem~C  is, in fact, a counterpart of  Theorem~1.1 from  \cite{GaHi2} because
the sums $\sum_{i>0}C_i$ and $\sum_{i>r}C_i'$ from the statement of Theorem~C are algebraic analogues of the combinatorial sets of the so-called  ``expected occurrences'' of certain symbolic words used implicitly  in the statement of  \cite{GaHi2, Theorem~1.1}.
The two theorems deal only with the case of $\Bbb Z$-actions.
However,  while \cite{GaHi2, Theorem~1.2} deals with the (classic) rank-one $\Bbb Z$-actions,
Theorem~C is valid for a much larger class of  {\it funny} rank-one $\Bbb Z$-actions. Theorem~D is an  extension of   from the class rank-one $\Bbb Z$-actions to the class of (funny) rank-one actions of arbitrary linearly ordered Abelian groups.
Both Theorem~E and \cite{GaHi2, Theorem~1.4} solve the inverse problem in the class of rank-one $\Bbb Z$-actions. 
An advantage of Theorem~E is that it provides the solution in terms of the original $(C,F)$-parameters, while \cite{GaHi2, Theorem~1.4} gives the solution only in terms of the so-called ``canonical generating sequence'' which thus adds an extra computational step of passing from the original parameters to ``canonical'' ones.
 Theorem~B is new even in the case of rank-one $\Bbb Z$-actions.
Theorem~A is a universal isomorphism criterion for the rank-one actions of arbitrary countable groups.
Theorems~B--E follow from it.
The conditions (countably many of them) for isomorphism in Theorem~A  look ``heavier'' than a single condition from the statement of Theorem~C or its analogue \cite{GaHi2, Theorem~1.1}.
It is however a price we pay for  the generality of Theorem~A.
Finally, we emphasize that every of these conditions in Theorem~A is {\it finitary}, i.e. involves only finitely many finite sets $C_i$, while the condition from Theorem~C (and  \cite{GaHi2, Theorem~1.1}) is {\it infinitary}, i.e. involves infinitely many parameters.  

In the rest of the paper (\S 3--\S 5) we study {\it measure theoretical properties} of the  $(C,F)$-actions according to the most general definition (in which 
 \thetag{1-4} holds instead of the more restrictive \thetag{1-3} considered in \S 2).
 Then the $(C,F)$-actions are defined on Polish  but not necessarily locally compact  spaces
 and they have a canonical ergodic invariant $\sigma$-finite measure. 
 Thus we regard  the $(C,F)$-actions  in \S 3--\S 5 as    standard {\it measure preserving} dynamical systems and study them by modulo measure theoretical isomorphism.
In particular, in the case where $G=\Bbb Z$, the adapted rank-one transformations considered by Gao and Hill in \cite{GaHi1}, \cite{GaHi3}, \cite{GaHi4} and \cite{Hi2} are all in this class of
$(C,F)$-systems.

In \S 3, $G$ is a {\it monotileable}  amenable group \cite{We}.
We note that each Abelian, nilpotent or solvable countable group is monotileable \cite{We}.
It is currently unclear whether or not every  countable amenable group is monotileable: no example of non-monotileable amenable group is known so far.
Let $\Cal  F$ be a  F{\o}lner sequence of finite sets that tile $G$.
We denote by $\Cal A_G$ the set of all Lebesgue measure preserving actions of $G$ on the unit interval $[0,1)$.
Endow this set with the natural (Polish) weak topology.
We  generalize the concept of  a {\it model} for $\Bbb Z$-actions in the sense of Foreman-Rudolph-Weiss \cite{Fo} to the case of $G$-actions~(Definition~3.5).

\definition{Definition} A {\it model} for  $\Cal A_G$  is a pair $(W,\pi)$, 
where $W$ is a Polish space and $\pi:W\to\Cal A_G$ is a continuous map such that for a comeager set $\Cal M\subset\Cal A_G$ and each $A\in\Cal M$, the set $\{w\in W\mid \pi(w)\text{ is isomorphic to } A\}$ is dense in $W$.
\enddefinition

Denote by $\goth R_1^{\text{fin}}$ the set of all possible $(C,F)$-parameters, i.e. sequences $(C_n)_{n\ge 1}$ and $(F_n)_{n\ge 1}$, satisfying the conditions for $(C,F)$-actions and for which the corresponding  canonical measure is finite. 
Thus $\goth R_1^{\text{fin}}$ is a subset of $(\goth F\times\goth F)^{\Bbb N}$, where $\goth F$ is the countable set of all finite subsets in $G$.
Then we introduce a certain Polish topology on $\goth R_1^{\text{fin}}$ which is stronger than the product topology inherited from $(\goth F\times\goth F)^{\Bbb N}$.
We also construct a continuous map $\Psi:\goth R_1^{\text{fin}}\ni \Cal S\mapsto\Psi(\Cal S)\in\Cal A_G$ such that
$\Psi(\Cal S)$ is isomorphic to the $(C,F)$-action associated with $\Cal S$ for each $\Cal S\in \goth R_1^{\text{fin}}$.
The following is the main result of \S 3 (see Proposition~3.6 and Corollary~3.7).

\proclaim{Theorem F} The subset of  $G$-actions which are of rank one  along $\Cal F$ is a dense $G_\delta$ in $\Cal A_G$.
The pair
$(\goth R_1^{\text{fin}},\Psi)$ is a model for $\Cal A_G$.
\endproclaim

 
In the particular case where $G=\Bbb Z$ and $\Cal F=\{[0,\dots,n)\mid n\in\Bbb N\}$, 
the first claim of Theorem~F is a folklore in ergodic theory; the second claim 
 is an alternative version of the main result from \cite{GaHi1}.
 Thus, Theorem~F is an extension of \cite{GaHi1} to the class of monotileable amenable groups.

In \S 4, we consider  rank-one $\Bbb Z$-actions with {\it bounded parameters}.
This means that the number of cuts and the total number of spacers added on each step of the inductive cutting-and-stacking construction  are both  bounded.
Equivalently, in the language of the $(C,F)$-construction, the parameters $(C_n)_{n\ge 1}$ and $(F_n)_{n\ge 0}$  of a $(C,F)$-action of $\Bbb Z$ are bounded if the sequence $(\#C_n)_{n\ge 1}$ is bounded
and there is a finite subset $K\subset \Bbb Z$ such that $K+F_n+C_{n+1}\supset F_{n+1}$ and $F_n\in \Cal F:=\{[0,\dots,n)\mid n\in\Bbb N\}$ for each $n\ge 0$.
We note that each rank-one $\Bbb Z$-action with bounded parameters 
is finite measure preserving.
The interest to such systems grew up after Bourgain's work \cite{Bo} where it was shown that they
 satisfy the M{\"o}bius orthogonality property\footnote{See also subsequent works \cite{EALedR}, \cite{Ry}.}.  
We now state   the main results of  \S 4 (see Theorems~4.5, 4.8  and Corollary~4.9).

\proclaim{Theorem G} Let  $T$ be a $(C,F)$-action of $\Bbb Z$ associated with  bounded parameters $(C_n,F_{n-1})_{n\ge 1}$ and let $
F_n\in\Cal F$ for each   $n>0$\footnote{The condition $F_n\in\Cal F$ for each $n$ means that $T$ is of rank one.}.
Then $T$ is rigid if and only if for each $N>0$, there are integers $n,m$ such that $m>n+N>n>N$  and  the set
$C_n+\cdots+C_m$ is an arithmetic sequence.
\endproclaim

\proclaim{Theorem H} Let $T=(T_g)_{g\in\Bbb Z}$ be a $(C,F)$-action of $\Bbb Z$ associated with   a sequence $(C_n,F_{n-1})_{n\ge 1}$ and let $
F_n\in\Cal F$ for each   $n>0$.
\roster 
\item"(i)"
If $T_d$ is  ergodic
then for each divisor $p$ of $d$, there are infinitely many $n>0$ such that some $c\in C_n$ is not divisible by $p$.
\item"(ii)"
If the sequence $(\# C_n)_{n=1}^\infty$ is bounded  
and for each divisor $p$ of a positive integer $d$,  there are infinitely many  $n$ such that $p$ does not divide some $c\in C_n$ then $T_d$ is ergodic.
\item"(iii)"
If the sequence $(\# C_n)_{n=1}^\infty$ is bounded 
then $T$ is totally ergodic if and only if for each $d>1$, there are infinitely many $n>0$ such that
some element $c$ of $ C_n$ is not divisible by $d$.
\endroster
\endproclaim

Comparing the   above two theorems with the main results of  \cite{GaHi3},
 we note that Theorem G extends \cite{GaHi3, Theorem~1.2} from the class of {\it adapted} rank-one transformations  with bounded parameters to the class of {\it all} rank-one transformations with bounded parameters.
In Theorem~H(iii) we consider both finite measure preserving and infinite measure preserving systems, in particular, extending \cite{GaHi3, Theorem~1.3}, where the finite measure preserving case is considered. 
Our claims (i) and (ii) seem to be considered for the first time.
We also note that claim (i) is proved for arbitrary rank-one $\Bbb Z$-actions while 
\cite{GiHi3} studied  transformations with {\it bounded cuts}, i.e. the sequence 
$(\# C_n)_{n=1}^\infty$ was assumed bounded there.


In \cite{Ry}, Ryzhikov stated a theorem that the totally ergodic and non-rigid rank-one transformations with bounded parameters have the property of {\it minimal self-joinings (MSJ)}.
We refer to \cite{dJRu}  and \cite{Ru}  for the definition of MSJ.
The theorem generalizes  the well known result from \cite{dJRaSw} that  the Chacon transformation with 3 cuts has MSJ.
Ryzhikov provided a sketch of a proof that is based on the limit properties of joinings and the weak limits of powers.
In \cite{GaHi4}, Gao and Hill present a different proof of this result via tools of symbolic dynamics.
However they do only a particular case of Ryzhikov's theorem because they consider only  adapted rank-one transformations (see Theorem~L and the remark just below it). 

In \S 5, we provide a detailed proof of the {\it full version} of Ryzhikov's theorem in the framework of $(C,F)$-construction (see~Theorem~5.3 and Corollary~5.4):

\proclaim{Theorem I} Let $T$ be a  rank-one $\Bbb Z$-action with bounded parameters.
Suppose that $T$ is not rigid and that $T$ is totally ergodic. 
Then $T$ has  MSJ.
Hence $T_n$ and $T_m$ are disjoint\footnote{A stronger result that they are spectrally disjoint was proved in \cite{EALedR}.} for all $n\ne m\in\Bbb N$.
\endproclaim

Our proof is based on standard analysis   of generic points for the self-joinings of the rank-one maps under question (neither  weak limits of powers nor limit properties of joinings appear in our proof).
As a byproduct, we find a criterion for isomorphism and disjointness (in  Furstenberg sense \cite{Fu}) for  commensurate non-rigid rank-one transformations with bounded parameters (see Corollary~5.5).

\proclaim{Theorem J}  Let  $T$  and $T'$ be two $(C,F)$-action of $\Bbb Z$ associated with  bounded parameters $(C_n,F_{n-1})_{n\ge 1}$ and $(C_n',F_{n-1}')_{n\ge 1}$ and let $
F_n=F_n'\in\Cal F$ eventually.
Let $T$ not be rigid.
\roster
\item"$(i)$"
Then
$T$ and $T'$ are isomorphic if and only if $C_n=C_n'$ eventually.
\item"$(ii)$" If  $C_n\ne C_n'$ for infinitely many $n$
and
for each $n>0$, either $T_n$ or $T'_n$ is ergodic then
$T$ and $T'$ are disjoint.
\endroster
\endproclaim

As an application we obtain  a satisfactory solution of the {\it measure theoretical inverse} problem
within the class of  non-rigid rank-one transformations with bounded parameters (see Corollary~5.7 and cf. it with Theorem~D above).

\proclaim{Theorem K} Let $T$ be a $(C,F)$-action of $\Bbb Z$ associated with bounded parameters $(C_n,F_{n-1})_{n\ge 1}$
and let $
F_n\in\Cal F$ eventually.
If $T$ is not rigid then
 $T_1$ and $T_{-1}$ are 
 isomorphic if and only if  $C_n=\{\max C_n-c\mid c\in C_n\}$ eventually.
 If, moreover, $T$ is totally ergodic and $T_1$ is not isomorphic to $T_{-1}$ then $T_1$ and $T_{-1}$ are disjoint.
 \endproclaim

We note that  Theorems J and K  extend  the main results from \cite{GaHi4} and \cite{Hi2}, where only adapted rank-one transformations  were under consideration.
In this connection (see also our remark just below Theorem~H) a natural question arises: is the class of adapted rank-one transformations with bounded parameters (considered up to measure theoretical isomorphism)  is really less than the class of all rank-one transformations with bounded parameters?
If we drop the boundedness restriction then  the two classes coincide (see~Lemma~1.10).
  The affirmative answer follows from the next theorem (see Theorem~5.9).

\proclaim{Theorem L} Let  $T$ be  an adapted  rank-one action of $\Bbb Z$ with bounded parameters\footnote{In fact, it suffices to claim that only  the sequence of spacers is bounded.}. 
Then there is  a sequence $n_m\to +\infty$ and  a polynomial $P(Z)=\nu_0+\nu_1 Z+\cdots+\nu_K Z^K$ with non-negative real coefficients  $\nu_i$ such that  $\sum_{0\le i\le K}\nu_i=1$ and
$
T_{-{n_m}}\to P(T_1)
$
as $m\to\infty$ in the weak operator topology\footnote{We identify $T_n$  here with the unitary Koopman  operators in $L^2(X,\mu)$ generated by it, $n\in\Bbb Z$.}.
It follows that $T$ is not lightly mixing.
 \endproclaim

 Since the Chacon map with 2 cuts is lightly mixing \cite{FrKi}, it follows that it is not isomorphic to any adapted rank-one transformation with bounded parameters.

 The following quadchotomy theorem for rank-one transformations with bounded parameters  refines (with a different proof) the trichotomy  theorem \cite{EALeRu, Theorem~3} (see Theorem~5.10).

\proclaim{Theorem M} Let  $T$ be a rank-one $\Bbb Z$-action with bounded paprameters.
Let $K$ denote the upper bound for  the number of spacers put on a subcolum over all subcolums and all steps of the inductive cutting-and-stacking construction.
Then one of the four possibilities takes place:
\roster
\item"\rom(i)" $T$ has MSJ (in particular, $T$ is weakly mixing and $C(T)=\{T_n\mid n\in\Bbb Z\})$,
\item"\rom(ii)"
$T$ is non-ridid, the group $\Lambda_T\subset\Bbb T$ of eigenvalues of $T$ is nontrivial but finite and
 the order of each $\lambda\in\Lambda_T$ does not exceed  $K$.
 For each  ergodic 2-fold self-joining of $T$ which is neither a graph-joining nor $\mu\times\mu$, there is $\lambda\in\Lambda_T\setminus\{1\}$
 and $n>0$ such that $\lambda^n=1$ and
 $\frac 1n\sum_{i=0}^{n-1}\rho\circ(I\times T_i)=\mu\times\mu$.

\item"\rom(iii)" $T$ is rigid,
the group $\Lambda_T\subset\Bbb T$ of eigenvalues of $T$ is finite and
 the order of each $\lambda\in\Lambda_T$ does not exceed  $K$.
\item"\rom(iv)" $T$ is an odometer of bounded type\footnote{For the definition of odometers of bounded type see several lines above Proposition~4.6.}.
\endroster
\endproclaim
  
  After this paper had been submitted   I learned about a recent work \cite{AdFePe} by Adams, Ferenczi and Petersen.
  The main result of their work is that {\it every} probability preserving rank-one map defined by the cutting-and-stacking construction process admits a {\it  symbolic presentation} as a uniquely ergodic binary subshift.
  Thus the {\it constructive symbolic definition} of rank one given in the famous survey \cite{Fe2} (and in a sense conflicting with the usual  definition of odometers) is equivalent to the standard cutting-and-stacking one.
  Then I wrote  Appendix, where I use  the $(C,F)$-machinery\footnote{Adams, Ferenczi and Petersen use the Bratteli-Vershik models of rank-one maps to prove the main result of \cite{AdFePe}.} to give a short alternative proof of this result.
  It follows from the next theorem (see Theorem~A.3):   
  
  \proclaim{Theorem N} Each finite measure preserving rank-one transformation is isomorphic to a rank-one transformation  $(X,\goth B,\mu, T)$   that is essentially $0$-expansive, i.e. the smallest $T$-invariant sub-$\sigma$-algebra containing the initial level of the cutting-and-stacking inductive construction of $T$ is $\goth B\pmod \mu$.
  \endproclaim

Theorem N follows from the  fact that each finite measure preserving rank-one transformation is isomorphic to a rank-one transformation defined by the cutting-and-stacking inductive construction 
 in such a way that on each step of this construction, the number of spacers put over the last subtower is strictly greater than the number of spacers put over each other subtower.
This fact is, in turn,  a slight refinement of a particular case of \cite{Da6, Theorem~2.8}.

\subsubhead Acknowledgement
\endsubsubhead
The author is grateful to M. Lema{\'n}czyk for useful discussions and valuable suggestions related to  this work.
I also thank the referee  for his (her) remarks and finding a gap in the original version of Lemma~2.1.

 \head 1. $(C,F)$-construction of rank-one actions
 \endhead
 
 \subhead 1.1. Actions of rank one
 \endsubhead
 Let $G$ denote a discrete countable infinite group.
 Fix an infinite sequence $\Cal F=( \Cal F_n)_{n=0}^\infty$ of finite subsets in $G$.
Let $T=(T_g)_{g\in G}$ be a measure preserving  action of   $G$ on a standard $\sigma$-finite measure space $(X,\goth B,\mu)$.

 \definition{Definition 1.1}
If there exist a sequence  $(B_n)_{n\ge 0}$ of subsets of finite measure  in $X$ and an increasing  sequence  $(l_n)_{n\ge0}$ of non-negative integers such that  
\roster
\item"(i)" 
for each $n\ge 0$, the subsets $T_gB_n$, $g\in \Cal F_{l_n}$, are pairwise disjoint and 
\item"(ii)" 
for each subset $B\in\goth B$ with $\mu(B)<\infty$,
$$
\lim_{n\to\infty}\inf_{F\subset \Cal F_{l_n}}\mu\bigg(B\triangle\bigsqcup_{g\in F}T_gB_n\bigg)=0
$$
\endroster
then we say that $T$ is  {\it of  rank one\footnote{Sometimes $T$ is  called of {\it funny} rank one.} along $\Cal F$}.
\enddefinition

If $G=\Bbb Z$ and  $\Cal F=(\{0,1,\dots, n\})_{n> 0}$ then  we obtain the standard definition of {\it rank-one transformations} (or $\Bbb Z$-actions). 
If $G=\Bbb Z$ but $\Cal F$ is the set of all finite subsets of $\Bbb Z_+$ then we obtain the definition of  transformations of {\it funny rank one} \cite{Fe1}.
If $G=\Bbb Z^d$  for $d>1$ and  $\Cal F=(\{0,1,\dots, n\}^{ d})_{n>0}$ then  we obtain the definition of $\Bbb Z^d$-actions  of {\it rank one along cubes}. 
 
 \remark{Remark 1.2} We note that  if  $T=(T_g)_{g\in G}$  is of rank one along $\Cal F$ then we can assume without loss of generality that the following property holds in addition to~(i) and (ii):
 \roster
 \item"(iii)" for each $n\ge0$, there is a subset $C_n\subset F_{l_{n+1}}$ such that $B_n=\bigsqcup_{g\in C_n}T_gB_{n+1}$.
 \endroster
 \endremark
 We leave the proof of this standard claim as an exercise to the reader.
 
 \subhead 1.2. $(C,F)$-spaces,  canonical measures and tail equivalence relations
 \endsubhead
For a detailed exposition of  the $(C,F)$-concepts we refer to \cite{Da1}, \cite{Da5} and \cite{Da6}.
 Let $(F_n)_{n\ge 0}$  and $(C_n)_{n\ge 1}$ be two sequences of finite subsets in $G$ such that for each $n>0$,
\roster
\item"(I)"
  $1\in F_0\cap C_n$,  $\# C_n>1$,
\item"(II)"
$F_nC_{n+1}\subset F_{n+1}$,
\item"(III)"
$F_nc\cap F_nc'=\emptyset$ if $c,c'\in C_{n+1}$ and $c\ne c'$.
\endroster
We  let $X_n:=F_n\times C_{n+1}\times C_{n+2}\times\cdots$ and endow this set with the infinite product topology.
Then $X_n$ is a compact Cantor (i.e. totally disconnected  perfect metric) space.
The  mapping
$$
X_n\ni (f_n,c_{n+1},c_{n+2},\dots)\mapsto(f_nc_{n+1}, c_{n+2},\dots)\in X_{n+1}
$$
is a topological embedding
of $X_n$ into $X_{n+1}$.
Therefore an inductive limit $X$ of the sequence $(X_n)_{n\ge 0}$ furnished with these embeddings  is a well defined locally compact Cantor  space.
We call it the {\it $(C,F)$-space associated with the sequence} $(C_n,F_{n-1})_{n\ge 1}$.
It is easy to see that the $(C,F)$-space is compact if and only if there is $N>0$ with $F_{n+1}=F_nC_{n+1}$ for all $n>N$.
Given $n\ge 0$ and a subset $A\subset F_n$, we let 
$$
[A]_n:=\{x=(f_n,c_{n+1},\dots)\in X_n\mid f_n\in A\}
$$ 
and call this set an  {\it $n$-cylinder} in $X$.
It is open and compact in $X$.
The collection of all cylinders coincides with the family of  all compact open subsets  in $X$.
It is easy to see that
$$
\aligned
[A]_n\cap[B]_n &=[A\cap B]_n, \quad
[A]_n\cup[B]_n =[A\cup B]_n\quad \text{and}\\
 [A]_n &=[AC_{n+1}]_{n+1}
 \endaligned
 \tag1-1
 $$
  for all $A,B\subset F_n$ and  $n\ge 0$.
For brevity, we will write $[f]_n$ for $[\{f\}]_n$, $f\in F_n$.

Let $\Cal R$ denote the {\it tail equivalence relation} on $X$.
This means that for each $n\ge 0$, the restriction of $\Cal R$ to $X_n$ is the tail equivalence relation on $X_n$, i.e. two points $(f_n,c_{n+1},\dots)$  and $(f_n',c_{n+1}',\dots)$ from $X_n$ are equivalent if  there is $N>0$ such that $c_m=c_m'$ for all $m>N$.
We note that $\Cal R$ is {\it minimal}, i.e. the $\Cal R$-class of every point  is dense in $X$ and {\it uniquely ergodic}, i.e.
there exists a unique up to scaling non-zero $\sigma$-finite $\Cal R$-invariant\footnote{$\mu$ is called $\Cal R$-invariant if $\mu$ is invariant under each Borel transformation whose graph is contained in $\Cal R$.}   Radon\footnote{i.e. it is finite on every compact subset. Every Radon measure on $X$ is $\sigma$-finite.}  measure $\mu$ on $X$.
Moreover,  $\mu$ is strictly positive on every non-empty open subset.
We note that the {\it $\Cal R$-invariance} of $\mu$ is equivalent to the following property:
$$
\mu([f]_n)=\mu([f']_n)\quad\text{ for all }f,f'\in F_n, n\ge 0.
$$
Using this property and \thetag{1-1} we can compute  that 
$$
\mu([A]_n)=\frac{\mu([1]_0)\# A}{\# C_1\cdots\# C_n} \quad\text{for each subset }A\subset F_n, n>0.
$$
In what follows we normalize $\mu$ by the condition $\mu([1]_0)=1$.
We see that the restriction of $\mu$ to $X_0$ is the product of the ``counting'' measure on $F_0$ with the infinite product of equidistributions on $C_n$, $n\in\Bbb N$.
We call $\mu$ the {\it canonical measure associated with $(C_n, F_{n-1})_{n\ge 1}$}.
It is finite if and only if\footnote{In view of (I)--(III), the sequence $(\frac{\# F_n}{\# C_1\cdots\#C_n})_{n=1}^\infty$ is non-decreasing and bounded by $\#F_0$ from below.} 
$$
\lim_{n\to\infty}\frac{\# F_n}{\# C_1\cdots\#C_n}<\infty.
\tag1-2
$$
It is easy to see that  $\mu$   is {\it $\Cal R$-ergodic}, i.e. each Borel $\Cal R$-saturated subset of $X$ is either $\mu$-null or $\mu$-conull.

\subhead 1.3. $(C,F)$-actions
\endsubhead
We now define an action of $G$ on $X$ (or, more rigorously,  on a subset of $X$).
Given $g\in G$, let
$$
X_n^g:=\{(f_n,c_{n+1},c_{n+2}\dots)\in X_n\mid gf_n\in F_n\}.
$$
Then $X_n^g$ is a compact open subset of $X_n$ and $X_n^g\subset X_{n+1}^g$.
Hence the union $X^g:=\bigcup_{n\ge 0}X_n^g$ is  a well defined  open subset of $X$.
Let $X^G:=\bigcap_{g\in G}X^g$.
Then $X^G$ is a $G_\delta$-subset of $X$.
Hence $X^G$ is Polish in the induced topology.
Given $x\in X^G$ and $g\in G$, there is $n>0$ such that $x=(f_n,c_{n+1},\dots)\in X_n$ and $g f_n\in F_n$.
We now let 
$$
T_gx:=(gf_n, c_{n+1},\dots)\in X_n\subset X.
$$
It is standard to verify that 
\roster
\item"(i)"$T_gx\in X^G$,
\item"(ii)" the map $T_g:X^G\ni x\mapsto T_gx\in X^G$ is a homeomorphism of $X^G$ and 
\item"(iii)"$T_gT_{g'}=T_{gg'}$ for all $g,g'\in G$.
\endroster
Thus $T:=(T_g)_{g\in G}$ is a continuous action of $G$ on $X^G$.

\definition{Definition 1.3}
We call $T$ {\it the $(C,F)$-action of $G$ associated with the sequence $(C_n,F_{n-1})_{n\ge 0}$}.
\enddefinition

Each $(C,F)$-action is free (except for the trivial case where  $X^G=\emptyset$).
It is obvious that $X^G$ is $\Cal R$-invariant and the $T$-orbit equivalence relation is the restriction of $\Cal R$ to $X^G$.
It follows that  $T$ preserves  $\mu$.
Since $X^G$  is $\Cal R$-saturated  and  $\mu$ is $\Cal R$-ergodic, we have
either $\mu(X^G)=0$ or  $\mu(X\setminus X^G)=0$.
Each of the two cases can occur.

\proclaim{Lemma 1.4 \cite{Da6, Theorem 2.5}}
\roster
\item"$(i)$"
 $X^G=X$ if and only if for each $g\in G$ and $n>0$, there is $m>n$ such that
$$
gF_nC_{n+1}C_{n+2}\cdots C_{m}\subset F_m.\tag1-3
$$
\item"$(ii)$"
 $\mu(X\setminus X^G)=0$
 if and only if for each $g\in G$ and $n>0$,
$$
\lim_{m\to\infty}\frac{\#((gF_nC_{n+1}C_{n+2}\cdots C_{m})\cap F_m)}{\# F_n\#C_{n+1}\cdots\# C_m}=1.
\tag1-4
$$
\item"$(iii)$"
If $\mu(X)<\infty$ then $\mu(X\setminus X^G)=0$ if and only if 
$(F_n)_{n\ge 0}$ is a F{\o}lner sequence in $G$ and hence $G$ is amenable.
\endroster
\endproclaim

We also note that $T$ is of rank one along $(F_n)_{n\ge 0}$.
The converse assertion holds also.

\proclaim{Lemma 1.5 (\cite{Da6, Theorems~2.6, 2.8})}
If $T$ is a $\sigma$-finite measure preserving $G$-action of  rank one along $\Cal F$ then $T$ is (measure theoretically) isomorphic to a $(C,F)$-action of $G$ on the $(C,F)$-space  equipped with the canonical measure.
Moreover,  without loss of  generality we may assume that the corresponding $(C,F)$-sequence $(C_n,F_{n-1})_{n\ge 1}$ satisfy \thetag{1-3} and $(F_n)_{n\ge 0}$ is a subsequence of $\Cal F$.
\endproclaim

\remark{Remark \rom{1.6}}
We note that if $X^G=X$ but $X$ is not compact then  $T$ extends to the one-point compactification $X^*=X\sqcup\{\infty\}$  of $X$ by setting $T_g\infty=\infty$ for all $g\in G$.
We thus obtain a continuous action of $G$ on the compact Cantor space $X^*$.
This action is {\it almost minimal}, i.e. there is one fixed point and the orbit of any other point is dense.
This concept was introduced in \cite{Da2} in the case $G=\Bbb Z$. 
\comment
For the (topological) orbit classification of the almost minimal $\Bbb Z$-systems see \cite{Da2} and \cite{Ma}.
Some natural examples of such systems (subshifts arising from non-primitive substitutions) are given in \cite{Yu1}.
\endcomment
\endremark

\subhead 1.4. Telescoping \endsubhead
We now introduce an important concept of telescoping for  the $(C,F)$-sequences.
Given a non-decreasing sequence
 $(k_n)_{n\ge0}$ of non-negative integers, we let
$\widetilde F_n:=F_{k_n}$ and $\widetilde C_{n}:=C_{k_{n-1}+1}C_{k_{n-1}+2}\cdots C_{k_{n}}$.
We call the sequence $(\widetilde C_n,\widetilde F_{n-1})_{n>0}$  {\it the $(k_n)_{n\ge0}$-telescoping of $( C_n,F_{n-1})_{n>0}$.} 
It is easy to see that if  
$( C_n, F_{n-1})_{n>0}$ satisfies~(I)--(III) and~\thetag{1-3} (or \thetag{1-4}) then    $(\widetilde C_n,\widetilde F_{n-1})_{n>0}$   satisfies the same conditions.
Denote by  $\widetilde T$ 
 the $(C,F)$-action associated with $(\widetilde C_n,\widetilde F_{n-1})_{n>0}$.
  Then $\widetilde T$ is canonically isomorphic to $T$.
  Indeed, let $X$ and $\widetilde X$ denote the corresponding $(C,F)$-spaces
and
  $$
  X=\bigcup_{n\ge 0}X_{n}=\bigcup_{n\ge 0}X_{k_n}\quad\text{ and }\quad\widetilde X=\bigcup_{n\ge 0}\widetilde X_{n},
  $$ where $X_{n}=F_{n}\times C_{n+1}\times\cdots$
  and $\widetilde X_{n}=\widetilde F_{n}\times \widetilde C_{n+1}\times\cdots$.
Then the mappings
$$
 X_{k_n}\ni (f_{k_n}, c_{k_n+1},\dots)\mapsto (f_{k_n}, c_{k_n+1}\cdots c_{k_{n+1}},c_{k_{n+1}+1}\cdots c_{k_{n+2}},\dots)\in\widetilde X_n,
$$ 
$n\ge 0$, define a homeomorphism of $X$
onto $\widetilde X$.
It is easy to see that this homeomorphism  intertwines $T$ with $\widetilde T$ and the  tail equivalence relation on $X$ with the  tail equivalence relation on $\widetilde X$.

If $\sup_{n\ge 0}(k_{n+1}-k_n)<\infty$ we call the $(k_n)_{n\ge0}$-telescoping {\it bounded}.

\subhead 1.5. An application to topological group actions
\endsubhead
We now show how to use  the $(C,F)$-construction (plus Proposition~1.8 below) to answer  some  questions in the theory of topological group actions stated in \S 0.
\comment
\roster
\item"---" does $G$ admit a free minimal amenable (in the topological sense, see Definition~1.6 below) action on
a locally compact non-compact Cantor space $X$  \cite{KeMoR\o, Question~7.1}\footnote{A partial answer to this question was obtained recently in \cite{MaR\o} in the case where $G$ is {\it exact}, i.e. $G$ admits an amenable action on a compact set.}?
\item"---"  does $G$ admit a free minimal amenable action on
 $X$ which leaves invariant a non-zero Radon measure on $X$ \cite{KeMoR\o, Question~7.2}?
\endroster
\endcomment

\definition{Definition 1.7  \cite{An, Proposition 2.2}} A continuous action $T=(T_g)_{g\in G}$ of $G$ on a locally compact
second countable space $X$ is called {\it amenable} if
there exists a sequence $(g_i)_{i=1}^\infty$ of nonnegative continuous functions on $X \times G$  such
that  
\roster
\item"\rom{(a)}"
for every  $i\in\Bbb N$ and $x \in X$, $\int _Gg_i (x,t)\,dt = 1$;
\item"\rom{(b)}"
$\lim_{i\to\infty}\int_G |g_i(T_sx,st) - g_i(x,t)|dt = 0$ uniformly on compact subsets of $X \times G$.
\endroster
\enddefinition

\proclaim{Proposition 1.8} Let $T$ be a $(C,F)$-action of $G$ associated with a sequence $(C_n,F_{n-1})_{n\ge 1}$ satisfying (I)--(III) and \thetag{1-3}.
Then $T$ is amenable.
\endproclaim
\demo{Proof}
For each $i>0$, we 
define $g_i:X\times G\to\Bbb R_+$ by setting
$$
g_i(x,h):=\cases 
\frac1{\# F_i}1_{F_i}(h^{-1}f_i) &\text{if }x=(f_i,c_{i+1},\dots)\in X_i:=F_i\times C_{i+1}\times\cdots\\
\frac1{\# F_i}1_{F_i}(h^{-1}) &\text{if }x\not\in X_i.
\endcases
$$
Then $g_i$ is continuous and $\int_Gg_i(x,h)dh=1$ for each $x\in X$ and $i>0$.
Now fix $n>0$ and $s\in G$.
It follows from \thetag{1-3} that there is $i>n$ such that 
  $F_n\cup s F_nC_{n+1}\cdots C_i\subset F_i$ and hence $X_n\cup T_sX_n\subset X_i$.
Therefore for each $x=(f_n,c_{n+1},\dots)\in X_n$, we have 
$$
T_sx=T_s(f_i,c_{i+1},\dots)=(sf_i,c_{i+1},\dots),
$$
where $f_i=f_nc_{n+1}\cdots c_i$.
This yields
 $g_i(T_s x,sh)=g_i(x,h)
$
 and
$$
\lim_{i\to\infty}\max_{x\in X_n}\int_G|g_i(x,h)-g_i(T_s x,sh)|dh=0.
$$
Hence $T$ is amenable. 
\qed
\enddemo

Now, for each discrete countable infinite group $G$, take a sequence $(C_n, F_{n-1})_{n\ge 1}$ satisfying (I)--(III) and \thetag{1-3}.
This can be done easily via an inductive construction procedure.
Then the $(C,F)$-action $T$ associated with  $(C_n, F_{n-1})_{n\ge 1}$ 
 is a free minimal amenable action of $G$ on a locally compact Cantor space  and $T$ leaves invariant a unique (up to scaling) non-trivial Radon measure.
 This answers Questions 7.2(i) and (ii) from \cite{KeMoR\o}.

\subhead 1.6. $(C,F)$-concepts and the classical  ``cutting-and-stacking'' nomenclature in case of $\Bbb Z$-actions
\endsubhead
Suppose that $G=\Bbb Z$.
We recall the classical cutting-and-stacking construction of rank-one transformations (see, e.g. \cite{Ru}).
A  {\it tower} $A$ is an ordered finite collection of pairwise disjoint intervals (called the {\it levels} of $A$) in $\Bbb R$, each of the same Lebesgue measure. 
We think of the levels in a tower as being stacked on top of each other, so that the $(j + 1)$-st level is directly above the $j$-th level. 
Every tower $A = (I_j)_j$ is associated with a natural tower map $T_A$ sending each point in $I_j$ to the point directly above it in $I_{j+1}$. 
A rank-one cutting-and-stacking construction of a measure preserving transformation $T$ consists of a sequence of towers $(A_n)_{n\ge 0}$ such that $A_0$ is a single interval  $[0,1)$,
each tower $A_{n+1}$ is obtained from $A_n$ by cutting $A_n$ into
$r_n\ge 2$ subtowers of equal width, adding some number $\sigma_n(k)$ of new levels (called {\it spacers}) above the $k$-th subcolumn, $k=0,\dots,r_n-1$, and stacking every subtower under the subtower to its right.
 We note the spacers are intervals drawn from $\Bbb R$ that are disjoint from the levels of $A_n$ and the other spacers added to it.
  They are of the same length as the levels of the subcolumns of $A_n$.
  It is easy to see that $T_{A_{n+1}}\restriction A_n=T_{A_n}$ for each $n$.
  We now set $X:=\bigcup_{n\ge 0}\bigsqcup_{I\in A_n}I$, endow $X$ with the Lebesgue measure 
and define $T$ to be
 the pointwise limit of $T_{A_n}$ as $n\to\infty$.
 Then $T$ is a measure preserving invertible transformation of $X$.
We note that $T$ is completely defined by the sequence of  integers $(r_n)_{n\ge 1}$ and the sequence $(\sigma_n)_{n\ge 1}$ of maps $\sigma_n:\{1,\dots,r_n\}\to\Bbb Z_+$.
We denote this  fact by writing $T\sim(r_n,\sigma_n)_{n=1}^\infty$.
For example, if we let $r_n=2$, $\sigma_n(1)=0$ and $\sigma_n(2)=1$ for each $n\in\Bbb N$ then the rank-one transformation $T\sim(r_n,\sigma_n)_{n=1}^\infty$ is the {\it Chacon 2-cuts map}.
If $r_n=3$, $\sigma_n(1)=\sigma_n(3)=0$ and $\sigma_n(2)=1$ for each $n\in\Bbb N$ 
then the rank-one transformation $T\sim(r_n,\sigma_n)_{n=1}^\infty$ is the {\it Chacon 3-cuts map}.

We now show how to obtain $T$ via the $(C,F)$-construction.
For that we set $h_0:=0$, $h_{n+1}:=h_nr_n+\sum_{i=1}^{r_n}\sigma_n(i)$,
 $F_n:=\{0,1,\dots,h_n-1\}$, $s_n(0):=0$, $s_n(i):=\sum_{j\le i}\sigma_n(i)$ if  $i=1,\dots,r_n-1$ and 
$C_{n+1}:=\{ih_n+s_n(i)\mid i=0,\dots,r_n-1\}$.
In the sequel we will call $s_n$ {\it the integral} of $\sigma_n$.
Then the sequence $(C_n,F_{n-1})_{n\ge 1}$ satisfies (I)--(III).
Moreover, it satisfies \thetag{1-4} because $(F_n)_{n\ge 0}$ is a F{\o}lner sequence in $\Bbb Z$.
Hence the associated $(C,F)$-action of $\Bbb Z$ is well defined.
It is standard to see that this action is isomorphic to $(T^n)_{n\in\Bbb Z}$ by an isomorphism that identifies (uniquely, in accordance with the orders) the levels of $A_n$ with the cylinders $\{[f]_n\mid f\in F_n\}$ for each $n>0$
(the  order on $F_n$ is inherited from the standard linear order on $\Bbb Z$).

Conversely, let $\Cal F:=\{\{0,\dots,n\}\mid n\in\Bbb N\}$ and let
  $T=(T_g)_{g\in\Bbb Z}$ be a $(C,F)$-action associated with a sequence $(C_n,F_{n-1})_{n\ge 1}$ such that $F_n\in\Cal F$ for every $n$.
  Then $T_1\sim(r_n,\sigma_n)_{n=0}^\infty$, where $r_n:=\#C_{n+1}$, $h_n:=\# F_n$,
  $$
\sigma_n(i):=\cases 
 s_n(i)-s_n(i-1)&\text{if }1\le i<r_n,\\
 h_{n+1}-r_nh_n-s_n(r_n-1)&\text{if }i=r_n,
\endcases
\tag1-5
$$
and $s_n$ is the unique map $s_n:\{0,\dots,r_n-1\}\to\Bbb Z_+$ such that
$$
 C_{n+1}=\{ih_n+s_n(i)\mid i=0,\dots,r_n-1\}.
 \tag1-6
$$
We also note that for each $n\ge 0$, the pair $(C_{n+1},F_{n+1})$ uniquely defines (and, conversely, is uniquely defined by) the pair $(r_n,\sigma_n)$ via \thetag{1-5} and \thetag{1-6}.
We will denote this correspondence by $(C_{n+1},F_{n+1})\sim (r_n,\sigma_n)$.

\definition{Definition 1.9} We say that a rank-one transformation $T\sim(r_n,\sigma_n)_{n=1}^\infty$ is {\it adapted} if $\sigma(r_n)=0$ for each $n\ge 1$.
\enddefinition

The following claim is a folklore (at least, in the case of finite invariant measure). 
Unfortunately, we were unable to find a proof of this simple claim in the literature.
Therefore we
 provide below an idea of its proof.

\proclaim{Lemma 1.10} Each rank-one transformation $T$ is (measure theoretically) isomorphic to an adapted one.
\endproclaim
\demo{Idea of the proof}  Let $T\sim(r_n,\sigma_n)_{n=1}^\infty$ for some integers $r_n>1$ and maps $\sigma_n:\{1,\dots,r_n\}\to\Bbb Z_+$, $n\in\Bbb N$.
We now define a new sequence of maps $\widetilde\sigma_n:\{1,\dots,r_n\}\to\Bbb Z_+$, $n\in\Bbb N$ by setting
$$
\widetilde\sigma_1(i):=\cases\sigma_1(i)&\text{if } i<r_1,\\
0&\text{if } i=r_1
\endcases
$$
and
$$
\widetilde\sigma_n(i):=\cases\sigma_1(i)+\cdots+\sigma(r_{n-1})+\sigma_n(i)&\text{if } i<r_1,\\
0&\text{if } i=r_n
\endcases
$$
for $n>1$.
Let $\widetilde T$ stand for the rank-one transformation defined by the sequence $(r_n,\widetilde\sigma_n)_{n=0}^\infty$.
Of course, $\widetilde T$ is adapted.

As was shown above, we can assume without loss of generality that $(T^n)_{n\in\Bbb Z}$ is the $(C,F)$-action associated with a sequence $(C_n,F_{n-1})_{n\ge 1}$ such that $(C_{n+1},F_{n+1})\sim(r_n,\sigma_n)$ for each $n$ and 
$(\widetilde T^n)_{n\in\Bbb Z}$ is the $(C,F)$-action associated with a sequence $(\widetilde C_n,\widetilde F_{n-1})_{n\ge 1}$ such that $(\widetilde C_{n+1},\widetilde F_{n+1})\sim(r_n,\widetilde\sigma_n)$ for each $n$.
It is a routine to verify that $C_n=\widetilde C_n$ for each $n>0$.
Then if follows from Lemma~5.8 below that $T$ and $\widetilde T$ are isomorphic.
\qed

\comment
We then pass to a telescoping $(\widetilde C_n,\widetilde F_{n-1})_{n>0}$ of $(C_n,F_{n-1})_{n\ge 1}$ such that
$\sum_{k=0}^\infty(\#\widetilde C_n)^{-1}<\infty$.
Endow $\widetilde C_n$ $\widetilde F_{n}$ with the natural linear order inherited from $\Bbb Z$.
Put $c_n^\bullet:=\max C_n$ and $f_n^\bullet:=\max F_n$.
We now let
$F_n':=\widetilde F_n$ and $C_{n}':=(C_n\setminus\{c_n^\bullet\})\sqcup\{f_{n}^\bullet-\# F_{n-1}\}$.
Then the sequence $(C_n,F_{n-1})_{n\ge 1}$ satisfies (I)--(III) and \thetag{1-4}.
Let $T'=(T'_n)_{n\in\Bbb Z}$ stand for the associated $(C,F)$-action.
It is easy to see
that $T_1'$ is adapted.
To show that $T$ and $T'_1$  are isomorphic use the fact that for a.e. 
$x\in X_n'=F_n'\times C_{n+1}'\times C_{n+2}'\times\cdots$, if we write the expansion
$x=(f_n',c_{n+1}',c_{n+2}',\dots)$ for $x$ then there is $N=N(x)>0$ such that $c_n'\ne\max C_n'$ for all $n>N$.
Hence a measurable map
$$
x\mapsto(f_n'+c_{n+1}'+\cdots+c_N',c_{N+1}',c_{N+2}',\dots)\in \widetilde F_N\times \widetilde C_{N+1}\times \widetilde C_{N+1}\times\cdots = \widetilde X_N\subset\widetilde X
$$
is well defined.
It delivers an isomorphism (mod 0) of $T_1'$ and $\widetilde T_1$ which is, in turn, isomorphic to $T$.
\endcomment

\enddemo

\comment
We also define a {\it spacer map} $\sigma:\{1,\dots,r_n\}$ by setting
$$
\sigma(i)=\cases 
 s_n(i)-s_n(i-1)&\text{if }i<r_n,\\
 h_{n+1}-r_nh_n-s_n(r_n-1)&\text{if }i=r_n.
\endcases
$$
Then $\sigma(i)$ is the number of spacers between the $i$-th and $(i+1)$-th subtower if $i<r_n$ and $\sigma(r_n)$ in the number of spacers over the highest subtower.
We note also that $s_n(i)=\sum_{j\le i}\sigma(i)$, $i=1,\dots,r_n-1$.
\endcomment
 
 \comment

\head 1. Topological centralizer of a rank-one transformation.
\endhead

Let $T$ be a $(C,F)$-transformation on $(X,\goth B,\mu)$.
Let $\theta\in C(T)$.
Suppose that $\theta\in\text{Homeo}(X)$.
Since $\theta$ is a homeomorphism of $X$, the subset $\phi([0]_0)$ is compact and open.
Hence there is $l_1>0$ and a subset $L_1\subset F_{l_1}$ such that $\phi([0]_0)=[L_1]_{l_1}$.
Since the measure $\mu\circ\theta^{-1}$ is $T$-invariant and Radon, it follows that
$\mu\circ\theta^{-1}=c\mu$ for some $c>0$.
By \cite{Si}, $c=1$, i.e. $\theta$ preserves $\mu$.
It follows that $\# L_1=\# C_{1,l_1}$.
There is $l_2>l_1$ 
and  a subset $L_2\subset F_{l_2}$ such that $\theta^{-1}([0]_{l_1})=[L_2]_{l_2}$.
Without loss of generality we may assume that  $F_{l_1}+L_2\subset F_{l_2}$.
Since $\theta$ preserves $\mu$, we obtain that $\# L_2=\# C_{l_1+1,l_2}$.
Moreover, for all $f,f'\in F_{l_1}$,
$$
[f+L_2]_{l_2}\cap [f'+L_2]_{l_2}=T_f[L_2]_{l_2}\cap T_{f'}[L_2]_{l_2}=\theta^{-1}([f]_{l_1})\cap\theta^{-1}([f']_{l_1})=\emptyset.
$$
Hence $(f+L_2)\cap (f'+L_2)=\emptyset$.
Therefore $(F_{l_1}-F_{l_1})\cap (L_2-L_2)=\{0\}$.
We now have
$$
[0]_0=\theta^{-1}([L_1]_{l_1})=\bigsqcup_{g\in L_1}T_g\theta^{-1}([0]_{l_1})=
\bigsqcup_{g\in L_1}T_g[L_2]_{l_2}=[L_1+L_2]_{l_2}.
$$
This yields $C_{1,l_1}+C_{l_1+1,l_2}=L_1+L_2$.
Let $u:=\max C_{l_1+1,l_2}$ and $v:=\max L_2$.
It follows that $C_{1,l_1}+u=L_1+v$. 
This, in turn, implies that  $C_{l_1+1,l_2}-u=L_2-v$.
In particular,  $\theta([0]_0)=T_{u-v}[C_{1,l_1}]_{l_1}=T_{u-v}[0]_0$.
Next, do the same starting with $[0]_m$ for an arbitrary $m>0$.
Then we obtain that $\theta([0]_m)=T_{u_1-v_1}[0]_m$ for some $u_1,v_1$.
This implies that $\theta[A]_m=T_{u_1-v_1}[A]_m$ for each subset $A\subset F_m$.
Taking $A=F_m$, we obtain that $u-v=u_1-v_1$. 
Continuing, we finally obtain that  $\theta B=T_{u-v}B$ for each compact open subset $B\subset X$.
Therefore $\theta=T_{u-v}$.

Thus we have proved the following theorem.

\proclaim{Theorem} Let $T$ be a $(C,F)$-action by homeomorphisms on $X$.
Then $C(T)$ is just $T$.
\endproclaim

\endcomment

\head 2. Topological classification of $(C,F)$-actions
\endhead

\subhead 2.1. Topological isomorphism for general $(C,F)$-actions
\endsubhead
In this subsection we investigate when two $(C,F)$-actions defined on locally compact Cantor spaces are topologically isomorphic.
We first prove an auxiliary lemma.

\proclaim{Lemma 2.1} Let $T=(T_g)_{g\in G}$ and $T'=(T_g')_{g\in G}$ be two minimal  
free $G$-actions on locally compact Cantor spaces $X$ and $X'$ respectively.
Let $A$ be a compact open subset in $X$ and let $A'$ be a compact open subset in $X'$.
If there is a homeomorphism $\theta:A\to A'$ such that 
given $x\in A$ and $g\in G$, it follows that 
\roster
\item"(i)" $T_gx\in A$ if and only if $T_g'\theta x\in A'$
and 
\item"(ii)"
 $\theta T_g x=T'_g\theta x$ 
 \endroster
then $T$ and $T'$ are conjugate.
The homeomorphism  from $X$ to $X'$ intertwining  $T$ with $T'$ and extending $\theta$ is unique.
\endproclaim

\demo{Proof} 
Given $x\in X$, there is $g\in G$ such that $T_gx\in A$.
This follows from the fact that $T$ is minimal.
We now set $\widetilde\theta x:=T'_{g^{-1}}\theta T_gx\in X'$.
We first check that $\widetilde\theta$ is well defined, i.e.
if $T_hx\in A$ for some $h\in G$ then $T'_{h^{-1}}\theta T_hx=T'_{g^{-1}}\theta T_gx$.
Indeed, we have $T_{hg^{-1}}T_gx\in A$ and hence $\theta T_hx=\theta T_{hg^{-1}}T_gx=T'_{hg^{-1}}\theta T_gx$, as desired.

Next, we claim that $\widetilde\theta$ is one-to-one.
Take $x,y\in X$ with $\widetilde\theta x=\widetilde\theta y$.
Then there are $g,h\in G$ such that
 $T_gx\in A$ and $T_hy\in A$ and
 $T'_{g^{-1}}\theta T_gx=T'_{h^{-1}}\theta T_hy$.
 Hence
$T'_{hg^{-1}}\theta T_gx=\theta T_hy\in A'$.
Since $\theta T_hy=T'_h\theta y$ and $\theta T_gx=T'_g\theta x$, we obtain that 
 $\theta y=\theta x$ and hence
 $y=x$.

We now show that $\widetilde\theta$ is onto.
Take $y\in X'$.
Since $T'$ is minimal, there is $g\in G$ such that $T_g'y\in A'$.
We let $x:=T_{g^{-1}}\theta^{-1}T_g'y$.
Then $T_gx\in A$ and hence $\widetilde\theta x:=T_{g^{-1}}'\theta T_gx=y$.

It follows easily from the definition of $\widetilde\theta$ that $\widetilde\theta$ is continuous at every point of $X$.
Since $X$ is sigma-compact and $\theta$ is a bijection of $X$ onto $X'$, the mapping $\theta^{-1}$ is  continuous everywhere on $X'$.

It is straightforward to verify that $\widetilde \theta T_g=T_g'\widetilde \theta$ for each $g\in G$.

The final claim of the lemma is obvious.
\qed
\enddemo

\proclaim{Corollary 2.2} Let $T$ and $T'$ be two  $(C,F)$-actions of $G$ associated with  sequences $(C_n,F_{n-1})_{n\ge 1}$ and  $(C_n',F_{n-1}')_{n\ge 1}$ respectively
and let the sequences satisfy~(I)--(III) and \thetag{1-3}.
If $C_n=C_n'$ eventually
then $T$ and $T'$ are topologically isomorphic.
\endproclaim

\demo{Proof} Let $X$ and $X'$ stand for the $(C,F)$-spaces of $T$ and $T'$ respectively.
They are locally compact Cantor spaces.
Let $N>0$ be such that $C_n=C_n'$ for all $n>N$.
Then we set $A:=[1]_N\subset X$, $A':=[1]_N\subset X'$, $\theta x:=x$ for each $x\in A$
and apply Lemma~2.1.
\qed
\enddemo

\comment
Let $T$ and $T'$ be two $(C,F)$-actions associated with sequences $(F_n,C_{n+1})_{n\ge 0}$ and
$(F_n',C_{n+1}')_{n\ge 0}$ respectively. 
Let $\phi:X\to X'$ be an open  continuous  onto map such that $\phi T_g=T'_g\phi$.
Then we say that $T'$ is a topological factor of  $T$.
\endcomment

We now state and prove one of the main results of this section.

\proclaim{Theorem 2.3}
Let $T=(T_g)_{g\in G}$ and $T'=(T'_g)_{g\in G}$ be two $(C,F)$-actions of $G$ associated with some sequences $(C_n,F_{n-1})_{n\ge 1}$ and $(C_n',F_{n-1}')_{n\ge 1}$ respectively and the two sequences satisfy (I)--(III) and \thetag{1-3}.
Then $T$ and $T'$  are topologically isomorphic if and only if 
 there is 
an
 increasing sequence of integers
 $$0=l_0<l_1'<l_1<l_2'<l_2<\cdots$$ 
 and subsets
$A_n \subset  F_{l_n'}'$,  $B_n \subset  F_{l_n}$, 
such that
$$
\gathered
A_{n}B_{n}=C_{l_{n-1}+1}\cdots C_{l_{n}},\quad
B_{n} A_{n+1}=C'_{l_{n}'+1}\cdots C'_{l_{n+1}'},
\\
F'_{l_n'}B_n\subset F_{l_n},\quad F_{l_n}A_{n+1}\subset F'_{l_{n+1}},
\\
(F'_{l_n'})^{-1}F'_{l_n'}\cap B_nB_n^{-1}=F_{l_n}^{-1}F_{l_n}\cap A_{n+1}A_{n+1}^{-1}=\{1\}
\endgathered
\tag2-1
$$
 for each $n>0$.
\endproclaim
\demo{Proof}
$(\Rightarrow)$
Let $\phi:X\to X'$ be a  homeomorphism such that $\phi T_g=T_g' \phi$ for each $g\in G$.
 Then $\phi([1]_0)$ is a clopen subset of $X'$.
 Hence
 there are $l_1'\ge 0$ and a subset $A_{1}\subset F'_{l_1'}$ such that  $\phi([1]_0)=[A_1]_{l_1'}$.
It follows that
$$
\{\phi^{-1}([a]_{l_1'})\mid a\in A_1\}=\{T_a\phi^{-1}([1]_{l_1'})\mid a\in A_1\}
$$
is a finite partition of $[1]_0$ into compact open subsets.
In view of  \thetag{1-3}, there is $l_1\ge l_1'$ and a subset $B_1\subset F_{l_1}$ such that 
$$
\phi^{-1}([1]_{l_1'})=[B_1]_{l_1}, \quad F'_{l_1'}B_1\subset F_{l_1}, \quad [1]_0=[A_1B_1]_{l_1}
$$ 
and $a B_1\cap bB_1=\emptyset$ if $a,b\in F'_{l_1'}$ and $a\ne b$.
Since $[1]_0=[C_1\cdots C_{l_1}]_{l_1}$, we obtain that $A_1B_1=C_1\cdots C_{l_1}$.
In a similar way, it follows from the equality
 $\phi^{-1}([1]_{l_1'})=[B_1]_{l_1}$ that
 there is $l_2'>l_1$  and a subset $A_2\subset F'_{l_2'}$ such that 
 $$
 \phi([1]_{l_1})=[A_2]_{l_2'},\quad F_{l_1}A_2\subset F'_{l_2'}, \quad B_1A_2=C'_{l_1'+1}\cdots C'_{l_2'}
 $$
and $aA_2\cap bA_2=\emptyset$ if $a,b\in F_{l_1}$,  and $a\ne b$.
Continuing this process infinitely many times, we obtain an increasing sequence
$0=l_0<l_1'<l_1<l_2'<l_2<\cdots$ and subsets 
$A_{n} \subset F'_{l_{n}'}$ and $B_{n}\subset F_{l_{n}}$ for each $n>0$
such that \thetag{2-1} is satisfied.

$(\Leftarrow)$
Suppose that  \thetag{2-1} is satisfied.
Let $\widetilde T$ and $\widetilde T'$ denote the $(C,F)$-actions associated with  the $(l_{i})_{i\ge 0}$-telescoping of $(C_n,F_{n-1})_{n\ge 1}$ and the $(l'_{i})_{i\ge 1}$-telescoping of 
  $(C_n',F_{n-1}')_{n\ge 1}$ respectively.
  Let $\widetilde X$ and $\widetilde X'$ denote the $(C,F)$-spaces of these actions.
  Of course, $\widetilde T$ is isomorphic to $T$ and 
  $\widetilde T'$ is isomorphic to $T'$.
  Take $x\in [1]_0\subset\widetilde X$.
In view of \thetag{2-1}, we have a unique  expansion  $x=(1,a_1b_1,a_2b_2,\dots)\in[1]_0$ of 
$x$ such that $a_i\in A_i$ and $b_i\in B_i$ for all $i\ge 1$.
We now set
$$
\phi(x):=(a_1, b_1a_2,b_2a_3,\dots)\in[A_1]_1\subset \widetilde X'.
$$
It is standard to verify that  $\phi$ is a homeomorphism of  the cylinder $[1]_0\subset \widetilde X$ onto the cylinder $[A_1]_1\subset\widetilde X'$.
 Then 
 $$
 \align
 \{g\in G\mid  \widetilde T_gx\in [1]_0\}&=\bigcup_{n=1}^\infty  (A_1 B_1)\cdots (A_n B_n)(a_nb_n)^{-1}\cdots(a_1b_1)^{-1}\\
&=\bigcup_{n=1}^\infty A_1 (B_1A_2)\cdots (B_n A_{n+1})(b_na_{n+1})^{-1}\cdots(b_1a_2)^{-1}a_1^{-1}\\
& =\{g\in G\mid  \widetilde T_g'\phi(x)\in [A_1]_1\}.
 \endalign
 $$
 \comment

  such that $\widetilde T_gx\in [1]_0$.
 Since $x$ and $\widetilde T_gx$ have the same tails,
there is $n>0$ and $\widehat a_i\in A_i$, $\widehat b_i\in B_i$, $i=1,\dots,n$
 such that 
 $$
 \widetilde T_gx=(1,\widehat a_1\widehat b_1,\dots,\widehat a_n\widehat b_n,a_{n+1} b_{n+1},a_{n+2} b_{n+2}\dots )
 $$
and hence  $g=(\widehat a_1\widehat b_1)\cdots(\widehat a_n\widehat b_n)(a_nb_n)^{-1}\cdots(a_1b_1)^{-1}$.
 On the other hand,
 $$
 \phi(\widetilde T_gx)=(\widehat a_1, \widehat b_1\widehat a_2,\dots,\widehat b_{n-1}\widehat a_n, \widehat b_{n}a_{n+1},  b_{n+1}a_{n+2},\dots)
 $$

 \endcomment
The equality $ \phi(\widetilde T_gx)=\widetilde T_g\phi(x)$ whenever $\widetilde T_gx\in[1]_0$ is now verified in a straightforward way.
 It remains to apply Lemma~2.1.
\qed
\enddemo

\subhead 2.2. Topological isomorphism for $(C,F)$-actions of linearly ordered Abe\-lian groups
\endsubhead
We first recall the  definition of  linearly ordered countable  Abelian groups.

\definition{Definition 2.4}
If $A$ is a countable Abelian group and  $A_+$ is a subset of $A$ such that
$A_++A_+\subset A_+$,
$A_+-A_+=A$ and
$A_+\cap (-A_+)=\{0\}$
then  the pair  $(A,A_+)$ is called an {\it ordered group}.
If, in addition, $A=A_+\cup(-A_+)$, we say that $(A,A_+)$ is a {\it linearly ordered group}.
\enddefinition

For instance, $(\Bbb Z,\Bbb Z_+)$ is a linearly ordered group.
More generally, $\Bbb Z^d$ endowed with the lexicographical order is a linearly ordered group  for each $d>0$.
It was shown in \cite{Le}  that an Abelian group admits a linear order if and only if it is torsion free.

In case $G$ is linearly ordered we can strengthen Theorem~2.3.

\proclaim{Theorem 2.5} Let $(G,G_+)$ be a linearly  ordered discrete countable Abelian group.
Let $T=(T_g)_{g\in G}$ and $T'=(T'_g)_{g\in G}$ be two $(C,F)$-actions of $G$ associated with some sequences $(C_n,F_{n-1})_{n\ge 1}$ and $(C_n',F_{n-1}')_{n\ge 1}$ respectively and the two sequences satisfy (I)--(III) and \thetag{1-3}. 
Suppose that $C_n\cup C_n'\subset G_+$ for all $n$.
Then $T$ and $T'$ are topologically isomorphic if and only if 
 there are  an increasing sequence
 of integers
 $0=l_0<l_1'<l_1<l_2'<l_2<\cdots$ 
 and subsets
$A_n \subset  F_{l_n'}'\cap G_+$ and   $B_n \subset  F_{l_n}\cap G_+$
such that $0\in A_n\cap B_n$ and 
$$
\gathered
A_{n}+B_{n}=C_{l_{n-1}+1}+\cdots +C_{l_{n}},\quad B_{n}+ A_{n+1}=C'_{l_{n}'+1}+\cdots +C'_{l_{n+1}'},
\\
F'_{l_n'}+B_n\subset F_{l_n},\quad F_{l_n}+A_{n+1}\subset F'_{l_{n+1}},
\\
(F'_{l_n'}-F'_{l_n'})\cap (B_n-B_n)=(F_{l_n}-F_{l_n})\cap (A_{n+1}-A_{n+1})=\{0\}
\endgathered
$$
 for each $n>0$.
\endproclaim
\demo{Proof}
We start with a simple observation.
Let $A$ be a finite subset in $G$.
Since $G$ is linearly ordered, there exists $\min A\in G$.
Since $0\in C_n\subset G_+$, we have that 
$$
\min A=\min (A+C_n)=\min (A+C_n+C_{n+1})=\cdots
$$
for each $n>0$.
It follows from this and \thetag{1-3} that if $A\subset F_n$ then  for each sufficiently large $m>n$, we have that $ A'-\min A'\subset F_m$, where
$A':=A+C_{n+1}+\cdots+C_m$.

In view of that, we can modify in an obvious way the proof of Theorem~2.3 (increasing $l_n$ and $l_n'$ if necessary) so that 
$$
\widetilde A_n:=A_n-\min A_n\subset F'_{l_n'},\quad  \widetilde B_n:=B_n-\min B_n\subset F'_{l_n'}
$$ 
and the latter two lines of \thetag{2-1} hold for  $\widetilde A_n$ and $\widetilde B_n$
in place of $A_n$ and $B_n$ respectively
for each $n>0$.
Since $0\in \widetilde A_n$ and 
$$
\widetilde A_n+B_n+\min A_n=C_{l_{n-1}+1}+\cdots +C_{l_{n}}\subset G_+
$$ 
by the first line of \thetag{2-1} and the condition of the theorem, it follows that
$0\in B_n+\min A_n\subset G_+$.
This yields that $\min A_n=-\min B_n$. 
Hence $\widetilde A_n+\widetilde B_n=A_n+B_n=C_{l_{n-1}+1}+\cdots +C_{l_{n}}$.
In a similar way, $\min A_{n+1}=-\min B_n$ and 
$\widetilde B_{n}+ \widetilde A_{n+1}=C'_{l_{n}'+1}+\cdots +C'_{l_{n+1}'}$, i.e. the first line in \thetag{2-1} holds for 
$\widetilde A_n$ and $\widetilde B_n$
in place of $A_n$ and $B_n$ respectively for each $n>0$.
\qed
\enddemo

We now reformulate Theorem~2.5 in an equivalent way.

\proclaim{Theorem 2.6} Let $(G,G_+)$ be a linearly  ordered discrete countable Abelian group.
Let $T=(T_g)_{g\in G}$ and $T'=(T'_g)_{g\in G}$ be two $(C,F)$-actions of $G$ associated with some sequences $(C_n,F_{n-1})_{n\ge 1}$ and $(C_n',F_{n-1}')_{n\ge 1}$ respectively and the two sequences satisfy (I)--(III) and \thetag{1-3}. 
Suppose that $C_n\cup C_n'\subset G_+$ for all $n$.
Then $T$ and $T'$ are topologically isomorphic if and only if 
 there are  an increasing sequence
 of integers
 $0=l_0<l_1'<l_1<l_2'<l_2<\cdots$ 
 and subsets
$A_n \subset  F_{l_n'}'\cap G_+$ and   $B_n \subset  F_{l_n}\cap G_+$
such that $0\in A_n\cap B_n$ and 
$$
\gathered
\sum_{i>l_{n-1}} C_i = A_{n}+\sum_{i>l_n'}C'_{i},\quad
\sum_{i>l_{n}'} C_i' = B_{n}+\sum_{i>l_n}C_{i},
\\
F'_{l_n'}+B_n\subset F_{l_n},\quad F_{l_n}+A_{n+1}\subset F'_{l_{n+1}},
\\
(F'_{l_n'}-F'_{l_n'})\cap (B_n-B_n)=(F_{l_n}-F_{l_n})\cap (A_{n+1}-A_{n+1})=\{0\}
\endgathered
\tag2-2
$$
 for each $n>0$.
\endproclaim

\demo{Proof} The ``only if'' part follows immediately from Theorem~2.5.
We now prove the ``if'' part.
Indeed, it follows from the first line in \thetag{2-2} that
$$
\sum_{i>l_{n-1}} C_i = A_{n}+B_n+\sum_{i>l_n}C_{i}.\tag2-3
$$
Utilizing the definition of $A_n$ and the second line in \thetag{2-2}, we obtain that
$A_n+B_n\subset F_{l_n'}'+B_n\subset F_{l_n}$ and hence $(A_n+B_n)\cap F_{l_n}=A_n+B_n$.
Intersecting the two sides of \thetag{2-3} with $F_{l_n}$ we now get that
$\sum_{i=l_{n-1}+1}^{l_n} C_i = A_{n}+B_n$.
The equality
$B_{n}+ A_{n+1}=\sum_{i=l_{n}'+1}^{l_{n+1}'}C'_i$ is obtained in a similar way.
It remains to apply Theorem~2.5.
\qed
\enddemo

We are now ready to show that the conditions for topological isomorphism of $(C,F)$-actions are especially simple in the case of $\Bbb Z$-actions.

\proclaim{Theorem 2.7} Let $G=\Bbb Z$. Let $T=(T_g)_{g\in G}$ and $T'=(T'_g)_{g\in G}$ be two $(C,F)$-actions of $G$ associated with some sequences $(C_n,F_{n-1})_{n\ge 1}$ and $(C_n',F_{n-1}')_{n\ge 1}$ respectively and the two sequences satisfy (I)--(III) and \thetag{1-3}. 
Suppose that $C_n\cup C_n'\subset \Bbb Z_+$.
Then $T$ and $T'$ are topologically isomorphic if and only if 
there is  $r>0$\footnote{We can replace the condition $r>0$  with an equivalent  condition $r\ge i_0$ for each non-negative integer $i_0$ (passing to a corresponding telescoping).} and a subset $R\subset F_r'\cap\Bbb Z_+$ such that
$$
\sum_{i>0} C_i = R+\sum_{i>r}C'_{i}.\tag2-4
$$
\endproclaim

\demo{Proof} The ``only if'' part follows immediately from Theorem~2.6.
To prove the ``if'' part,
we start with some notation.
Given a finite subset $C\subset\Bbb Z_+$, we define 
a polynomial $P_C$ of a (single) variable  $t$ by setting
$P_C:=\sum_{c\in C}t^c$.
All non-zero coefficients of $P_C$ are equal to $1$.
If $D$ is another finite subset of $\Bbb Z_+$ such that $(C-C)\cap (D-D)=\{0\}$ then it is easy to verify that $P_CP_D=P_{C+D}$.
Therefore, if $D_1,D_2,\dots$ is an infinite sequence of finite subsets of $\Bbb Z_+$ such that
 $0\in\bigcap_{i>0}D_i$ and $(D_n-D_n)\cap\sum_{i=1}^{n-1}(D_i-D_i)=\{0\}$ for each $n>1$ then 
 the infinite product $P_{D_1}P_{D_2}\cdots$ is well defined as a formal power series.
 Moreover, it is easy to see that
 $$
 P_{D_1}P_{D_2}\cdots=\sum_{d\in D_1+D_2+\cdots}t^d.
 $$
 We also note that given three finite subsets $C,D,E$ of $\Bbb Z$, if $P_CP_D=P_E$ then
 $(C-C)\cap (D-D)=\{0\}$ and $E=C+D$.
 
 Suppose now that \thetag{2-4} holds.
 We are going to construct  sequences $(l_n)_{n\ge 0}$, $(l_n')_{n\ge 1}$, $(A_n)_{n\ge 1}$ and $(B_n)_{n\ge 1}$ satisfying \thetag{2-2}  to apply Theorem~2.6.
 This will be done in an inductive way.
 We first set $l_0:=0$, $l_1':=r$ and $A_1:=R$.
 Suppose that we have already constructed   integers $0=l_0<l_1'<\cdots<l_{n-1}<l_n'$ and subsets $A_1,\dots, A_n$, $B_1,\dots,B_{n-1}$     such that
 $
 \sum_{i>l_{n-1}}C_i=A_n+\sum_{i>l_n'}C_i'
 $.
 It follows  that
 $$
 P_{C_{l_{n-1}+1}}P_{C_{l_{n-1}+2}}\cdots=P_{A_n}P_{C_{l_n'+1}'}P_{C_{l_n'+2}'}\cdots.\tag2-5
 $$
 Hence, the polynomial $P_{A_n}$ divides the lefthand side of~\thetag{2-5}.
The factrorization property of polynomials with integer coefficients yields that there is $l_n>l_{n}'$
 and a polynomial $Q_n$ with integer coefficients such that
 $$
 Q_nP_{A_n}=P_{C_{l_{n-1}+1}}\cdots P_{C_{l_n}}.\tag2-6
 $$
 From this and \thetag{2-5} we deduce that
 $$
  Q_nP_{C_{l_n+1}}P_{C_{l_n+2}}\cdots=P_{C_{l_n'+1}'}P_{C_{l_n'+2}'}\cdots.\tag2-7
 $$
 Since $\deg Q_n<\deg (P_{C_{l_n+1}}\cdots P_{C_{l_n}})<\min (C_{l_n+1}\setminus\{0\})$,
 it follows from \thetag{2-7}\footnote{Decompose the two sides of \thetag{2-7} into formal power series.} that there is
 a  finite subset $B_n\subset C'_{l_n'+1}+C'_{l_n'+2}+\cdots$ such that $0\in B_n$ and $Q_n=P_{B_n}$.
 Hence $(B_n-B_n)\cap(F'_{l_n'}-F'_{l_n'})=\{0\}$.
 It follows from \thetag{2-6} that $A_n+B_n=C_{l_{n-1}+1}+\cdots +C_{l_n} $.
 Hence $B_n\subset F_{l_n'}+C_{l_n'+1}+\cdots+C_{l_n}$.
 Therefore, in view of~\thetag{1-3}, if $l_n$ is chosen sufficiently large then
 $F'_{l_n'}+B_n\subset F_{l_n}$.
 Thus we defined $l_n$ and $B_n$ satisfying \thetag{2-2}.
 In a similar way on the next step we construct $l_{n+1}'$ and $A_{n+1}$.
 Continuing this construction  infinitely many times we define the entire sequences
  $(l_n)_{n\ge 0}$, $(l_n')_{n\ge 1}$, $(A_n)_{n\ge 1}$ and $(B_n)_{n\ge 1}$ satisfying \thetag{2-2} as desired.
 \qed
 
 \comment

 Hence, for each $n>0$, the polynomial $P_{C_{l_{n-1}+1}}P_{C_{l_{n-1}+2}}\cdots P_{C_{l_n'}}$ divides the righthand side of~\thetag{2-5}.
 The factrorization property of polynomials with integer coefficients yields that there is $l_n>l_n'$
 and a polynomial $Q_n$ with integer coefficients such that
 $$
 Q_nP_{C_{l_{n-1}+1}}\cdots P_{C_{l_n'}}=P_{A_n}P_{C_{l_n'+1}'}\cdots P_{C_{l_n}'}.\tag2-6
 $$
 From this and \thetag{2-5} we deduce that
 $$
  P_{C_{l_n'+1}}P_{C_{l_n'+2}}\cdots=Q_nP_{C_{l_n+1}'}P_{C_{l_n+2}'}\cdots.\tag2-7
 $$
 Since $\deg Q_n<\deg (P_{A_n}P_{C_{l_n'+1}'}\cdots P_{C_{l_n}'})<\min C'_{l_n+1}\setminus\{0\}$,
 it follows from \thetag{2-7} that there is 
 a  finite subset $A_{n+1}\subset C_{l_n'+1}+C_{l_n'+2}+\cdots$ such that $0\in A_{n+1}$ and $Q_n=P_{A_{n+1}}$.
 Hence $(A_{n+1}-A_{n+1})\cap (F_{l_n'}-F_{l_n'})=\{0\}$. 
 Since $F_{l_n'}$ is an  interval (in $\Bbb Z$), it follows from~\thetag{2-6} that $A_n\subset F_{l_n}'$.
 \endcomment
\enddemo

The conditions for isomorphism of $(C,F)$-actions in Theorem~2.6 are especially simple in the case where $T$ and $T'$ are {\it commensurate},
i.e.  $F_n=F_n'$ eventually.

\proclaim{Theorem 2.8}
Let $(G,G_+)$ be a linearly  ordered discrete countable Abelian group.
Let $T=(T_g)_{g\in G}$ and $T'=(T'_g)_{g\in G}$ be two $(C,F)$-actions of $G$ associated with some sequences $(C_n,F_{n-1})_{n\ge 1}$ and $(C_n',F_{n-1}')_{n\ge 1}$ respectively and the two sequences satisfy (I)--(III) and \thetag{1-3}. 
Suppose that $C_n\cup C_n'\subset G_+$  for all $n$ and $F_n=F_n'$ for all $n>N$ (for some $N>0$).
Then $T$ and $T'$ are
topologically isomorphic if and only if
 $C_n=C_n'$ for all $n>M$
(for some $M>0$).
\endproclaim

\demo{Proof} Suppose that $T$ and $T'$ are
topologically isomorphic. 
By Theorem~2.5 and~\thetag{2-1}, for each $n>N$ and $b,\widetilde b\in B_n$, we have
$$
A_n+b\subset F_{l_n'}+b,\quad (F_{l_n'}+b)\cap(F_{l_n'}+\widetilde b)=\emptyset\quad\text{if}\quad b\ne\widetilde b\quad \text{and }0\in B_n.
$$
Therefore
$$
A_n=(A_n+B_n)\cap F_{l_n'}=(C_{l_{n-1}+1}+\cdots +C_{l_n})\cap F_{l_n'}=C_{l_{n-1}+1}+\cdots +C_{l_n'}.
$$ 
This implies, in turn,  that 
$$
\bigsqcup_{b\in B_n}(A_n+b)=\bigsqcup_{c\in C_{l_n'+1}+\cdots +C_{l_n}}(A_n+c).\tag2-8
$$
Hence $\# B_n=\#(C_{l_n'+1}+\cdots +C_{l_n})$.
Moreover, \thetag{2-8} yields that there are two maps $B_n\ni b\mapsto \gamma(b)\in C_{l_n'+1}+\cdots +C_{l_n}$ and $B_n\ni b\mapsto\alpha(b)\in A_n$ such that
$b=\alpha(b)+\gamma(b)$ for each $b\in B_n$.
We first note that the map $\gamma$ is one-to-one.
Indeed, if $\gamma(b_1)=\gamma(b_2)$ then $b_1-\alpha(b_1)=b_2-\alpha(b_2)$ and hence $b_1-b_2\in F_{l_n'}-F_{l_n'}$.
The later implies that $b_1=b_2$, as desired.
Hence $\gamma$ is a bijection of $B_n$ onto $C_{l_n'+1}+\cdots +C_{l_n}$ (we recall that these finite sets are of the same cardinality).
Comparing the sums of the elements in the lefthand side and the righthand side of \thetag{2-8} we obtain that 
$$
\#B_n\sum_{a\in A_n}a+\sum_{b\in B_n}(\alpha(b)+\gamma(b))=\#(C_{l_n'+1}+\cdots +C_{l_n})\sum_{a\in A_n}a+
\sum_{c\in C_{l_n'+1}+\cdots +C_{l_n}}c.
$$
Hence $\sum_{b\in B_n}\alpha(b)=0$.
Since $\alpha(b)\in A_n\subset G_+$, it follows that $\alpha(b)=0$ for each $b\in B_n$.
Hence $B_n=C_{l_n'+1}+\cdots +C_{l_n}$.
In a similar way we derive from the second equality in the first line of \thetag{2-1} that
$B_n=C_{l_n'+1}'+\cdots +C_{l_n}'$ and $A_{n+1}=C_{l_{n}+1}'+\cdots +C_{l_{n+1}'}'$
for all sufficiently large $n$.
Therefore
 $$
 \align
 C_{l_{n}'+1}+\cdots +C_{l_n}&=C'_{l_{n}'+1}+\cdots +C'_{l_n}\quad\text{and}\\
  C_{l_{n}+1}+\cdots +C_{l_{n+1}'}&=C_{l_{n}+1}'+\cdots +C_{l_{n+1}'}'.
  \endalign
 $$
It follows now from (III) and the equality $F_j=F_j'$ eventually that
 $C_j=C_j'$ eventually, as claimed.
 
 The converse implication follows from Corollary~2.2.
 \qed
\enddemo

\subhead 2.3. Applications to  topological centralizers and inverse actions 
\endsubhead
Given a topological action $T$ of $G$ on a locally compact Cantor space $X$, we let
$$
C_{\text{top}}(T):=\{\theta\in\text{Homeo}(X)\mid \theta T_g=T_g\theta\text{ for each }g\in G\}
$$
and call this set {\it the topological centralizator} of $T$.

\comment
\proclaim{Lemma 2.3} Let $A$ be a finite subset of $\Bbb Z_+$ and $0\in A$.
Then there is a unique (up to permutation) decomposition $A=A_1\dotplus\cdots\dotplus A_k$ of $A$ into a sum of independent subsets $A_1,\dots,A_k\subset A$ and $0\in A_1\cap\cdots \cap A_k$.
\endproclaim
\demo{Proof}
For each finite subset $B\subset \Bbb Z_+$, we define a polynomial  $p^B\in\Bbb Z_+[\boldsymbol X]$ by setting $p^B(\boldsymbol X):=\sum_{b\in B}\boldsymbol X^b$.
Since the ring $\Bbb Z_+[\boldsymbol X]$ is factorial, there is a unique (up to permutation) decomposition $p^B=p_1\cdots p_k$ of $p_B$ into product of irreducible polynomials $p_1,\dots,p_k\in\Bbb Z_+[\boldsymbol X]$.
Since each coefficient of  $p^B$  is either $0$ or $1$, it follows that each coefficient of $p_i$ is either $0$ or $1$ for $i=1,\dots,k$.
Hence there is a finite subset $A_i\subset\Bbb Z_+$ such that $p_i=p^{A_i}$  for each $i=1,\dots,k$.
Moreover, $p^B(1)=1$ implies that  $p_i(1)=1$ for each $i=1,\dots,k$, i.e. $0\in A_1\cap\cdots \cap A_k$.
\qed
\enddemo
\endcomment

\proclaim{Corollary 2.9} Let $(G,G_+)$ be a linearly  ordered discrete countable Abelian group.
Let $T$ be a $(C,F)$-action of $G$ associated with a sequence $(C_n, F_{n-1})_{n\ge 1}$ 
satisfying ~(I)--(III) and \thetag{1-3}.
If $C_n\subset G_+$ for each $n\ge 1$ then 
 $C_{\text{top}}(T)=\{T_g\mid g\in G\}$.
\endproclaim
\demo{Proof}
Let $\phi\in C_{\text{top}}(X)$.
Using $\phi$ we can construct the sequences
$0=l_0<l_1'<l_1<l_2'<l_2<\cdots$ and  $(A_n)_{n\ge 1}$ and $(B_n)_{n\ge 1}$ as in the proof of Theorem~2.3.
In particular, $\phi([0]_{l_n})=[A_{n+1}]_{l_{n+1}'}$ for each $n$.
It was shown in the proof of Theorem~2.5 that $\min A_n=-\min B_n=\min A_{n+1}$ for each $n>0$.
Hence there is $g_0\in G$ such that 
$$
g_0=\min A_1=\min A_2=\cdots=-\min B_1=-\min B_2=\cdots.
$$
We now set $\psi:=T_{-g_0}\circ\phi$.
Then
 $$
 \psi([0]_{l_n})=T_{-g_0}[A_{n+1}]_{l_{n+1}'}=[A_{n+1}-g_0]_{l_{n+1}'}=[A_{n+1}-\min A_{n+1}]_{l_{n+1}'}\supset[0]_{l_{n+1}'}
 $$
  for each $n$.
  We let $\bold{0}:=(0,0,\dots)\in X_0\subset X$.
  Hence 
  $$
  \psi(\bold{0})=\psi\left(\bigcap_{n=1}^\infty[0]_{l_n}\right)=\bigcap_{n=1}^\infty\psi([0]_{l_n})\supset\bigcap_{n=1}^\infty[0]_{l_{n+1}'}=\bold{0}.
  $$
  It follows that $ \psi(\bold{0})=\bold 0$.
Since $\psi$ is equivariant, we obtain that $\psi(x)=x$ for each $x$ from the $T$-orbit of $\bold 0$. 
Thus $\psi$ and Id coincide on a dense subset of $X$.
Therefore $\psi=\text{Id}$ and hence $\phi=T_{g_0}$.
\qed

\comment
By Theorem~2.2, there is an increasing sequence of integers
$0=l_0<l_1'<l_1<l_2'<l_2<\cdots$ and subsets $A_n \subset  F_{l_n'}$ and   $B_n \subset  F_{l_n}$, 
such that 
$$
A_{n}+B_{n}=C_{l_{n-1}+1}+\cdots+ C_{l_{n}},\quad
B_{n}+ A_{n+1}=C_{l_{n}'+1}+\cdots +C_{l_{n+1}'},\tag2-2
$$
$F_{l_n'}+B_n\subset F_{l_n}$,  $F_{l_n}+A_{n+1}\subset F_{l_{n+1}}$ and 
$$
(F_{l_n'}-F_{l_n'})\cap (B_n-B_n)=(F_{l_n}-F_{l_n})\cap (A_{n+1}-A_{n+1})=\{0\}\tag2-3
$$
 for each $n>0$.
We let   $a_n:=\min A_n$ and $b_n:=\min B_n$ and set
  $\widetilde A_n:=A_n-a_n$ and $\widetilde B_n:=B_n-b_n$ for each $n>0$.
 Of course, $\widetilde A_n\cup \widetilde B_n\subset G_+$.
 In view of the remark at the beginning of the proof, we may assume without loss of generality (increasing $l_n$ and $l_n'$ in the proof of Theorem~2.3, if necessary) that $\widetilde A_n\subset F_{l_n'}$ and  $\widetilde B_n \subset  F_{l_n}$. 
In follows from \thetag{2-2} and the condition of the theorem  that
$\widetilde A_n +(B_n+a_n)=A_n+B_n\subset G_+$.
Since $0\in \widetilde A_n$, we obtain that $B_n+a_n\subset G_+$.
Moreover, since $0\in A_n+B_n$, it follows that $0\in B_n+a_n$.
This yields that $a_n=-b_n$ because $G$ is linearly ordered.
In a similar way, exploiting the second equality in \thetag{2-2}, we obtain that $a_{n+1}=-b_n$ for each $n>0$.
Hence $a_1=a_2=\cdots=-b_1=-b_2=\cdots=:g_0$.
Moreover, \thetag{2-2} and~\thetag{2-3} hold if we replace $A_n$ and $B_n$ with $\widetilde A_n$ and $\widetilde B_n$.
Therefore we obtain that
  $$
  \align
  \widetilde A_n&=F_{l_n'}\cap ( \widetilde A_n+\widetilde B_n)=F_{l_n'}\cap(C_{l_{n-1}+1}+\cdots+ C_{l_{n}})=C_{l_{n-1}+1}+\cdots+ C_{l_{n}'}\text{ and}\\
  \widetilde B_n&=F_{l_n}\cap ( \widetilde B_n+\widetilde A_{n+1})=F_{l_n}\cap(C_{l_{n}'+1}+\cdots+ C_{l'_{n+1}})=C_{l_{n}'+1}+\cdots+ C_{l_{n}}.
  \endalign
  $$
  We let $\psi:=T_{-g_0}\circ\phi$.
  Since $[A_{n}]_{l_n'}=\phi([0]_{l_{n-1}})$ and $[B_n]_{l_n}=\phi^{-1}([0]_{l_n'})$, we obtain that  $[\widetilde A_{n}]_{l_n'}=\psi([0]_{l_{n-1}})$ and $[\widetilde B_n]_{l_n}=\psi^{-1}([0]_{l_n'})$ for each $n>0$.
  It follows from the proof of Theorem~2.2 that $\psi$ is the identity. (?)
  \endcomment
\enddemo

  \comment
  and hence 
  $$
   \widetilde A_n+\widetilde B_n= \widetilde A_n+(C_{l_{n}'+1}+\cdots+ C_{l_{n}}).\tag2-3
   $$
 This yields that for each $b\in \widetilde B_n$, there exist $a_b\in \widetilde A_n$ and $c_b\in C_{l_{n}'+1}+\cdots+ C_{l_{n}}$ such that $b=a_b+c_b$ (we use the fact that  $0\in \widetilde A_n$).
 Moreover, the map $ \widetilde B_n\ni b\mapsto c_b\in C_{l_{n}'+1}+\cdots+ C_{l_{n}}$
 is a bijection in view of \thetag{2-2}.
 Comparing the sums of elements from the lefthand side and the righthand side of \thetag{2-3},  we obtain that $\sum_{b\in\widetilde B_n}a_b=0$.
 This yields that $a_b=0$ for each $b\in\widetilde B_n$.
 Hence $\widetilde B_n=C_{l_{n}'+1}+\cdots+ C_{l_{n}}$ for each $n>0$.
 Thus 
 $$
 \widetilde A_{n}+\widetilde B_{n}=C_{l_{n-1}+1}+\cdots+ C_{l_{n}},\quad
 \widetilde B_{n}+  \widetilde A_{n+1}=C_{l_{n}'+1}+\cdots +C_{l_{n+1}'}
 $$
 
 $a_n=a_{n+1}?$.
 \endcomment

We note that, in particular, the odometers have trivial (topological) centralizer in their locally compact $(C,F)$-realizations.

\comment
We say that a subset $D$ in $G$ is symmetric if there is $w\in G$ such that $w-D=D$.
If $T$ is a $(C,F)$-action associated with a sequence $(C_n,F_{n-1})_{g\ge 1}$.
Suppose that each $F_n$ is symmetric, i.e. there is $d_n\in G$ such that $d_n-F_n=F_n$.
Denote by $T^*=(T^*_g)_{g\in G}$ the $G$-action (on the same space) which is inverse to $T$, i.e.
$T^*_g:=T_{-g}$. 
It is easy to verify that $T^*$ is a $(C,F)$-action of $G$ associated with the sequence $(d_n-d_{n-1}-C_n, F_{n-1})_{n\ge 1}$.
The isomorphism is
$$
(f_n,c_{n+1}, c_{n+2},\dots)\mapsto (d_n-f_n, d_{n+1}-d_n-c_{n+1},d_{n+2}-d_{n+1}-c_{n+2},\dots).
$$
\endcomment


Given a continuous   action $T=(T_g)_{g\in G}$ on  an Abelian group $G$ on  a topological space $X$, we let  $T^{-1}:=(T_{-g})_{g\in G}$.
Then $T^{-1}$ is also a continuous action of $G$ on $X$.
We call it  the {\it inverse to $T$}.

\proclaim{Lemma 2.10}
Let $(G,G_+)$ be a linearly  ordered discrete countable Abelian group.
Let  $T$ be a $(C,F)$-action of $G$ associated with a sequence $(C_n,F_{n-1})_{n\ge 1}$ satisfying (I)--(III) and \thetag{1-3}.
Then    $T^{-1}$ is (topologically) isomorphic a $(C,F)$-action
 associated with the sequence
$(C_n^*,F_{n-1}^*)_{n\ge 1}$, where $C_n^*:=\{-c+\max C_n\mid c\in C_n\}$ and
$F_n^*:=\{-f+\sum_{j=1}^n\max C_j\mid f\in F_n\}$.
\endproclaim
\demo{Sketch  of the proof}
It is straightforward to verify that (I)--(III) and \thetag{1-3} are satisfied for the sequence
$(C_n^*,F_{n-1}^*)_{n\ge 1}$.
The canonical  isomorphism of $T$ with $T^{-1}$ is given by the map 
$$
(f_n,c_{n+1},\dots)\mapsto\left(\sum_{j=1}^n\max C_j-f_n,\max C_{n+1}-c_{n+1},\max C_{n+2}-c_{n+2},\dots\right)
$$
from $X_n:=F_n\times C_{n+1}\times C_{n+2}\times\cdots$ onto $X_n^*:=F_n^*\times C_{n+1}^*\times C_{n+2}^*\times\cdots$  for each $n\ge 0$.
\enddemo

\proclaim{Corollary 2.11}
Let $(G,G_+)$ be a linearly  ordered discrete countable Abelian group.
Let $T$ be a  $(C,F)$-action associated with a sequence $(C_n,F_{n-1})_{n\ge 1}$ satisfying (I)--(III) and \thetag{1-3}.
Suppose that $\bigcup_{n\ge 1}C_n\subset G_+$ and  that  (II) and \thetag{1-3} are satisfied for  the sequence $(C_n^*,F_{n-1})_{n\ge 1}$.
Then $T$ is topologically isomorphic to $T^{-1}$ if and only if $C_n=C_n^*$ eventually.
\endproclaim

\demo{Proof}
Of course, (I) and (III) are also satisfied for $(C_n^*,F_{n-1})_{n\ge 1}$.
Hence a $(C,F)$-action associated with $(C_n^*,F_{n-1})_{n\ge 1}$ is well defined.
By Corollary~2.2,   this action  is isomorphic to $T^{-1}$.
 Now Theorem~2.8 yields that $T$ and $T^{-1}$ are isomorphic if and only if  $C_n=C_n^*$ eventually.
\qed
\enddemo

\remark{Remark 2.12} We note that  the conditions of Corollary~2.11 are satisfied for the most important example where $(G,G_+)=(\Bbb Z,\Bbb Z_+)$ and every set $F_n$ is an {\it order interval}  $\{j\in\Bbb Z\mid -\alpha_n\le j\le\beta_n\}$ for some integers $\alpha_n,\beta_n>0$.
This follows from the  following simple fact: if $A$ is a finite subset of $\Bbb Z$ such that $\min A\in F_n$ and $\max A\in F_n$ then $A\subset F_n$.
It remains to notice that $\min C_n=\min C_n^*=0$ and $\max C_n=\max C_n^*$ for each $n>0$.
\endremark

\comment

This set is meager?

Let $\goth R$ denote the set of all sequences $(C_n,F_n)_{n=1}^\infty$ of finite subsets in $G$ satisfying (I)-(II).
Moreover, we will assume that each $C_n$ is not decomposable.
Then $\goth R$ is  a parametrization of all rank one transformations.
Of course, $\goth R$ is a Polish space in the natural product  topology.
Moreover, it is totally disconnected.

We now want to understand the complexity of the isomorphism 
relation on
$\goth R$
as a Borel equivalence relation.

\proclaim{Proposition} There is a Borel map  $\Phi$ from $\{0,1\}^\Bbb N$ into $\goth R$ such that $\Phi(a)$ and $\Phi(b)$ are topologically conjugate if and only if  $a$ and $b$ are tail equivalent.
\endproclaim
\demo{Proof} We set $h_n:=3h_n+3$ and 
 $C_n^{(0)}:=\{0,h_n,2h_n+1\}$ and $C_n^{(1)}:=\{0,h_n+1,2h_n+1\}$, $F_n:=\{-n,-n+1,\dots, 3h_n-n+2\}$.
 Let $\Phi (a)$ be the sequence $(C_n^{(a_n)}, F_n)_{n=1}^\infty$.
 \qed
\enddemo

\proclaim{Proposition} The isomorphism equivalence relation on $\goth R$ is hyperfinite and $F_\sigma$.
\endproclaim
\demo{Proof} Given $N$, define $\goth R_N$ by setting $C_i+l=C_{i+m}'$  or  $C_i'=C_{i+m}$ for some $0\le l,m\le N$ and  all $i>N$.
Then $\goth R_1\subset\goth R_2\subset\dots$ and $\bigcup_{n>0}\goth R_N=\goth R$.
Each $\goth R_N$ is a finite equivalence relation.
\qed
\enddemo

It follows from \cite{DoJaKe}.

\proclaim{Corollary}
The isomorphism relation on
$\goth R$
is Borel bi-reducible with
$E_0$.
\endproclaim

\endcomment

\head 3. $(C,F)$-models for~measure preserving~actions~of monotileable amenable groups
\endhead

Let $G$ be a {\it monotileable} amenable discrete countable group \cite{We}. 
This means that there is a F{\o}lner sequence $(\Cal F_n)_{n\ge 0}$ in $G$ such that every set $\Cal F_n$ tiles $G$, i.e. there is a subset $\Cal C_n\subset G$ such that $\Cal F_n\Cal C_n=G$ and $\Cal F_nc\cap \Cal F_nc'=\emptyset$ whenever $c,c'\in \Cal C_n$ and 
$c\ne c'$.
Without loss of generality we may assume that $1\in \Cal F_n\cap \Cal C_n$ for each $n\ge 0$. 

Fix a standard nonatomic probability space $(X,\goth B,\mu)$.
Denote by Aut$(X,\mu)$ the group of  $\mu$-preserving transformations. 
It is a Polish group when endowed with the weak operator topology.
Let $\Cal A_G$ stand for the  set of  $\mu$-preserving actions of $G$ on $X$.
An element  of $\Cal A_G$ is a group homomorphisms from  $G$ to Aut$(X,\mu)$. 
Thus 
 $\Cal A_G$ is a subset of the infinite product space Aut$(X,\mu)^G$.
Endow the latter space with the infinite product of the weak operator topologies on Aut$(X,\mu)$.
Then  Aut$(X,\mu)^G$ is a Polish space and  $\Cal A_G$ is a closed subset of it.
Hence $\Cal A_G$ is Polish in the induced topology.

Let $\goth F$ denote the set of all finite subsets of $G$.
We now let
$$
\multline
\goth R_1:=\{(C_n,F_{n-1})_{n\ge 1}\in(\goth F\times\goth F)^\Bbb N\mid (F_n)_{n\ge 0}\text{ is a subsequence of $(\Cal F_n)_{n\ge 0}$  } \\
\text{ and (I)--(III) and \thetag{1-3} are satisfied}\}.
\endmultline
$$
Endow $\goth F$ with the discrete topology.
Then the space $(\goth F\times\goth F)^\Bbb N$ furnished with the infinite product topology---we denote it by
$\tau$---is a Polish 0-dimensional space.

\proclaim{Lemma 3.1}
The subset $\goth R_1$ is a $G_\delta$-subset  of  $(\goth F\times\goth F)^\Bbb N$.
\endproclaim
\demo{Proof}
It is clear that each of the conditions (I)--(III) determines a closed subset in $(\goth F\times\goth F)^\Bbb N$.
The condition  \thetag{1-3} determines a $G_\delta$-subset in $(\goth F\times\goth F)^\Bbb N$.
\qed
\enddemo

Hence $\goth R_1$ endowed with $\tau$ is  Polish and 0-dimensional.
The function
$$
\phi_m:\goth R_1\ni (C_n,F_{n-1})_{n\ge1}\mapsto \frac{\# F_m}{\#C_1\cdots \#C_m}\in\Bbb R
$$
is continuous on $\goth R_1$ for each $m>0$.
Moreover, $\phi_1\le\phi_2\le\cdots$. 
Hence the extended function $\phi:=\sup_{n>0}\phi_n$ taking values in $\Bbb R\cup\{+\infty\}$ is lower semicontinuous. 
We note that $\phi((C_n,F_{n-1})_{n\ge1})= \lim_{m\to\infty}\frac{\# F_m}{\#C_1\cdots \#C_m}$.
We now set
$$
\goth R_1^{\text{fin}}:=\{\Cal S \in\goth R_1\mid \phi(\Cal S)<\infty \}.
$$
Then $\goth R_1^{\text{fin}}$ is an $F_\sigma$-subset of $\goth R_1$.

\definition{Definition 3.2} Denote  by $\tau^{\text{fin}}$ the weakest topology on $\goth R_1^{\text{fin}}$ which is stronger than $\tau$ and such that $\phi$ is continuous in it.
\enddefinition

It is well known that  $\tau^{\text{fin}}$ is Polish (see, e.g. \cite{Ke}).
We note that a sequence $(\Cal S_n)_{n=1}^\infty$ in $\goth R_1^{\text{fin}}$ converges
in $\tau^{\text{fin}}$ to an element $\Cal S\in \goth R_1^{\text{fin}}$   if and only if $(\Cal S_n)_{n=1}^\infty$ converges in $\tau$ to $\Cal S$  as $n\to \infty$ and  $\lim_{n\to\infty}\phi(\Cal S_n)=\phi(\Cal S)$.

Our purpose now is to  construct a continuous map 
$$
\Psi:\goth R_1^{\text{fin}}\ni\Cal S\mapsto\Psi^{\Cal S}\in\Cal A_G.
$$
Without loss of generality we may assume that $X=[0,1)$ and $\mu$ is the Lebesgue measure on $[0,1)$. 
Fix a linear order $\succ$ on $G$.
Given $\Cal S:=(C_n,F_{n-1})_{n\ge 1}\in\goth R_1^{\text{fin}}$, we denote by  $T^{\Cal S}=(T^{\Cal S}_g)_{g\in G}$
 the $(C,F)$-action of $G$ associated with $\Cal S$.
 Let $X^\Cal S$ and $\mu^\Cal S$ stand for  the $(C,F)$-space and the $(C,F)$-measure of $T^\Cal S$ respectively.
 For each $n> 0$,
we set 
$$
\alpha_0:=\frac 1{\phi(\Cal S)}\quad\text{and}\quad\alpha_n:=\frac{\alpha_0\# F_n}{\# C_1\cdots\# C_n}.
$$ 
Then the sequence $(\alpha_n)_{n=1}^\infty$ increases and converges to $1$ as $n\to\infty$.
For each $n\ge 0$, we partition the interval $[0,\alpha_n)$ into subintervals $I^{(n)}_f$, $f\in F_n$, of length $\alpha_n/\# F_n$ such that
\roster
\item"$\bullet$"
every subinterval $I^{(n)}_f$, $f\in F_n$, is the union of  pairwise disjoint subintervals $I^{(n+1)}_{fc}$, $c\in C_{n+1}$,
\item"$\bullet$" if $f\in F_n$, $c,c'\in C_{n+1}$ and $c\succ c'$ then $I^{(n)}_{fc}$ is on the right of $I^{(n)}_{fc'}$,
\item"$\bullet$" 
if $f,f'\in F_{n+1}\setminus(F_nC_{n+1})$ and $f\succ f'$ then $I^{(n)}_{f}$ is on the right of $I^{(n)}_{f'}$.
\endroster
It is obvious that these conditions determine the partition $[0,\alpha_n)=\bigsqcup_{g\in F_n}I^{(n)}_g$ and the enumeration of its atoms by elements of $F_n$ in a unique way for each $n\ge 0$.
By the standard properties of Lebesgue spaces, there is a unique (mod 0) Borel bijection 
$\beta:X^S\to [0,1)$ such that
 $\beta([f]_n)=I^{(n)}_f$ (mod 0) for each $f\in F_n$ and $n\ge 0$
and   $\mu^{\Cal S}(X^S)^{-1}\mu^\Cal S=\mu\circ\beta$.
We now define the $G$-action $\Psi^{\Cal S}=(\Psi^{\Cal S}_g)_{g\in G}$ on $X$ by setting
$\Psi^{\Cal S}_g:=\beta T^{\Cal S}_g\beta^{-1}$, $g\in G$.
Thus, if $g\in G$ and $f\in F_n$ with $gf\in F_n$  then $\Psi^{\Cal S}_gI^{(n)}_f=I^{(n)}_{gf}$,
$n\in\Bbb N$.

\proclaim{Lemma 3.3} $\Psi$ is continuous.
\endproclaim
\demo{Proof} If $\Cal S=(C_n,F_{n-1})_{n\ge 1}$ and $\Cal S'=(C_n',F_{n-1}')_{n\ge 1}$ are $\tau^{\text{fin}}$-close, there is $\epsilon>0$ and $N>0$ such that 
$(C_n,F_{n-1})=(C_n',F_{n-1}')$ for each $n=1,\dots,N$ and $\phi(\Cal S)=(1\pm\epsilon)\phi(\Cal S')$.
Let $I^{(n)}_f$, $f\in F_n$, denote the subintervals of $[0,1)$ used in the definition of $\Psi^{\Cal S}$ and 
let $J^{(n)}_f$, $f\in F_n$, denote the  similar subintervals of $[0,1)$ used to define
 $\Psi^{\Cal S'}$ for $n=0,1,\dots,N$. 
 Then for each $n=1,\dots, N-1$ and $f\in F_n$, we have $\mu(I^{(n)}_f\triangle J^{(n)}_f)\le2\epsilon\mu(I^{(n)}_f)$.
 Hence for
 each $g\in G$ and $f\in F_n$ such that $gf\in F_n$, we obtain that
 $$
 \mu(\Psi^{\Cal S}_gI^{(n)}_f\triangle\Psi^{\Cal S'}_gI^{(n)}_f)\le \mu(I^{(n)}_{gf}\triangle J^{(n)}_{gf})+ 2\epsilon\mu(I^{(n)}_f)\le 4\epsilon\mu(I^{(n)}_f).
 $$
A standard argument yields now that $\Psi^{\Cal S}$ and $\Psi^{\Cal S'}$ are close in $\Cal A_G$ if $\epsilon$ is small enough and $N$ is large enough. \qed
\enddemo

Denote by $\boldsymbol  R$ the {\it tail equivalence relation} on $(\goth F\times\goth F)^\Bbb N$.
This means that two sequences $(D_n)_{n\ge 1}$ and $(D_n')_{n\ge 1}$ from $(\goth F\times\goth F)^\Bbb N$ are $\boldsymbol  R$-equivalent if $D_n=D_n'$ eventually.
 It is easy to see that  $\boldsymbol R$ is an $F_\sigma$-subset of $(\goth F\times\goth F)^\Bbb N\times(\goth F\times\goth F)^\Bbb N$.
 It is straightforward to verify that each $\boldsymbol  R$-class is dense in $(\goth F\times\goth F)^\Bbb N$.
 Of course, the subsets $\goth R_1$ and $\goth R_1^{\text{fin}}$ are  $\boldsymbol R$-saturated.
 Hence $\boldsymbol  R$ restricted to $\goth R_1^{\text{fin}}$ is also $F_\sigma$ in $\tau^{\text{fin}}$.
 It is easy to see that if $\Cal S$, $\Cal S'\in \goth R_1^{\text{fin}}$ and $(\Cal S$, $\Cal S')\in\boldsymbol  R$ then the $G$-actions $\Psi^{\Cal S}$ and $\Psi^{\Cal S'}$ are isomorphic.

\proclaim{Lemma 3.4} For each  $\Cal S\in \goth R_1^{\text{fin}}$,  the  $\boldsymbol  R$-class of $\Cal S$ is $\tau^{\text{fin}}$-dense in 
$\goth R_1^{\text{fin}}$.
\endproclaim

\demo{Proof} Let   $\Cal S=(C_{n},F_{n-1})_{n\ge 1}$.
 Take an arbitrary $\widehat{\Cal S}=(\widehat C_{n},\widehat F_{n-1})_{n\ge 1}\in\goth R^{\text{fin}}_1$.
We will construct   $\Cal S'\in\goth R^{\text{fin}}_1$ which is $\tau^{\text{fin}}$-close (as close as we wish)   to $\widehat{\Cal S}$ and such that $(\Cal S',\Cal S)\in \boldsymbol  R$.

Fix $\epsilon>0$.
Select  $j>0$ such that $\phi_j(\widehat{\Cal S})>\phi(\widehat{\Cal S})-\epsilon$.
For $i>0$, let 
$$
F_i^\circ:=\{f\in F_i\mid \widehat F_j\widehat F_j^{-1}f\subset F_i\}.
$$
Since $(F_n)_{n\ge 1}$ is a F{\o}lner sequence in $G$,
 there exists
 $i>0$  such  that $\# F_i^\circ >(1-\epsilon)\# F_i$.
 Since $\phi(\Cal S)<\infty$, we can also assume without loss of generality that
 $$
 \frac{\# F_{m}}{\# F_i\#C_{i+1}\cdots \# C_{m}}<1+\epsilon\tag3-1
 $$
 for each $m>i$.
 Let $C:=\{c\in\Cal C_j\mid F_i^\circ\cap \widehat F_jc\ne\emptyset\}$.
 Then $F_i\supset \widehat F_jC\supset F_i^\circ$ and $\widehat F_jc\cap \widehat F_jc'=\emptyset$ if $c,c'\in C$ and $c\ne c'$.
 It follows that
 $$
{\# F_i}<(1+\epsilon){\# C\# \widehat F_j}.\tag3-2
 $$
 We now set
 $$
 (C_{a+1}',F_a'):=
 \cases
 (\widehat C_{a+1},\widehat F_a) &\text{if }a<j,\\
 (C,\widehat F_a) &\text{if }a=j,\\
 (C_{a-j+i},F_{a-j+i-1})& \text{if }a>j.
 \endcases 
 \tag3-3
 $$
 Let $\Cal S':=(C_{a+1}',F_a')_{a\ge 0}$.
 It is obvious that (I)--(III) are satisfied for $\Cal S'$ for each $a>0$.
 We  note that \thetag{2-3} holds for $\Cal S'$ because it holds for $\Cal S$.
 Hence $\Cal S'\in\goth R_1$.
 For each $m>j$, we have
 $$
 \phi_m(\Cal S'):=\frac{\# F'_m}{\#C'_1\cdots\#C'_m}=\phi_j(\widehat{\Cal S})\cdot \frac{\# F_i}{\# C\#\widehat F_j}\cdot\frac{\# F_{m+i-j}}{\# F_i\#C_{i+1}\cdots \# C_{m+i-j}}.
 $$
 Applying \thetag{3-1} and \thetag{3-2} we obtain that
  $\phi_m(\Cal S')=(\phi(\Cal S)\pm\epsilon)(1\pm \epsilon)^2$.
It follows that $\Cal S'\in\goth R_1^{\text{fin}}$ and $\phi(\Cal S')$ is close to $\phi(\Cal S)$.
  The first line in  \thetag{3-3} yields that $\Cal S'$ is $\tau$-close to $\Cal S$ (if $j$ is chosen large).
  Hence $\Cal S'$ is $\tau^{\text{fin}}$-close to $\Cal S$.
  It follows from the third line in \thetag{3-3} that $(\Cal S',\Cal S)\in\boldsymbol R$, as desired.
  \qed
  \enddemo

We now generalize  the concept of model for Aut$(X,\mu)$ (viewed as the $\Bbb Z$-actions on $(X,\mu)$) introduced in \cite{Fo, Definition~10} to the $G$-actions on $(X,\mu)$.

\definition{Definition 3.5} A {\it model} for  $\Cal A_G$  is a pair $(W,\pi)$, 
where $W$ is a Polish space and $\pi:W\to\Cal A_G$ is a continuous map such that for a comeager set $\Cal M\subset\Cal A_G$ and each $A\in\Cal M$, the set $\{w\in W\mid \pi(w)\text{ is isomorphic to } A\}$ is dense in $W$.
\enddefinition

We let $\Cal O_G:=\{T\in\Cal A_G\mid\text{$T$ is of rank one along $(\Cal F_n)_{n\ge 0}$}\}$.

\proclaim{Proposition 3.6} $\Cal O_G$ is a dense $G_\delta$ (and hence comeager) in
$\Cal A_G$.
\endproclaim
\demo{Proof}
We first prove that $\Cal O_G$ is a $G_\delta$.
Without loss of generality we may assume that $X$ is a compact Cantor space.
Let $\Cal K$ denote the class of clopen subsets in $X$.
Since $\Cal K$ is countable, we can write it as $\Cal K=\{K_i\mid i\in\Bbb N\}$.
Given $A\in \Cal K$, a finite subset $F\subset G$, $m>0$ and an action  $T=(T_g)_{g\in G}\in\Cal A_G$, we set
$$
\align
\alpha_{m,A,F}(T)&:=\max_{1\le i\le m}\min_{F'\subset F}\mu\left(K_i\triangle \bigcup_{g\in F'}T_gA\right)\\
\beta_{A,F}(T)&:=\max_{g\in F}\mu(T_gA\cap A)
\endalign
$$
Then $ \alpha_{m,A,F}:\Cal A_G\ni T\mapsto \alpha_{m,A,F}(T)\in \Bbb R$ and 
$\beta_{A,F}:\Cal A_G\ni T\mapsto \beta_{A,F}(T)\in \Bbb R$ are continuous functions.
It is straightforward to verify that $\Cal O_G$ equals 
$$
\bigcap_{r=1}^\infty\bigcap_{m=1}^\infty\bigcap_{n=1}^\infty\bigcap_{l=1}^\infty\bigcup_{s>r}\bigcup_{A\in\Cal K}
\left\{T\in\Cal A_G\,\bigg|\,\alpha_{m,A,\Cal F_s}(T)<\frac{\mu(A)}{l\#\Cal F_s}\text{ and }\beta_{A,\Cal F_s}(T)<\frac{\mu(A)}{l\#\Cal F_s}\right\}.
$$
Hence it is a $G_\delta$.

By \cite{FoWe, Claim~18} (see also \cite{Da3, just below Proposition~1.2} for the proof),
given a free  $T\in \Cal A_G$, the class of $G$-actions on $(X,\mu)$ that are isomorphic to $T$
 is dense in $\Cal A_G$.
Hence it remains to show that $\Cal O_G$ contains a free $G$-action.
For that we will utilize the $(C,F)$-construction.
As in the proof  of Lemma~3.4,
we can construct inductively the sequence $(C_n,F_{n-1})_{n\ge 1}$ such that
$F_n=\Cal F_{l_n}$, $C_{n+1}\subset \Cal C_{l_n}$ for some increasing sequence $l_1<l_2<\cdots$ and the conditions (I)--(III), \thetag{1-2}
and~\thetag{1-3} are satisfied.
Then the associated $(C,F)$-action of $G$ belongs to $\Cal O_G$.
It is free because every $(C,F)$-action is free.
\qed
\enddemo

Since $\Psi$ takes values in $\Cal O_G$, we deduce from  Lemmata~3.3 and  3.4 and Proposition~3.6 the following corollary.

\proclaim{Corollary 3.7} $(\goth R_1^{\text{fin}},\Psi)$ is a model for $\Cal A_G$.
\endproclaim

\comment

If $T$ is rank-one map, there is a $(C,F)$-data $(F_n,C_{n+1})_{n=0}^\infty$ such that $T$ is isomorphic to the $(C,F)$-map generated by this data.
Notice that $F_n=\{\alpha_n,\dots,0,\dots,\beta_n\}$ and $C_{n+1}=\{0,c^{(n+1)}_1,\dots, c^{(n+1)}_{r_n-1}\}$.

Given $A,B\subset \Bbb Z$, let 
$$
C^A_B:=\{D\subset\Bbb Z\mid \# D\ge 2, D+B\subset A, (d+B)\cap (d'+B)=\emptyset\text{ if }D\ni d\ne d'\in D\}.
$$
Then $C_{n+1}\in C_{F_n}^{F_{n+1}}$.
We assume that $\# F_{n+1}\ge 2\# F_{n}$.

We say that $(F_n,C_{n+1})_{n=1}^\infty$ is isomorphic to $(F_n',C_{n+1}')_{n=1}^\infty$ if the corresponding $(C,F)$-actions are topologically conjugate.
Let $\goth R^*$ be the subset of $\goth R$ consisting of systems with finite invariant Radon measure.
Then
$$
\goth R^*=\{\Cal S=(F_n,C_{n+1})_{n=1}^\infty\in\goth R\mid \lim_{n\to\infty}\frac{\# F_n}{\#C_1\cdots \#C_n}<\infty\}.
$$

We let $\Psi(\Cal S):=\lim_{n\to\infty}\frac{\# F_n}{\#C_1\cdots \#C_n}\in\Bbb R\cup\{\infty\}$.
It is a Borel map.
Hence there is a weakest Polish topology on $\goth R$ finer than the original topology such that $\Psi$ is continuous.
Denote it by $\tau^*$.
Of course, $\Cal R^*$ is an $\tau^*$-open and hence the space $(\goth R^*,\tau^*)$ is Polish.
Hence $\Cal S_m\to\Cal S$ in $\tau^*$ if and only if  $\Cal S_m\to\Cal S$ in $\tau$ and
$\Psi(\Cal S_m)\to\Psi(\Cal S)$.

\endcomment

\head 4. Measure preserving rank-one actions with bounded parameters
\endhead

In this section $G$ is Abelian.
We prove here 
Theorems~G and H.

\subhead 4.1. Rigidity and rigidity along asymptotically invariant subsequence of subsets
\endsubhead
Let  $S=(S_g)_{g\in G}$ stand for an ergodic free action of $G$ on a  standard probability measure space $(Y,\nu)$.
$S$ is called {\it rigid}  if there
is a sequence $(g_n)_{n=1}^\infty$ of elements from $G$ such that $g_n\to\infty$ and $S_{g_n}\to\text{Id}$ as $n\to\infty$ in the Polish group Aut$(Y,\nu)$.
The sequence $(g_n)_{n=1}^\infty$ is called a {\it rigidity sequence for} $S$.
It follows straightforwardly  that 
\roster
\item"$\bullet$"
 for each $0\ne p\in\Bbb Z$, the sequence $(pg_n)_{n=1}^\infty$ is also a rigidity sequence for $S$;
 \item"$\bullet$" 
 if $(g_n')_{n=1}^\infty$ is  another rigidity sequence in $G$  then either  $g_n=g_n'$ eventually or the difference $(g_n-g_n')_{n=1}^\infty$ is  a rigidity sequence for $S$;
  \item"$\bullet$" a subsequence of the rigidity sequence for $S$ is a rigidity sequence for $S$.
\endroster

We need a generalization of the rigidity concept. 
Let $(W_n)_{n=1}^\infty$ be an {\it asymptotically $S$-invariant} sequence of subsets in $Y$, i.e.
$\nu(S_gW_n\triangle W_n)\to 0$ for each $g\in G$ and $\mu(W_n)\to\delta>0$ as $n\to\infty$.
A standard reasoning implies that then   $\nu(A\cap W_n)\to\delta\nu(A)$ for each subset $A\subset Y$ (see, for instance, \cite{Da6, Lemma~5.6} for the case where $G=\Bbb Z$; in the general case the proof is similar).

The following definition appeared in \cite{BeFr} in the case $G=\Bbb Z$.

\definition{Definition 4.1 }
If there is a sequence $(g_n)_{n=1}^\infty$ of elements from $G$ such that for each subset $A\subset Y$,
$$
\nu((S_{g_n}A\triangle A)\cap W_n)\to 0\quad\text{as }n\to\infty\tag4-1
$$
then $S$ is called {\it rigid} along $(W_n)_{n=1}^\infty$ and  $(g_n)_{n=1}^\infty$ is called {\it a rigidity sequence for $S$ along $(W_n)_{n=1}^\infty$}.
\enddefinition

We note that \thetag{4-1} is equivalent to
$$
\nu(S_{g_n}A\cap A\cap W_n)\to \delta\nu(A)\quad\text{as }n\to\infty\tag4-2
$$
Of course, each rigidity sequence for $S$ is a rigidity sequence for $S$ along every asymptotically $S$-invariant sequence.

\proclaim{Lemma 4.2} Let  $(g_n)_{n=1}^\infty$ be  a rigidity sequence for $S$ along $(W_n)_{n=1}^\infty$.
If  $(g+g_n)_{n=1}^\infty$ is also a rigidity sequence for $S$ along $(W_n)_{n=1}^\infty$ for some $g\in G$ then $g=0$.
\endproclaim

\demo{Proof} For each subset $A\subset Y$, 
$$
\align
\lim_{n\to\infty}\nu((S_gS_{g_n}A\triangle A)\cap W_n)&=\lim_{n\to\infty}\nu((S_{g_n}A\triangle S_{-g}A)\cap S_{-g}W_n)\\
&=\lim_{n\to\infty}\nu((S_{g_n}A\triangle S_{-g}A)\cap W_n).
\endalign
$$
It follows from this, \thetag{4-1} and the condition of the lemma that 
$$
\lim_{n\to\infty}\nu((A\triangle S_{-g}A)\cap W_n)=0.
$$
Since $(W_n)_{n=1}^\infty$  is asymptotically $S$-invariant, we obtain that 
$\nu(A\triangle S_{-g}A)=0$.
Hence $S_{-g}$ is the identity and we are done.
\qed
\enddemo

\subhead 4.2. Bounded $(C,F)$-constructions. Rigidity criterion for the bounded rank-one actions of $\Bbb Z$
\endsubhead
Let $T=(T_j)_{j\in\Bbb Z}$ be a $(C,F)$-action of $G$ associated with a sequence $(C_n,F_{n-1})_{n\ge 1}$ satisfying (I)--(III) and \thetag{1-4}.
The following definition generalizes naturally  Ryzhikov's definition  given in the case of  rank-one actions of $\Bbb Z$  \cite{Ry}.

\definition{Definition 4.3} We say that the sequence of parameters $(C_n,F_{n-1})_{n\ge 1}$ is  {\it  bounded} if  there is $R>0$ and
a finite subset $K\subset G$ such that $\# C_n\le R$ and $K+F_n+C_{n+1}\supset F_{n+1}$ for each $n\ge 0$.
\enddefinition 

It is easy to see that each $(C,F)$-action associated with bounded sequence of parameters is {\it finite} measure preserving.
We note also  that if a sequence $(\widetilde C_n,\widetilde F_{n-1})_{n>0}$ is a bounded telescoping of 
 a bounded sequence  $(C_n,F_{n-1})_{n\ge 1}$ then $(\widetilde C_n,\widetilde F_{n-1})_{n>0}$ is also bounded.
In particular, for each
$p>0$, the
$(pn)_{n\ge 0}$-telescoping of  $(C_n,F_{n-1})_{n\ge 1}$  yields a  bounded sequence of parameters.

From now on let $G=\Bbb Z$ and
$
F_n=\{0,1,\dots,h_n-1\}$
  for some  integers $h_n\ge 0$.
  In other words, $T$ is a rank-one action of $\Bbb Z$.
Then 
$(C_{n+1},F_{n+1})\sim(r_n,\sigma_n)$ for certain $r_n$ and $\sigma_n$ via   \thetag{1-5} and \thetag{1-6} for each $n\ge0$.
\comment

$C_{n+1}=\{ih_n+s_n(i)\mid i=0,\dots,r_n-1\}$
for some
 $r_n>1$ and an increasing map $s_n:\{0,\dots,r_n-1\}\to\Bbb Z_+$ with $s_n(0)=0$,
$n=0,1,\dots$.
Let $\sigma_n$ be determined by \thetag{1-5}.

We define a map $\sigma_n:\{1,\dots,r_n\}\to\Bbb Z_+$ by settting
$$
\sigma_n(i):=\cases
s_n(i)-s_n(i-1)& \text{if } i=1,\dots,r_n-1\text{ and}\\
h_{n+1}-r_nh_n-s_n(r_n-1)&\text{if }i=r_n.
\endcases
$$
It is called the {\it spacer map} on the $n$-th step of the cutting-and-stacking construction.
The sequences $(r_n)_{n=1}^\infty$ and $(\sigma_n)_{n=1}^\infty$
determine completely (and in a unique way) the sequence  $(C_n,F_{n-1})_{n\ge 1}$.

\endcomment
We note that $(C_n,F_{n-1})_{n\ge 1}$ is   bounded if and only if 
there is $R>0$ such that $r_n\le R$  and   $\max_{1\le i\le r_n}\sigma_n(i)\le R$
for each $n>0$.
We also note that if $(C_n,F_{n-1})_{n\ge 1}$ is   bounded then the canonical measure $\mu$ is finite, i.e. \thetag{1-2} is satisfied.

We need a notation. 
Given two maps 
$$
a:\{1,\dots,r\}\ni i\mapsto a(i)\in\Bbb Z_+\quad\text{and}\quad  b:\{1,\dots,t\}\ni i\mapsto b(i)\in\Bbb Z_+,
$$
 we define a new map $a\diamond b:\{1,\dots,rt\}\ni i\mapsto a\diamond b\in \Bbb Z_+$ by setting
$$
a\diamond b(i)=\cases
a(j) &\text{if }i\equiv j\pmod r\text{ and }1\le j<r,\\
a(r)+b(k) &\text{if }i=rk\text{ and }1\le k\le t.
\endcases
$$
Our interest to this ``operation'' is due to the following property  of the parameters of  rank-one  $(C,F)$-actions of $\Bbb Z$:
if $(C_{n+1},F_{n+1})\sim(r_n,\sigma_n)$ and $(C_{n+2},F_{n+2})\sim(r_{n+1},\sigma_{n+1})$ then
$(C_{n+1}+C_{n+2},F_{n+2})\sim(r_nr_{n+1},\sigma_n\diamond\sigma_{n+1})$ for each $n\in\Bbb N$. 
We leave verification of this property to the readers because it is a routine.

Given $i\in\{1,\dots,r-2\}$, we say that {\it $a$ is $i$-periodic} if $a(i+j)=a(j)$ fir all $j\in\{1,\dots,r-1-i\}$.

\proclaim{Lemma 4.4} Let $a$ and $b$ be as above.
\roster
\item"$(i)$"
If $b$ is $i$-periodic for some $i\in\{1,\dots,t-2\}$ then
$a\diamond b$ is $ri$-periodic.
\item"$(ii)$" If $a$ is both $i$-periodic and $j$-periodic,  $j>i$ and $i+j<r$ then $a$ is also $(j-i)$-periodic.
 \item"$(iii)$"
 If $a\diamond b$ is $i$-periodic and $i<r$ then $b(1)=b(2)=\cdots=b(t-1)$.
\endroster
\endproclaim

\demo{Proof}
$(i)$ Take $j\in\{1,\dots, rt-1-ri\}$. Then 
$(a\diamond b)(ri+j)$ equals
$$
\cases
a(r)+b(i+s)=a(r)+b(s)&\text{if }j=rs\text{ for some }s>0,\\
a(v) &\text{if }j=rs+v\text{ for some $s\ge 0$ and $1\le v<r.$}
\endcases
$$ 
Thus in the very of the two cases we obtain that $(a\diamond b)(ri+j)=(a\diamond b)(j)$.

$(ii)$ Take  $s\in\{1,\dots, r-1-j+i\}$.
If $i+j<r$ then either $j+s<r$ or $s-i>0$.
In the first case we have that $a(j-i+s)=a(j+s)=a(s)$.
In the second case we have that 
$a(j-i+s)=a(-i+s)=a(s)$.

$(iii)$ If $1\le j<t$ then $a(i)=(a\diamond b)(jr+i)=(a\diamond b)(jr)=a(r)+b(j)$.
Hence $b(1)=b(2)=\cdots=b(t-1)$.
\qed
\enddemo

We now state the main result of the subsection.

\proclaim{Theorem 4.5} Let  $T$ be a $(C,F)$-action of $\Bbb Z$ associated with a bounded sequence $(C_n,F_{n-1})_{n\ge 1}$ and let $
F_n=\{0,1,\dots,h_n-1\}$
  for some  integers $h_n\ge 1$ for each  $n\ge 0$.
Then $T$ is rigid if and only if for each $N>0$, there are integers $n,m$ such that $m>n+N>n>N$  and  the set
$C_n+\cdots+C_m$ is an arithmetic sequence.
\endproclaim

\demo{Proof}
The ``if" part is easy.
Indeed, given $N>0$ we can find integers $n,m,l$ and $a_N$  such that
$m>n+N>n>N$ and $C_n+\cdots+C_m=\{0,a_N,2a_N,\dots, la_N\}$.
We set $P:=\{0,a_N,2a_N,\dots, la_N\}$,  $Q:=\{0,a_N,2a_N,\dots, (l-1)a_N\}$ and $R:=\{a_N,2a_N,\dots,l a_N\}=a_N+Q$.
Given two subsets $A,B\subset F_N$, we now have that
$$
\align
\mu(T_{a_N}[A]_N\cap[B]_N)&=\mu\left(T_{a_N}\biggl[A+\sum_{i=N+1}^{n-1} C_i+P\bigg]_m\cap
[B]_N\right)\\
&=\mu\left(T_{a_N}\biggl[A+\sum_{i=N+1}^{n-1} C_i+Q\bigg]_m\cap
[B]_N\right)\pm\frac 1{2^N}\\
&=\mu\left(\biggl[A+\sum_{i=N+1}^{n-1} C_i+R\bigg]_m\cap
[B]_N\right)\pm\frac 1{2^N}\\
&=\mu([A]_N\cap[B]_N)\pm\frac 2{2^N}.
\endalign
$$
Since the set of cylinders is dense in the entire Borel $\sigma$-algebra on $X$,
it follows that $T$ is rigid along $(a_N)_{n=1}^\infty$.

As for the ``only if", we first note that it is enough to prove  the following claim: 
\roster
\item"$\bullet$" for each $N>0$, there is $n>N$ such that  $C_n$ is an arithmetic sequence.
\endroster
Indeed, if this claim is true then for each $k>0$, the action $T$ is associated also with a   $(kn)_{n=1}^\infty$-telescoping of $(C_n,F_{n-1})_{n\ge 1}$.
The telescoping  yields a bounded  sequence of parameters.
It remains to apply the claim to this sequence of parameters.

We now prove the claim.
Assume that $T$ is rigid.
Let $(F_{n+1},C_{n+1})\sim(r_n,\sigma_n)$ for each $n\ge0$.
 Since
 $(C_n,F_{n-1})_{n\ge 1}$ is bounded, there is $R>0$  such that $\sup_{n\ge 1}r_n\le R$  and   $\sup_{n\ge 1}\max_{1\le i\le r_n}\sigma_n(i)\le R$.
Passing to the  $(2n)_{n=1}^\infty$-telescoping we may assume without loss of generality that  $\# C_n\ge 4$
for each $n>0$.
Choose  a rigidity sequence  $(g_n)_{n=1}^\infty$ for $T$ such that $g_n\to+\infty$ as $n\to\infty$.
Then for each $n>0$,
 there is a unique $l_n>0$ such that $g_n\in F_{l_{n}+1}\setminus F_{l_n}$.
 In view of the boundedness condition we may assume without loss of generality that
 there exists  $f_n\in F_{l_n}+\{0,1,\dots, R^2\}$ and $i_n\in\{1,\dots,r_{l_n}-1\}$ such that
 $
 g_n=f_n+i_nh_{l_n}$.
\comment

 If  $i_n>0.6r_{l_n}$ for some $n$, then 
 $$
 2g_n=2f_n+2i_nh_{l_n}\in F_{l_n+2}\setminus F_{l_n+1}.
 $$
 Hence passing to a subsequence of  $(g_n)_{n=1}^\infty$ (or replacing  $(g_n)_{n=1}^\infty$ with   $(2g_n)_{n=1}^\infty$  and passing to a subsequence of  $(2g_n)_{n=1}^\infty$) we may assume without loss in generality that $i_n\le0.6r_{l_n}$ for all $n>0$.

 \endcomment
 Using again the boundedness of $(C_n,F_{n-1})_{n\ge 1}$   and passing (if necessary) to a further subsequence of  $(g_n)_{n=1}^\infty$ we may assume in addition  that
 there are integers $r,r'\le R$, $i\in\{1,\dots,r-1\}$, a real $\delta\in[0,1]$ and  maps $\sigma:\{1,\dots,r\}\to\{0,1,\dots,R\}$ and
 $\sigma':\{1,\dots,r'\}\to\{0,1,\dots,R\}$
 such that  $r_{l_n}=r$, $r_{l_n+1}=r'$, $i_n=i$,
 $\sigma_{l_n}=\sigma$ and  $\sigma_{l_n+1}=\sigma'$ for all $n>0$ and $\lim_{n\to\infty}f_n/h_{l_n}=\delta$.
 Let $s$ stand for the integral of $\sigma$.
 We consider separately three cases.
 
 {\it Case A.} Let $\delta<1$.
\comment
 We then set 
$$
s(j):=
\cases
\sum_{k=1}^j\sigma(k)&\text{if }j=1,\dots, r-1,\\
0&\text{if } j=0.
\endcases
$$
\endcomment
Take a cylinder  $A$ in $X$.
Then for each sufficiently large $n$, there is a subset $A_n\subset F_n$ such that $A=[A_n]_n$.
We now set 
$$
\align
W_n&:=[F_{l_n}\cap (F_{l_n}-f_n-2R^2)]_{l_n}\quad\text{and}\\
V_n:&=\left[\bigsqcup_{j=0}^{r-1-i}(F_{l_n}\cap (F_{l_n}-f_n-2R^2)+ jh_{l_n}+s(j)\right]_{l_n+1}.
\endalign
$$
It is easy to verify that the two sequences $(W_n)_{n=1}^\infty$ and $(V_n)_{n=1}^\infty$ are asymptotically  $T$-invariant, 
$$
\lim_{n\to\infty}\mu(W_n)=1-\delta>0\quad\text{and}\quad
\lim_{n\to\infty}\mu(V_n)=\frac{(1-\delta)(r-i)}{r}>0.\tag4-3
$$
\comment
$$
\aligned
A_{l_n}^\circ&:=A_{l_n}\cap (F_{l_n}-f_n-2R^2)\quad\text{and}\\
A_{l_n+1}'&:=\bigsqcup_{j=0}^{r-1-i}(A_{l_n}^\circ+ jh_{l_n}+s(j)]_{l_n+1})\subset A_{l_n+1}.
\endaligned
\tag4-3
$$
Then we have 
$$
\frac{\# A_{l_n+1}'}{\# A_{l_n+1}}=\frac{\#(A_{l_n}\cap (F_{l_n}-f_n-R^2))}{\# A_{l_n}}\cdot\frac{r-i}r\to\frac{(1-\delta)(r-i)}r
\tag4-4
$$
as $n\to\infty$.
Since $(g_n)_{n=1}^\infty$
 is a rigidity sequence for $T$,
 \thetag{4-4} yields that
 $$
 \lim_{n\to\infty}\mu([A'_{l_n+1}]_{l_n+1}\cap A)=
\frac{(1-\delta)(r-i)}r
\mu(A).\tag4-5
 $$
 \endcomment
 We set $A_{l_n}^\circ:=A_{l_n}\cap (F_{l_n}-f_n-2R^2)$.
Since $g_n=f_n+ih_{l_n}$, it follows  that
$$
T_{g_n}(A\cap V_n)=\bigsqcup_{j=i}^{r-1}[(A_{l_{n}}^\circ+f_n+s(j-i)-s(j))+jh_{l_n}+s(j)]_{l_n+1},
$$
because  $A_{l_{n}}^\circ+f_n+s(j-i)-s(j)\subset F_{l_n}$ for each $j=i,\dots,r-1$.
\comment
$$
T_{g_n}[A_{l_n+1}']_{l_{n}+1}=\bigsqcup_{j=i}^{r-1}[(A_{l_{n}}^\circ+f_n+s(j-i)-s(j))+jh_{l_n}+s(j)]_{l_n+1}
$$
and $A_{l_{n}}^\circ+f_n+s(j-i)-s(j)\in F_{l_n}$ for each $j=i,\dots,r-1$.
\endcomment
Therefore
$$
\align
\mu(T_{g_n}(A\cap V_n)\cap A)&=\frac 1r\sum_{j=i}^{r-1}\mu([(s(j-i)-s(j)+f_n+A_{l_{n}}^\circ)\cap A_{l_n}]_{l_n})\\
&=\frac 1r\sum_{j=i}^{r-1}\mu([(s(j-i)-s(j)+f_n+A_{l_{n}}^\circ)]_{l_n}\cap [A_{l_n}]_{l_n})\\
&=
\frac 1r\sum_{j=i}^{r-1}\mu(T_{s(j-i)-s(j)+f_n}[A_{l_{n}}^\circ]_{l_n}\cap A)\\
&=
\frac 1r\sum_{j=i}^{r-1}\mu(T_{s(j-i)-s(j)+f_n}(A\cap W_n)\cap A).
\endalign
$$
Passing to the limit and utilizing \thetag{4-3}  and the rigidity of $(g_n)_{n\ge 1}$ we obtain that
$$
\align
\frac{(1-\delta)(r-i)\mu(A)}{r}&=\lim_{n\to\infty}\mu((A\cap V_n)\cap T_{-g_n}A)\\
&=
\lim_{n\to\infty}\frac 1r\sum_{j=i}^{r-1}\mu(A\cap W_n\cap T_{s(j)-s(j-i)-f_n}A)\\
&\le\frac{(1-\delta)(r-i)\mu(A)}{r}.
\endalign
$$
The equality is only possible if  for each $j=i,\dots,r-1$, there exists a limit
$$
\lim_{n\to\infty}\mu(A\cap W_n\cap T_{s(j)-s(j-i)-f_n}A)=(1-\delta)\mu(A).
$$
Applying the standard approximation argument we conclude that the limit  exists for an arbitrary Borel set $A\subset X$, not only for the cylinders.
Therefore, in view of  \thetag{4-2}, the sequence $(s(j)-s(j-i)-f_n)_{n=1}^\infty$
is a rigidity sequence for $T$ along $(W_n)_{n=1}^\infty$. 
Hence Lemma~4.2
yields that 
$$
s(i)=s(i+1)-s(1)=s(i+2)-s(2)=\cdots=s(r-1)-s(r-1-i).
$$
This is equivalent to the fact that $\sigma$ is $i$-periodic.
On the other hand,  we see that $g_n\in F_{l_n+2}\setminus F_{l_n}$.
Since $(C_{l_n+1},F_{l_n+1})\sim(r,\sigma)$ and 
$(C_{l_n+2},F_{l_n+2})\sim(r',\sigma')$, it  follows that $(C_{l_n+1}+C_{l_n+2},F_{l_n+2})\sim(rr',\sigma\diamond\sigma')$ for each $n\ge 0$.
Hence replacing  $C_{l_n+1}$ in the above argument with $C_{l_n+1}+C_{l_n+2}$ 
(i.e. passing to the corresponding bounded telescoping), we obtain that 
 the map $\sigma\diamond\sigma'$ is also $i$-periodic.
\comment

By the definition of $\sigma\diamond\sigma'$,
$$
\sigma\diamond\sigma'(r-i)=\sigma\diamond\sigma'(2r-i)=\cdots=
\sigma\diamond\sigma'(r(r'-1)-i)=\sigma(r-i).
$$
Since $\sigma\diamond\sigma'$ is  $i$-periodic, it follows that
$$
\sigma\diamond\sigma'(r)=\sigma\diamond\sigma'(2r)=\cdots=
\sigma\diamond\sigma'(r(r'-1))
$$
and hence, by the definition of $\sigma\diamond\sigma'$,

\endcomment
By Lemma~4.4(iii),
$
\sigma'(1)=\sigma'(2)=\cdots=\sigma'(r'-1).
$
This is equivalent to the fact that  $C_{l_n+2}$ is an arithmetic sequence for each $n\ge 0$.

{\it Case B.} Let  $\delta=1$ but $i<r-1$.
Then we  write $g_n=(f_n-h_n)+(i+1)h_{l_n}$ for each $n>0$.
Since $i+1<r$ and $\lim_{n\to\infty}\frac{f_n-h_{l_n}}{h_{l_n}}=0<1$, we prove the claim in~Case~B in the very same way as in Case A just by replacing
 $f_n$ with $f_n-h_n$ and $i$ with $i+1$ everywhere in the above argument.
 
 {\it Case C.} Let  $\delta=1$ and $i=r-1$.
 Then $\lim_{n\to\infty}g_n/h_{l_n+1}=1$.
 Hence $2g_n\in F_{l_n+2}\setminus F_{l_n+1}$ for each $n>0$.
 The sequence $(2g_n)_{n=1}^\infty$ is also a rigidity sequence for  $T$.
 Since $\# C_{l_n+2}\ge 4$, it follows that $2g_n/h_{l_n+2}\le 0.5$  for each $n>0$.
 Thus considering $(2g_n)_{n=1}^\infty$ instead of $(g_n)_{n=1}^\infty$ we reduce Case~C to Case~A or Case~B.
\qed

\enddemo

As a corollary we obtain a  characterization of the discrete spectrum of rigid rank-one $\Bbb Z$-actions with bounded parameters.
To state this assertion we need  a new concept.
First, we recall that each odometer $Q=(Q_n)_{n\in\Bbb Z}$ can be represented as a minimal rotation on a monothetic compact totally disconnected Abelian group $H$ which is the inverse limit of a sequence of cyclic groups:  
$
\Bbb Z/l_1\Bbb Z@<<< \Bbb Z/(l_1l_2\Bbb Z) @<<< \cdots
$
associated with a sequence of  integers $l_i>1$ for each $i\in\Bbb N$.
An element of $q\in Q$ is a sequence $q=(q_i)_{i\in\Bbb N}$ such that $0\le q_i<l_1\cdots l_i$
and $q_i\equiv q_{i+1}\pmod{ l_1\cdots l_{i}}$ for each $i\in\Bbb N$.
Then $(Q_nq)_i=q_i+n\pmod{ l_1\cdots l_{i}}$ for each $q\in Q$ and $n\in\Bbb Z$.
We say that $Q$ {\it is of bounded type} if the set $\{l_i\mid i\in\Bbb N\}$ is finite.
 It is easy to verify that the  bounded type is well defined by $Q$, i.e. it does not depend on the choice  of the 
 ``approximating'' sequence $(l_i)_{i\in\Bbb N}$\footnote{In fact, it is a spectral property of $Q$, i.e. a property of  the Koopman operator generated by $Q$.}.

\proclaim{Proposition 4.6} Let  $T$ be a rigid $(C,F)$-action of $\Bbb Z$ associated with a bounded sequence $(C_n,F_{n-1})_{n\ge 1}$ and let $
F_n=\{0,1,\dots,h_n-1\}$
  for some  integers $h_n\ge 1$ for each  $n\ge 0$.
Then either $T$ is an odometer of bounded type  or the group $\Lambda_T\subset\Bbb T$ of eigenvalues of $T$ is finite.
Moreover, in the latter case, if $(C_{n+1},F_{n+1})\sim(r_n,\sigma_n)$ for each $n\ge 0$ and 
$K:=\sup_{n\ge 1}\sup_{1\le j\le r_n}\sigma_n(j)$ then the order of each $\lambda\in\Lambda_T$ does not exceed  $K$.
\endproclaim

\demo{Proof}
Suppose first that there are infinitely many $n\in\Bbb N$ for which there exists $m>2n$ such that the set $C_n+\cdots +C_m$ is an arithmetic sequence but $C_n+\cdots +C_m+C_{m+1}$ is not.
We will call such $n$ {\it good}.
Let $\lambda\in\Lambda_T$  and  let $\xi:X\to\Bbb C$ be a corresponding non-zero measurable eigenfunction for $T$, i.e.
$\xi\circ T_j=\lambda^j\xi$ almost everywhere on $X$ for each $j\in\Bbb Z$.
Since $T$ is ergodic, we may assume that $|\xi|=1$ almost everywhere on $X$.
Fix $\epsilon>0$.
Then we can  select a good $n>0$, an element $f\in F_{n-1}$ and a complex number $z\in\Bbb T\subset\Bbb C$ such that the set  $A:=\{x\in X\mid |\xi(x)-z|<\epsilon\}$ is of positive measure  and  the cylinder $[f]_{n-1}$ is $(1-\epsilon^2)$-{\it full of $A$}, i.e. 
$$
\mu(A\cap[f]_{n-1})>(1-\epsilon^2)\mu([f]_{n-1}).
$$
By the definition of good numbers,  there are  integers $a>0$, $l\in\{0,\dots,\# C_{m+1}-2\}$ and $v\ne 0$ such that 
$$
\align
C_{n}+\cdots+C_m&=\{0,a,\dots,(p-1)a\}\quad\text{and }\\
C_{m+1}&=\{0, pa, 2pa,\dots,lpa,(l+1)pa+v,\dots\},
\endalign
$$
where $p:=\# C_n\cdots\# C_m$.
Since $pa=h_m+\sigma_m(1)$ and $pa+v=h_m+\sigma_m(l+1)$, it follows that $v=\sigma_m(l+1)-\sigma_m(1)$ and hence $|v|\le K$.
Let $p'$ stand for the integer part of $p/3$.
Since $[f]_{n-1}=\bigsqcup_{c\in C_n+\cdots+C_{m+1}}[f+c]_{m+1}$,
it follows that more then $(1-\epsilon)\#C_n\cdots\#C_{m+1}$ cylinders $[f+c]_{m+1}$, $c\in C_n+\cdots+C_{m+1}$,
are $(1-\epsilon)$-full of $A$.
 Hence there is $1<j<p'$ such that the cylinders
$[f+ja]_{m+1}$, $[f+ja+p'a]_{m+1}$, $[f+ja+(p-1+lp)a-p'a]_{m+1}$ and
$[f+(j-1)a+(l+1)pa+v]_{m+1}$
 are $(1-\epsilon)$-full of $A$.
 Since $[f+ja+p'a]_{m+1}=T_{p'a}[f+ja]_{m+1}$ and 
 $$
 [f+(j-1)a+(l+1)pa+v]_{m+1}=T_{p'a+v}[f+ja+(p-1+lp)a-p'a]_{m+1},
 $$
 It follows that there are $x,y\in A$ such that $T_{p'a}x\in A$ and $T_{p'a+v}y\in A$.
 Hence $|z-\lambda^{p'a}z|<2\epsilon
 $
  and
 $|z-\lambda^{p'a+v}z|<2\epsilon$.
 This implies that $|1-\lambda^v|<4\epsilon$.
 Since $\epsilon$ can be chosen arbitrarily small, $\lambda^v=1$.
 
 Suppose now that the set of good numbers is finite.
 Then it follows from Theorem~4.5 that there is $n>0$ such that the  set
 $C_n+C_{n+1}+\cdots$ is an infinite arithmetic sequence.
 Let $a:=\sigma_{n-1}(1)$.
 Then, as it is well known, $T_1$ is the integral transformation constructed over the odometer determined by the following approximating sequence of integers: $\# F_{n-1}$, $\# F_{n-1}\#C_n$, $\# F_{n-1}\#C_n\#C_{n+1}$,\dots and under the map that takes constant value $a$. 
 Hence $T_1$ is also an odometer.
 It is determined by the  sequence of integers $a\# F_{n-1}$, $a\# F_{n-1}\#C_n$, $a\# F_{n-1}\#C_n\#C_{n+1}$,\dots.
 Thus we see that  $T_1$ is an odometer of bounded type.
 \qed
\enddemo

\comment

Next, we see that
$$
\bigsqcup_{j=0}^{r-1-i}T_{g_n}[A_{l_{n}}+jh_{l_n}+s(j)]_{l_n+1}=\bigsqcup_{j=0}^{r-1-i}
[A_{l_{n}}+f_n+s(j+i)+(j+i)h_{l_n}-\sum_{k=j+1}^{j+i}\sigma(k)]_{l_n+1}
$$
In a similar way,
$$
\align
\bigsqcup_{i=r-k}^{r-1}T_{g_n}[A_{l_{n}-1}+ih_{l_n-1}+\sum_{j=0}^{i}s(j)]_{l_n}
&=
\bigsqcup_{i=r-k}^{r-1}[g_n-...+A_{l_{n}-1}+ih_{l_n-1}+\sum_{j=0}^{i}s(j)]_{l_n}\\
&=?
\bigsqcup_{i=0}^{r-1}
[A_{l_{n}-1}-\sum_{j=i}^{i+k}s(j)+(i+k)h_{l_n-1}+\sum_{j=0}^{i+k}s(j)]_{l_n}
\endalign
$$
Hence
$\mu(T_{g_n}A\cap A)=\sum_{i=0}^r\mu(T_{g_n}[A_{l_n-1}+c_{l_n}(i)]_{l_n}\cap A)\to \mu(A)$ if and only if  $\mu(T_{g_n}[A_{l_n-1}+c_{l_n}(i)]_{l_n}\cap A)\to \mu(A)/r$ for each $i$.
If $i<r-1-k$, we have
$$
\align
\mu(T_{g_n}[A_{l_n-1}+c_{l_n}(i)]_{l_n}\cap A)
&\approx\mu([(A_{l_n-1}-\sum_{j=i}^{i+k}s(j))\cap A_{l_n-1}+c_{l_n}(i+k)]_{l_n})\\
&=\frac 1r\mu(T^{-\sum_{j=i}^{i+k}s(j)}[A_{l_n-1}]_{l_n}\cap[A_{l_n-1}]_{l_n}).
\endalign
$$
We obtain that $\mu(T^{-\sum_{j=i}^{i+k}s(j)}[A_{l_n-1}]_{l_n}\cap[A_{l_n-1}]_{l_n})\to\mu([A_{l_n-1}]_{l_n})$.
Since $A$ was arbitrary, $\mu(T^{-\sum_{j=i}^{i+k}s(j)}B\cap B)\to\mu(B)$ for each $B$.
This is only possible if $\sum_{j=i}^{i+k}s(j)=0$.
Hence $s(i)=s(i+1)=\dots=s(i+k)=0$.
\qed

\endcomment

\subhead 4.3. Ergodicity of powers for rank-one transformations
\endsubhead
The following lemma is a particular case of \cite{Da1, Lemma~2.4} (see also \cite{Da3, Lemma~2.1}).
However for completeness of our argument,  we state  it here with a  proof.

\proclaim{Lemma 4.7} Let  $T=(T_g)_{g\in G}$ be a $(C,F)$-action of $G$ associated with a sequence $(C_n,F_{n-1})_{n\ge 1}$ satisfying (I)--(III) and \thetag{1-4}.
Let $\delta:G\to (0,1)$ be a map.
If for each $n$ and every pair $f,f'\in F_n$, there is a subset $A\subset [f]_n$ and $l\in\Bbb Z$ such that
$\mu(A)\ge\delta(f'-f)\mu([f]_n)$ and $T_{lg}A\subset[f']_n$
 then the transformation $T_g$ is ergodic.
  \endproclaim

 \demo{Proof}
 Let $A_1,A_2\subset X$ be two subsets of positive finite measure.
 Then there is $n>0$ and elements 
 $f_1,f_2\in F_n$ such that the cylinder $[f_i]_n$ is 0.99-full of $A_i$, i.e. $\mu(A_i\cap[f_i]_n)>0.99\mu([f_i]_n)$, $i=1,2$.
 For $m>n$, consider the partition of $[f_i]_n$ into $m$-cylinders: $[f_i]_n=\bigsqcup_{c\in C}[f_i+c]_m$, where $C:=C_{n+1}+\cdots+C_m$. 
 When $m$ increases, these partitions refine and generate the entire Borel $\sigma$-algebra on $[f_i]_n$.
 Fix $\epsilon>0$.
 By  a standard measure theoretical fact, the total measure of those $m$-cylinders in $[f_i]_n$ that are $(1-\epsilon)$-full of $A_i$ goes to $\mu(A_i\cap[f_i]_n)$ as $m\to\infty$.
Since  $\mu(A_i\cap[f_i]_n)>0.5\mu([f_i]_n)$, there are subsets $D_1,D_2\subset C$ such that $\# D_i>0.5\# C$ and $\mu(A_i\cap[f_i+c]_m)>(1-\epsilon)\mu([f_i+c]_m)$ for all $c\in D_i$.  
Now let $\epsilon:=0.5\delta(f_1-f_2)$ and $m$ be large.
Take $d\in D_1\cap D_2$ (the intersection is not empty) and apply the condition of the lemma to $f_1+d$ and $f_2+d$, which belong both to $F_m$.
Then there are a subset $A\subset [f_1+d]_m$ and $l\in\Bbb Z$ such that $T_{lg}A\subset [f_2+d]_m$ and $\mu(A)>\delta(f_1-f_2)\mu([f_1+d]_m)$.
By the choice of $\epsilon$, it is easy to deduce that
$\mu(A\cap A_1)>0$ and $\mu(T_{lg}(A\cap A_1)\cap A_2)>0$.
Hence $T_g$ is ergodic.
\qed
 \enddemo

We now state and prove the main result of this subsection.

\proclaim{Theorem 4.8} Let $T$ be as in Lemma~4.7, $G=\Bbb Z$ and
$
F_n=\{0,1,\dots,h_n-1\}$
  for some  integers $h_n\ge 1$ for each  $n\ge 0$.
\roster 
\item"(i)"
If $T_d$ is  ergodic
then for each divisor $p$ of $d$, there are infinitely many $n>0$ such that some $c\in C_n$ is not divisible by $p$.
\item"(ii)"
If the sequence $(\# C_n)_{n=1}^\infty$ is bounded  
and for each divisor $p$ of a positive integer $d$,  there are infinitely many  $n$ such that $p$ does not divide some $c\in C_n$ then $T_d$ is ergodic.
\endroster
\endproclaim

\demo{Proof}
(i) 
Suppose that there are a divisor $p$ of $d$ and a positive integer $N$ such that for each $n>N$ and each $c\in C_n$, $c$ is divisible by $p$. 
Since $T_d$ is ergodic,
there is $m>0$ such that $T_{dm}[0]_N\cap[1]_N\ne\emptyset$.
Hence $dm-1\in\sum_{n>N}(C_n-C_n)$.
It follows that  $dm-1$ is divisible by $p$, a contradiction.

(ii) Let  $K$ be a positive integer such that $\sup_{n>0}\# C_n<K$.
We let 
$$
D:=\{j\in\{1,\dots,d-1\}\mid \exists q\in C_n-C_n, q\equiv j\bmod d\text{ for  infinitely many }n  \}.
$$
Identifying naturally the set $\{0,1,\dots,d-1\}$ with the cyclic group $\Bbb Z/d\Bbb Z$, we consider $D$ as a subset of $\Bbb Z/d\Bbb Z$.
It follows easily from the condition of the theorem that $D$ is non-empty.
We claim that $D$ generates   $\Bbb Z/d\Bbb Z$.
Indeed, otherwise there is  $p>1$ such that $p$ divides $d$ and $D\subset\{0,p,2p,\dots,(p_1-1)p \}\subset\Bbb Z/d\Bbb Z$, where $p_1=d/p$.
Moreover,  it follows now from the definition of $D$ that there is $N>0$ such that $p$ divides each element $q\in C_n-C_n$ for every $n>N$.
Since $C_n\subset C_n-C_n$, we obtain  that $p$ divides every element $c\in C_n$ whenever $n>N$. 
This contradicts the condition of the theorem.

It follows from the claim that
 there are a subset $D_0\subset D$ and non-zero integers $\alpha_j$, $j\in D_0$,
such that 
$$
\sum_{j\in D_0}\alpha_jj\in 1+ d\Bbb Z.
\tag4-4
$$
Now we can find mutually disjoint infinite subsets $\Cal L_j\subset\Bbb N$, $j\in D_0$, such that for each $l\in\Cal L_j$, there are $c_l,c_l'\in C_l$ such that $c_l-c'_l-\frac{\alpha_j}{|\alpha_j|} j\in d\Bbb Z$.
\comment
We recall that for $\alpha\in\Bbb R$,
$$
\text{sign}(\alpha):=
\cases \alpha/|\alpha| &\text{if }\alpha\ne 0\\
0&\text{if }\alpha\ne 0
\endcases
$$
\endcomment

Suppose that we are given $N>0$ and $f, f'\in F_N$ such that
$m:=f-f'>0$. 
It follows from \thetag{4-4} that
 $\sum_{j\in D_0}m\alpha_jj\in m+ d\Bbb Z$.
 Choose finite subsets $L_j\subset\Cal L_j$ such that $N<\min L_j$ and $\# L_j=m|\alpha_j|$ for each $j\in D_0$.
 We now set
 $$
 L:=\max_{j\in D_0}  \max L_j\quad\text{and}\quad A:=\left[f+\sum_{j\in D_0}\sum_{l\in L_j}c_l\right]_{L}.
 $$
 Then
 $$
 \sum_{j\in D_0}\sum_{l\in L_j}(c_l'-c_l)\in
  \sum_{j\in D_0}\sum_{l\in L_j}\frac{\alpha_jj}{|\alpha_j|}+d\Bbb Z=
   \sum_{j\in D_0}{m\alpha_jj}+d\Bbb Z=m+d\Bbb Z.
 $$
 Hence there is $l\in\Bbb Z$ such that
 $
dl+m+ \sum_{j\in D_0}\sum_{l\in L_j}c_l=\sum_{j\in D_0}\sum_{l\in L_j}c_l'.
 $
This yields
$$
T_{dl}A=\left[f-m+\sum_{j\in D_0}\sum_{l\in L_j}c_j'\right]_{L}\subset[f']_N.
$$
Moreover,
$$
\mu(A)=\frac{\mu([f]_N)}{\prod_{j\in D_0}\prod_{l\in L_j}\#C_l}\ge \frac{\mu([f]_N)}{\prod_{j\in D_0}K^{m|\alpha_j|}}= \frac{\mu([f]_N)}{\big(K^{\sum_{j\in D_0}|\alpha_j|}\big)^{f-f'}}.
$$
By Lemma~4.7,  $T_d$ is ergodic.
\qed
\enddemo

\proclaim{Corollary 4.9} Let $T$ be as in Theorem~4.8 and let the sequence $(\# C_n)_{n=1}^\infty$ be bounded.
Then $T$ is totally ergodic if and only if for each $d>1$, there are infinitely many $n>0$ such that
some element $c$ of $ C_n$ is not divisible by $d$.
\endproclaim

We now apply the main results of this section to the Chacon maps.

\example{Example  4.10}
(i) Let  $T$ be the Chacon map with 2 cuts.
Then $T$ is associated with the sequence $(C_n,F_{n-1})_{n=1}^\infty$ such that
$F_n=\{0,\dots,h_n-1\}$, $C_{n+1}=\{0, h_n\}$ and $h_{n}=2h_{n-1}+1=2^{n+1}-1$ for each $n\ge 1$.
Hence $T$ is a non-adapted rank-one transformation with bounded parameters.
Since $C_{n}+C_{n+1}=\{0, h_{n-1}, 2h_{n-1}+1, 3h_{n-1}+1\}$ is not an arithmetic sequence for any $n>0$, it follows from Theorem~4.5 that $T$ is not rigid.
Since $h_n=2h_{n-1}+1$, no divisor of $h_{n-1}$ is a divisor of $h_n$ for any $n\in\Bbb N$.
Then Corollary~4.9 yields that $T$ is totally ergodic.

(ii) Let  $T$ be the Chacon map with 3 cuts.
Then $T$ is associated with the sequence $(C_n,F_{n-1})_{n=1}^\infty$ such that
$F_n=\{0,\dots,h_n-1\}$, $C_{n+1}=\{0, h_n,2h_n+1\}$ and $h_{n}=3h_{n-1}+1=\frac{3^{n+1}-1}2$ for each $n\ge 1$.
Hence $T$ is an adapted rank-one transformation with bounded parameters.
Since $C_n$ is not an arithmetic sequence for any $n>0$, 
it follows from Theorem~4.5 that $T$ is not rigid.
Of course, no divisor of $h_{n-1}$ is a divisor of $h_n$ for any $n\in\Bbb N$.
Hence $T$ is totally ergodic in view of ~Corollary~4.9.
\endexample

\head 5. Disjointness and MSJ for~rank-one~transformations with bounded parameters
\endhead

\subhead 5.1 Joinings and disjointness
\endsubhead
We first recall definitions of joinings and disjointness in the sense of Furstenberg (see \cite{dJRu},  \cite{Ru} and  \cite{Fu} for details).
Let $T=(T_n)_{n\in\Bbb Z}$ and $R=(R_n)_{n=1}^\infty$ be two ergodic $\Bbb Z$-actions on standard probability spaces $(X,\mu)$ and $(Y,\nu)$ respectively.
A {\it joining} of $T$ and $R$ is a $(T_1\times R_1)$-invariant probability measure on $X\times Y$ whose pullbacks on $X$ and $Y$ are $\mu$ and $\nu$ respectively.
The set of all joinings of $T$ and $R$ is denoted by $J(T,R)$.
The subset of ergodic joinings of $T$ and $R$ is denoted by $J^e(T,R)$.
We note that $J^e(T,R)$ is the set of extreme points of the convex set $J(T,R)$.
If $J^e(T,R)=\{\mu\times\nu\}$ then $T$ and $R$ are called {\it disjoint}.
In the case where $R=T$, we reduce the notation $J^e(T,R)$ to $J_2^e(T)$.

The {\it centralizer} $C(T)$ of $T$ is the group of all $\mu$-preserving transformations
of $X$ commuting with $T_1$.
Given $\theta\in C(T)$, we define a measure $\mu_\theta$  on $X\times X$ by setting
$\mu_\theta(A\times B)=\mu(A\cap \theta^{-1}B)$, where $A$ and $B$ are Borel subsets in $X$.
Of course, $\mu_\theta\in J^e_2(T)$ and $\mu_\theta$ is supported on the graph of $\theta$.
It is called a {\it graph-joining} of $T$.
More generally, each measure preserving isomorphism $\theta:X\to Y$  intertwining $T$ with $R$ generates a joining from $J^e(T,R)$ that is supported on the graph of $\theta$.
If $J^e_2(T)=\{\mu_{T_n}\mid n\in\Bbb Z\}\cup\{\mu\times\mu\}$ then $T$ is said to have {\it  minimal self-joinings of order 2 (MSJ$_2)$}.
It follows, in particular, that  $C(T)=\{T_n\mid n\in\Bbb Z\}$ and hence $T$ is not rigid if $T$ has MSJ$_2$.
The property of  MSJ of  higher  orders is defined in a similar way (see \cite{dJRu}).

\subhead 5.2. Joinings of rank-one transformations with bounded parameters
\endsubhead
In this subsection $T=(T_i)_{i\in\Bbb Z}$ and 
 $T'=(T'_i)_{i\in\Bbb Z}$ are the $(C,F)$-actions of $\Bbb Z$
 associated with  sequences of parameters
$(C_n,F_{n-1})_{n\ge 1}$  and
 $(C_n',F_{n-1}')_{n\ge 1}$ respectively.
The two sequences satisfy (I)--(III), \thetag{1-2} and \thetag{1-4}.
We will also assume that $F_n=\{0,\dots,h_n-1\}$ and $F_n'=\{0,\dots,h_n'-1\}$ for some $h_n>0$ and $h_n'>0$ for each $n\ge 0$.
Let $(X,\mu)$ and $(X',\mu')$ be the measure spaces of  $T$ and $T'$ respectively.
Let $(C_{n+1},F_{n+1})\sim(r_n,\sigma_n)$ and $(C_{n+1}',F_{n+1}')\sim(r_n',\sigma_n')$
for some maps $\sigma_n:\{1,\dots r_n\}\to\Bbb Z_+$ and $\sigma_n':\{1,\dots r_n'\}\to\Bbb Z_+$
for each $n>0$.
Denote by $s_n:\{0,\dots,r_n-1\}\to\Bbb Z_+$ and $s_n':\{0,\dots,r_n-1\}\to\Bbb Z_+$
the integrals of $\sigma_n$ and $\sigma_n'$ respectively, $n\ge 0$.

\comment
If for some $n>0$, we have that $F_{n-1}=F_{n-1}'=:h_{n-1}$ and $\# C_n=\# C_n'=:r_n$
then we can write $C_n=\{s_n(i)+ih_{n-1}\mid i=0,1,\dots,r_n-1\}$
and $C_n'=\{s_n'(i)+ih_{n-1}\mid i=0,1,\dots,r_n-1\}$ for some mappings $s_n,s_n':\{0,1,\dots,r_n-1\}\to\Bbb Z_+$ with $s_n(0)=0$.
We extend $s_n$ and $s_n'$ to the point $r_n$ by setting $s_n(r_n):=h_{n}-r_nh_{n-1}-s_n(r_n-1)$ and $s_n'(r_n):=h_{n}-r_nh_{n-1}-s_n'(r_n-1)$.

We write $s_n\perp s_n'$ if for each $i_0=0,\dots,r_n-1$, there is $i<r_n-i_0$ such that
$s'(i+i_0)-s'(i_0)\ne s(i)$.
This means the following. 
Let $s(i)=\sigma(i)+\sigma(i-1)+\cdots+\sigma(0)$, $\sigma(0)=0$.
Then $s'(i+i_0)-s'(i_0)=s(i)$ for all $i<r-i_0$ if and only if $\sigma(i)=\sigma'(i+i_0)$ for all $i<r-i_0$.
May be we should add a symmetric  condition (just reverse $s$ and $s'$).

Given two maps $s,s':\{0,1,\dots,r-1\}\to\Bbb Z_+$, we write $s\sim s'$ if there is $i_0$ such that
$$
s'(i)=
\cases
s(i+i_0) &\text{if }1\le i< r-i_0\\
s(i+i_0-r) &\text{if }r-i_0< i<r
\endcases
$$

\endcomment
\comment
Given $C_n$ and $C_{n+1}$, if $C_n=\{ih_{n-1}+s_n(i)\mid i=0,\dots,r_n-1\}$
and 
$C_{n+1}=\{ih_{n+1}+s_{n+1}(i)\mid i=0,\dots,r_{n+1}-1\}$,
let $\sigma_n:\{1,\dots,r_n\}\to\Bbb Z_+$ is given by $s_n(i)=\sigma_n(1)+\dots+\sigma_n(i)$,
$i=1,\dots,r_n-1$ and $\sigma_n(r_n)=h_{n+1}-r_nh_{n-1}-s_n(r-1)$.
Then the corresponding map $s$ for $C_n+C_{n+1}$ is the following:
$s(ir_{n}+j)=s_{n}(j)$ if $j=1,\dots,r_n-1$ or $s(ir_{n}+j)=s_n(j)+s_{n+1}(i)$ if $j=r_n$, $i=0,\dots,r_{n+1}$.
\endcomment

\comment
\proclaim{Theorem 5.1} Suppose that  $F_n=F_n'$ for all $n$ and there is $r>1$ and an infinite subset $D\subset\Bbb N$ such that $\#C_n=\#C_n'=r$, $s_n=s$, $s_n'=s'$  and $s\perp s'$ for all $n\in D$.
Then $T$ and $S$ are disjoint.
\endproclaim
\endcomment

\proclaim{Proposition 5.1} Suppose that there exist  an infinite subset $D\subset\Bbb N$, an integer $r\ge 3$ and two maps $s,s':\{0,1,\dots,r-1\}\to\Bbb Z_+$ such that
  $h_n=h_n'$,   
 $r_n=r_n'=r$, 
  $s_{n}=s$ and $s_{n}'=s'$  whenever $n+1\in D$.
    Then for each $\lambda\in J^e(T,T')$ and every $i_0\in\{0,1,\dots,r-1\}$,
 there is $i_0'\in\{0,1,\dots,r-1\}$ such that 
 $$
 \lambda=\lambda\circ (T_{s(i+i_0)-s(i_0)}\times T'_{s'(i+i_0')-s'(i_0')})
 $$
  for every $i$ such that $-\min(i_0,i_0')\le i\le r-1-\max(i_0,i_0')$.
\endproclaim

\demo{Proof}
Let $\lambda\in J^e(T,T')$.
A point $(x,x')\in X\times X'$ is called a {\it generic point} for $\lambda$ if 
$$
\lim_{n\to\infty}\frac 1{b_n-a_n}\sum_{i=a_n}^{b_n-1}1_{A\times A'}(T_ix,T'_ix')=\lambda(A\times A')
\tag 5-1
$$
for all cylinders $A\subset X$ and $A'\subset X'$ whenever $b_n>0$, $a_n\le 0$ and $b_n-a_n\to+\infty$ as $n\to\infty$.
By the individual ergodic theorem, $\lambda$-a.e. point of $X\times X'$ is generic.
Fix $m>0$ such that $x=(f_m,c_{m+1},c_{m+2},\dots)\in X_m$ and $x'=(f_m',c_{m+1}',c_{m+2}'\dots)\in X_m'$.
Fix $i_0\in\{0,\dots,r-1\}$.
Since $\lambda$ projects onto $\mu$ and the measure $\mu\restriction X_m$ is proportional to the infinite product of the equidistributions on $F_m$ and on $C_j$ for all $j>m$, 
 we may assume without loss of generality (passing to a subsequence of $D$ if necessary) that  $3h_n/4\ge f_{n}\ge h_{n}/4$ and $c_{n+1}=s(i_0)+i_0h_n$ whenever $n+1\in D$. 
 Moreover, passing to a further subsequence in $D$, we may also assume
 that
 there is $ i_0'\in\{0,\dots,r-1\}$ such that $c_{n+1}'=s'(i_0')+i_0'h_n$ if $n+1\in D$ and 
 there exist 
 $$
 \lim_{D\ni n+1\to\infty}f_n/h_n=\delta\in[ 0.25,0.75]\quad\text{and }
\lim_{D\ni n+1\to\infty}f_n'/h_n=\delta'.\tag5-2
$$
Of course, $0\le\delta'\le 1$.
\comment
 We now let $\widetilde F_m:=F_m-f_m$, $\widetilde C_{m+1}:=C_{m+1}-c_{m+1}$, etc. and 
 $\widetilde F_m':=F_m'-f_m'$, $\widetilde C_{m+1}':=C_{m+1}'-c_{m+1}'$...
 Passing to the new coordinates, we may assume that $x=(0,0,\dots)$ and $x'=(0,0,\dots)$.
 
\endcomment
Let $A$ and $A'$ be two cylinders in $X$ and $X'$ respectively.
We represent  them as $A=[A_n]_{n}$ and $A'=[A_n']_n$ with $A_n\cup A_n'\subset F_n$ for all sufficiently large $n$.
 We will also assume that  $A_n\pm\max_{0<i<r} s(i)\subset F_n$ and $A_n'\pm\max_{0< i<r} s'(i)\subset F_n$ if $n+1\in D$.
 We now let 
 $$
 a_n:=\cases -f_n & \text{if $f_n'\ge f_n$}\\
 -f_n' & \text{if $f_n'\le f_n$}
 \endcases
\quad\text{and}\quad
b_n:=\cases 
h_n-f_n' & \text{if $f_n'\ge f_n$}\\
h_n -f_n& \text{if $f_n'\le f_n$}.
 \endcases
 $$
 Then $a_n\le 0$, $b_n> 0$ and $b_n-a_n\to+\infty$ as $n\to\infty$.
If $i\in[a_n,b_n)$ then the point
 $$
 T_ix=T_i(f_n,c_{n+1},\dots)=(i+f_n,c_{n+1},\dots)
 $$
 belongs to $A$
 if and only if  $i+f_n\in A_n$.
 In a similar way, $T_i'x'\in A'$ if and only if $i+f_n'\in A_n'$.
 The following four cases are possible: $f_n'\ge f_n$ and $f_{n+1}'\ge f_{n+1}$, 
 $f_n'\ge f_n$ and $f_{n+1}'\le f_{n+1}$, $f_n'\le f_n$ and $f_{n+1}'\ge f_{n+1}$, and $f_n'\le f_n$ and $f_{n+1}'\le f_{n+1}$.
 We consider only the first case because the other ones are similar.
 It follows from~\thetag{5-1} that
 $$
 \aligned
 \lambda(A\times A')&=
\frac{\#((A_n-f_n)\cap (A_n'-f_n')\cap[a_n,b_n))}{h_n+f_n-f_n'}+\overline o(1)\\
&= \frac{\#((A_n-f_n)\cap (A_n'-f_n'))}{h_n+f_n-f_n'}+\overline o(1)
\endaligned
\tag5-3
 $$
 eventually in $n$.
 As usual $\overline o(1)$ denotes a sequence going to $0$ as $n\to\infty$.
Take $n+1\in D$.
  Since $f_{n+1}=f_n+c_{n+1}$ and $f_{n+1}'=f_n'+c_{n+1}'$, it follows that
  $$
  A_{n+1}-f_{n+1}=A_n-f_n+C_{n+1}-c_{n+1}
  =\bigsqcup_{i=-i_0}^{r-1-i_0}(A_n-s(i_0)+s(i+i_0)-f_n+ih_n)
  $$
  and, in a similar way,
  $$
  A_{n+1}'-f_{n+1}'= \bigsqcup_{i=-i_0'}^{r-1-i_0'}(A_n'-s'(i_0')+s'(i+i_0')-f_n'+ih_n).
  $$
  We let $\widehat i_0:=\min(i_0,i_0')$,  $\widetilde i_0:=\max(i_0,i_0')$, $B_{n,i}:=A_n-s(i_0)+s(i+i_0)$ and $B_{n,i}':=A_n'-s'(i_0')+s'(i+i_0')$. 
  Then $B_{n,i}\cup B_{n,i}'\subset F_n$ and
 $$
( A_{n+1}-f_{n+1})\cap (A_{n+1}'-f_{n+1}')
\supset
\bigsqcup_{i=-\widehat i_0}^{r-1-\widetilde i_0}
((B_{n,i}-f_n)\cap (B_{n,i}'-f_n')+ih_n).
 $$
Now \thetag{5-3} yields that  
$$
\align
\lambda(A\times A' )
&\ge
\frac{\sum_{i=-\widehat i_0}^{r-1-\widetilde i_0}\#((B_{n,i}-f_n)\cap (B_{n,i}'-f_n'))
}{h_{n+1}+f_{n+1}-f_{n+1}'}+\overline{ o}(1)\\
&=\frac{(h_{n}+f_{n}-f_{n}')\sum_{i=-\widehat i_0}^{r-1-\widetilde i_0}\lambda([B_{n,i}]_n\times [B_{n,i}']_n)}{h_{n+1}+f_{n+1}-f_{n+1}'} +\overline{ o}(1)
\endalign
$$
as $D\ni n+1\to\infty$.
Applying \thetag{5-2} we obtain that
$$
\lambda(A\times A')+\overline{ o}(1)\ge
\frac{1+\delta-\delta'}{r+\delta-i_0}
\sum_{i=-\widehat i_0}^{r-1-\widetilde i_0}
\lambda\circ (T_{s(i+i_0)-s(i_0)}\times T'_{s'(i+i_0')-s'(i_0)})([A_n]_n\times[A_n']_n)
$$
as $D\ni n+1\to\infty$.
It follows that 
$$
\lambda\gg\sum_{i=-\widehat i_0}^{r-1-\widetilde i_0}\lambda\circ (T_{s(i+i_0)-s(i_0)}\times T'_{s'(i+i_0')-s'(i_0')}).
$$
Of course, $\lambda\circ (T_{s(i+i_0)-s(i_0)}\times T'_{s'(i+i_0')-s'(i_0')})\in J^e(T,T')$ for each 
$i=-\widehat i_0,-\widehat i_0+1,\dots,r-1-\widetilde i_0$.
By the extremal  property of the ergodic joinings, 
$$
\lambda\circ (T_{s(i+i_0)-s(i_0)}\times T'_{s'(i+i_0')-s'(i_0')})=\lambda
$$
 for each 
$i=-\widehat i_0,\widehat i_0+1,\dots,r-1-\widetilde i_0$.
\qed
\enddemo

We recall a classical lemma by  Furstenberg.

\proclaim{Lemma 5.2}
Given two standard probability  spaces $(X,\mu)$ and $(Y,\nu)$  and a measure $\lambda$ on $X\times Y$ whose pullbacks onto $X$ and $Y$ are $\mu$ and $\nu$ respectively, if there is  an ergodic $\nu$-preserving invertible transformation $R$ of $Y$
such that 
 $\lambda\circ(\text{Id}\times R)=\lambda$ then
$\lambda=\mu\times\nu$.
\endproclaim

\comment
there is  an infinite subset $D\subset\Bbb N$, an integer $r>1$ and two maps $s,s':\{0,1,\dots,r-1\}\to\Bbb Z_+$ such that
  $F_n=F_n'$,   $\#C_{n+1}=\#C_{n+1}'=r$, $s_{n+1}'=s$ and $s_{n+1}'=s'$  for each $n+1\in D$.
  \roster
  \item"(i)"
  If   $s\perp s'$  then 
  $T$ and $S$ are not isomorphic.
  \item"(ii)"
  If   $s\perp s'$ and  for each $k\in\{1,2,\dots,\max\{\max s,\max s'\}\}$, either  $T_k$ is ergodic or  $S_k$ is ergodic then 
  $T$ and $S$ are disjoint.
  \endroster
\endproclaim

\demo{Proof?} 
Let $\lambda\in J^e(T,S)$.
Since $s\perp s'$, it follows from Lemma~5.1 that $\lambda\circ (T_k\times\text{Id})=\lambda$ for some non-zero integer $k\le \max\{\max s,\max s'\}$.

(i) If $T$ and $S$ are isomorphic then the graph of this isomorphism is an ergodic joining
of $T$ and $S$.
Denote it by $\lambda$.
The fiber measures of this measures are singletones.
This contradicts tothe fact $\lambda\circ (T_k\times\text{Id})=\lambda$.

(ii)  By Lemma~5.2, $\lambda=\mu\times\mu'$.
\qed
\enddemo

\endcomment

Now we can prove  Theorem~I.

\comment
Since $s\perp s'$, it follows that $\lambda\circ (T^k\times\text{Id})$ for some $k\ne 0$.
Therefore $\lambda=\mu\times\mu'$.

No! That is not all.
Set now $a_n:=h_n-f_n'$, $b_n:=h_n-f_n$.
Without loss of generality we may assume that $c_{n+1}'$ is not maximal? (Or for infinitely many $n$).
If $i\in[a_n,b_n]$ then
$S^ix'=(f_n'+i-h_n+\widehat s, c_{n+1}'(i_0+1), c_{n+2},\dots)$.
Hence
$$
A_{n+1}'-f_{n+1}'+h_{n+1}+\widehat s=A_{n}'-f_n'+...
$$
We get $s(i)=s'(i+i_0-r)-s'(i_0)+\widetilde s$ for $i=r-i_0,\dots,r$.
\endcomment

\proclaim{Theorem 5.3} Let $T$ be a rank-one  $\Bbb Z$-action with  bounded parameters.
Suppose that $T$ is not rigid and that $T$ is totally ergodic. 
Then $T$ has MSJ$_2$ and hence MSJ.
\endproclaim
\demo{Proof}  Without loss of generality (due to Lemma~1.5) we may assume without loss of generality that $T$ is the $(C,F)$-action associated with bounded parameters $(C_n,F_{n-1})_{n\ge 1}$ and $F_n=\{0,\dots,h_n-1\}$ for some $h_n>0$.
Let $r_n$, $\sigma_n$ and $s_n$ be the same objects as in the beginning of \S 5.2.
Passing, if necessary,  to a bounded telescoping, we may assume that $r_n\ge 4$ for all $n\ge 1$.
Take a joining $\lambda\in J_2^e(T)$ which is not a graph-joining.
Given a generic point  
$$
(x,x')=((f_m,c_{m+1},\dots),(f_m',c_{m+1}',\dots))\in X\times X
$$ 
of $\lambda$, we let $D_{x,x'}:=\{d>m\mid c_d= c_d'\}$.
Let  $\widetilde D$ denote the set of $\lambda$-generic points $(x,x')$ for which $D_{x,x'}$ is finite.
Since $\lambda$ is ergodic and  $\widetilde D$  is $(T_g\times T_g)_{g\in \Bbb Z}$-invariant, it follows that either $\lambda(\widetilde D)=1$ or $\lambda(\widetilde D)=0$.

For each $n>0$, we denote by $c_n^\circ$ the only element of $C_n$ such that
$\# C_n/2-1\le\#(\{c\in C_l\mid c< c_n^\circ\})< \# C_l/2$.
Since the projection of $\lambda$ from $X\times X$ to each of the two coordinates is $\mu$,
 the  subset 
$$
\widehat D:=\{x=(f_m,c_{m+1},\dots)\in X_m\mid \#\{l>m\mid c_l=c_l^\circ \text{ and } c_{l-1}=c_{l-1}^\circ\}\text{  is infinite}\}
$$
 is of full measure $\mu$ in $X_m$.

Consider first the case where $\lambda(\widetilde D)=1$.
Then  there exists a $\lambda$-generic point $(x,x')\in\widetilde D\cap\widehat D$.
Since the cutting-and-stacking parameters are bounded,
 there exist  $r,t\ge 4$, an infinite subset $D\subset\{m+1,m+2,\dots\}$ and  maps 
 $\sigma:\{1,\dots,r\}\to\Bbb Z_+$  and $\tau:\{1,\dots,t\}\to\Bbb Z_+$ such that $\sigma_{n-1}=\tau$, $\sigma_n=\sigma$,  $c_{n}=c_{n}^\circ\ne c_{n}'$ and $c_{n+1}=c_{n+1}^\circ\ne c_{n+1}'$ whenever $n+1\in D$ (as above, $c_{n+1}$ and $c_{n+1}'$ are coordinates of $x$ and $x'$ respectvely).
 Let $s$ be the integral of $\sigma$ and let $u$ be an integral of $\tau$.
Then $s_n=s$ whenever $n+1\in D$.
It follows from the proof of Proposition~5.1 that there are $i_0,i_0'\in\{0,\dots,r-1\}$ such that
$c_{n+1}=s(i_0)+i_0h_n$ and $c_{n+1}'=s(i_0')+i_0'h_n$  whenever $n+1\in D$ and 
$$
\lambda=\lambda\circ(T_{s(i+i_0)-s(i_0)}\times T_{s(i+i_0')-s(i_0')})\
   \tag5-4
  $$
for each $i\in\{-\widehat i_0,\dots,r-1-\widetilde i_0\}$.
Since $c_{n+1}\ne c_{n+1}'$, it follows that $i_0\ne i_0'$.
The condition $c_{n+1}=c_{n+1}^0$ yields that $|i_0-i_0'|<r/2$.
For concreteness, we will assume that  $i_0>i_0'$.
Since $T$ is totally ergodic and~\thetag{5-4} holds, we deduce form Lemma~5.2 that $\lambda\ne\mu\times\mu$ only if  $s(i+i_0)-s(i_0)=s(i+i_0')-s(i_0')$ for all $i\in\{-\widehat i_0,\dots,r-1-\widetilde i_0\}$.
The latter condition 
is equivalent to the following one:
  $\sigma(i)=\sigma(i+(i_0-i_0'))$ for each $i=1,\dots, r-1-(i_0-i_0')$, i.e.
 $\sigma$ is $(i_0-i_0')$-periodic.
In a similar way we obtain that $\tau$ is $|j_0'-j_0|$-periodic for some $j_0,j_0'\in\{0,\dots,t-1\}$ such that $c_{n}=u(j_0)+j_0h_{n-1}$ and $c_{n}'=u(j_0')+j_0'h_{n-1}$ and $j_0\ne j_0'$
and $|j_0-j_0'|<t/2$.
It follows from Lemma~4.4(i) that 
$$
\tau\diamond\sigma\quad\text{is \ $t(i_0-i_0')$-periodic.}\tag 5-5
$$
Utilizing a  bounded telescoping we now pass from $(C_n,F_{n-1})_{n=1}^\infty$
to a new $(C,F)$-sequence $(\widetilde C_n,\widetilde F_{n-1})_{n=1}^\infty$ in such a way that 
$$
\{\widetilde C_n\mid n\in\Bbb N\}\supset\{C_n+C_{n+1}\mid n+1\in D\}.
$$
Of course, the ``telescoping isomorphism'' maps generic points into generic points.
We recall that if $(C_{n},F_{n})\sim (t,\tau)$ and $(C_{n+1},F_{n+1})\sim (r,\sigma)$ 
then $(C_n+C_{n+1},F_{n+1})\sim (tr,\tau\diamond\sigma)$. 
Denote by $v$ the integral of $\tau\diamond
\sigma$.
Let $k:=j_0+ti_0$ and $k':=j_0'+ti_0'$.
A direct computation yields that
$c_n+c_{n+1}=v(k)+kh_{n-1}$ and $c_n'+c_{n+1}'=v(k')+k'h_{n-1}$. 
Hence repeating the  above argument  for the ``partial'' periodicity of $\sigma$ we obtain that  
if  $\lambda\ne \mu\times\mu$ then $\tau\diamond
\sigma$ is $(k-k')$-periodic.
Now~\thetag{5-5} and Lemma~4.4(ii) yield that  $\tau\diamond
\sigma$ is $|j_0-j_0'|$-periodic.
By Lemma~4.4(iii), $\sigma(1)=\cdots=\sigma(r-1)$.
As was shown in the proof of Theorem~4.5, this  implies that $T$ is rigid, a contradiction.

Now consider the second case where $\lambda(\widetilde D)=1$.
If there is a $\lambda$-generic point $(x, x')$ such that $\Bbb N\setminus D_{x,x'}$ is finite then $x'$ belongs to the $T$-orbit of $x$.
It is well known that then $\lambda$ is a graph-joining.
Thus it remains to consider the case where $D_{x,x'}$ and $\Bbb N\setminus D_{x,x'}$ are both infinite.
Since the cutting-and-stacking parameters are bounded,
 there exist  $r,t\ge 4$, an infinite subset $D\subset\{m+1,m+2,\dots\}$ and  maps 
 $\sigma:\{1,\dots,r\}\to\Bbb Z_+$  and $\tau:\{1,\dots,t\}\to\Bbb Z_+$ such that $\sigma_{n-1}=\tau$, $\sigma_n=\sigma$,  $c_{n}\ne c_{n}'$ and $c_{n+1}= c_{n+1}'$ whenever $n+1\in D$ (as above, $c_{n+1}$ and $c_{n+1}'$ are coordinates of $x$ and $x'$ respectively).
  Passing to a telescoping (as in the first case) one can show that if $\lambda\ne\mu\times\mu$ then $\tau\diamond\sigma$ is $|j_0-j_0'|$-periodic for some $j_0,j_0'\in\{1,\dots,t\}$ and $j_0\ne j_0'$.
By Lemma~4.4(iii),  $\sigma(1)=\cdots=\sigma(r-1)$ and hence $T$ is rigid, a contradiction.
Thus in every  case we have that $\lambda$ is either $\mu\times\mu$ or a graph-joining.
Hence $T$  has MSJ$_2$.

\comment
In the latter case there exists a $\lambda$-generic point $(x,x')$ such that $D_{x,x'}$ is finite.
Then $x'$ belongs to the $T$-orbit of $x$. 
A standard reasoning implies that $\lambda$ is a graph-joining, a contradiction.
Hence from now on we may assume without loss of generality that $\lambda(\widetilde D)=0$.


  Indeed, otherwise, i.e. in the situation that $i_0=i_0'$ for each $i_0\in\{0,\dots,r-1\}$ for every choice of $D$, due to the boundedness condition, we would have  that $c_n=c_n'$ eventually in $n$ and hence $x'$ belongs to the $T$-orbit of $x$.
  Then a standard reasoning  implies that  $\lambda$ is a graph-joining.
  
  Next, we may not to  consider  generic points $(x,x')$ which are ``extreme'', i.e. such points  that $\{c_{n+1},c_{n+1}'\}=\{0, s_n(r_n-1)+(r_n-1)h_{n}\}$ 
   eventually in $n$, because the
 $\lambda$-measure of the set of all extreme points is $0$\footnote{This follows from the fact that the projection of the set of extreme points on each coordinate of $X\times X$ is of $\mu$-measure 
 $0$.}.
  Therefore, choosing $D$ in an appropriate way, we may assume that 
  $
  \{c_{n+1},c_{n+1}'\}\ne\{0, s(r-1)+(r-1)h_{n}\}  $
   for each sufficiently large $n+1\in D$.
   This means that
   $$
   \{i_0,i_0'\}\ne\{0,r-1\}.\tag5-5
  $$
  Since $T$ is totally ergodic and \thetag{5-4} holds, we deduce form Lemma~5.2 that $\lambda\ne\mu\times\mu$ only if  $s(i+i_0)-s(i_0)=s(i+i_0')-s(i_0')$ for all $i\in\{-\widehat i_0,\dots,r-1-\widetilde i_0\}$.
The latter condition (in view of  \thetag{5-5}) is equivalent to the following one:
  $\sigma(i)=\sigma(i+|i_0'-i_0|)$ for each $i=1,\dots, r-1-|i_0'-i_0|$, i.e.
 $\sigma$ is $|i_0'-i_0|$-periodic.
 Arguing as in the proof of Theorem~4.5 we obtain that $T$ is rigid, a contradiction.
 Hence $\lambda=\mu\times\mu$ and thus $T$ has MSJ$_2$.
 \endcomment

 Since $T$ is a rank-one action with bounded parameters, it is a routine to show that $T$ is {\it partially rigid}, i.e. there exist a sequence $n_k\to\infty$ and a real  parameter $\eta>0$ such that $\liminf_{k\to\infty}\mu(T_{n_k}A\cap A)\ge \eta\mu(A)$ for some for each Borel subset $A\subset X$.
Therefore every factor\footnote{A {\it factor} of $T$ is an invariant sub-$\sigma$-algebra of Borel subsets in $X$ or, more rigorously, the restriction of $(T,\mu)$ to this sub-$\sigma$-algebra.} of $T$ is also partially rigid (with the same parameter $\eta$) and hence
$T$ has no factor with Lebesgue spectrum.
Since $T$ has MSJ$_2$,
 it follows now from \cite{GlHoRu, Theorem~4} that $T$ has MSJ of all orders. 
\qed
\enddemo

It follows from this theorem and Example~4.10 that the Chacon 3-cuts map and the Chacon 2-cuts maps have MSJ.

The next claim follows from Theorem~5.3 and the properties of transformations with MSJ (see, for instance, \cite{dJRu, Corollary~6.5}).

\proclaim{Corollary 5.4} If  $n,m>0$ and $n\ne m$ then $T_n$ and $T_m$ are disjoint.
\endproclaim

We now deduce one more corollary from Proposition 5.1.
We first note that if $T$ and $T'$ are commensurate and  the parameters $(C_n,F_{n-1})_{n=1}^\infty$ and $(C_n',F_{n-1}')_{n=1}^\infty$ are bounded then $r_n=r_n'$ (i.e. $\# C_n=\# C_n'$) eventually.

\proclaim{Corollary 5.5} Let $T$ and $T'$ be two commensurate $(C,F)$-actions of $\Bbb Z$ associated with bounded parameters $(C_n,F_{n-1})_{n=1}^\infty$ and $(C_n',F_{n-1}')_{n=1}^\infty$ respectively.
Let either $T$ or $T'$ be non-rigid.
\roster
\item"$(i)$"
Then $T$ and $T'$ are isomorphic if and only if $s_n=s_n'$ eventually.
\item"$(ii)$" If 
$s_n\ne s_n'$ for infinitely many $n$
and
for each $n>0$, either $T_n$ or $T'_n$ is ergodic then
$T$ and $T'$ are disjoint.
\endroster
\endproclaim

We preface the proof of the corollary with a simple auxiliary lemma.

\proclaim{Lemma 5.6} Given  positive integers $r$ and $q$,  let $\sigma,\omega:\{1,\dots,r\}\to\Bbb Z_+$ and $\alpha,\beta:\{1,\dots,q\}\to\Bbb Z_+$ be four maps.
Suppose that there is $i_0>0$ such that $i_0\le rq/2$, $r$ does not divide $i_0$, and
$\sigma\diamond\alpha(i)=\omega\diamond\beta(i+i_0)$ for each $i\in\{1,\dots,rq-i_0-1\}$.
Then there is $p\le q/2$ such that  $\alpha(1)=\alpha(2)=\cdots=\alpha(q-p)$ and $\beta(p)=\beta(p+1)=\cdots=\beta(q)$.
\endproclaim
\demo{Proof}
There is $z\in\{1,\dots, r-1\}$ such that $i_0=j_0r+z$ for some non-negative integer $j_0< q/2$.
Then for each $j\in\{1,\dots,q-j_0-1\}$,
$$
\sigma(r)+\alpha(j)=\sigma\diamond\alpha(jr)=\omega\diamond\beta(jr+i_0)=\omega\diamond\beta((j+j_0)r+z)=\omega(z).
$$
Hence $\alpha(1)=\alpha(2)=\cdots=\alpha(q-j_0-1)$.
In a similar way, for each $j\in\{j_0+1,\dots,q\}$,
$$
\omega(r)+\beta(j)=\omega\diamond\beta(jr)=\sigma\diamond\alpha(jr-i_0)=\omega\diamond\beta((j-j_0-1)r+r-z)=\omega(r-z).
$$
Hence $\beta(q)=\beta(q-1)=\cdots=\beta(j_0+1)$.
\qed
\enddemo

\demo{Proof of Corollary 5.5} (i) The ``if'' part is trivial. 
We now prove the ``only if'' part.
Thus suppose that $T$ and $T'$ are isomorphic.
Let $T$ be rigid.
The case where $T'$ is rigid is considered in a similar way.
Let $\lambda\in J^e(T,S)$ denote the graph joining generated by this isomorphism.
Since the parameters $(C_n,F_{n-1})_{n=1}^\infty$ and 
$(C_n',F_{n-1}')_{n=1}^\infty$ are bounded,
there exist $r>0$, an infinite subset $D\subset\Bbb N$ 
 and  maps $s,s':\{0,\dots,r-1\}\to\Bbb Z_+$ such that
$s_n=s$ and   $s_n'=s'$ whenever $n+1\in D$.
Let $i_0$ stand for the integer part of $(r-1)/2$.
By Proposition~5.1,
there is $i_0'\in\{0,1,\dots,r-1\}$ such that 
 $$
 \lambda=\lambda\circ (T_{s(i+i_0)-s(i_0)-s'(i+i_0')+s'(i_0')}\times \text{Id})\tag5-6
 $$
  for every $i$ such that $-\min(i_0,i_0')\le i\le r-1-\max(i_0,i_0')$.
  Since the conditional measures of $\lambda$ corresponding to the projection $X\times X'\ni(x,x')\mapsto x'\in X'$ are delta-measures almost everywhere on $X'$ and $T$ is free, it follows from \thetag{5-6} that $T_{s(i+i_0)-s(i_0)-s'(i+i_0')+s'(i_0')}=\text{Id}$ and hence
  $s(i+i_0)-s(i_0)-s'(i+i_0')+s'(i_0')=0$
  for every $i$ such that $-\min(i_0,i_0')\le i\le r-1-\max(i_0,i_0')$.
  This is equivalent to
   $$
   \sigma(i)=\sigma'(i+|i_0'-i_0|)\text{ or  }\sigma'(i)=\sigma(i+|i_0'-i_0|)\tag5-7
   $$
   (depending on the sing of  $i_0-i_0'$)
    for each $i=1,\dots, r-1-|i_0'-i_0|$.
   We note that $|i_0'-i_0|\le r/2$.
   Passing to a subsequence in $D$ we may assume also that the 
   sequences  $(\sigma_{n})_{n\in D}$, $(\sigma_{n+1})_{n\in D}$, $(\sigma_{n}')_{n\in D}$ and $(\sigma_{n+1}')_{n\in D}$ are all constant, i.e.
    there are $p,q>1$ and maps $\alpha,\alpha':\{1,\dots,q\}\to\Bbb Z_+$ and 
   $\beta,\beta':\{1,\dots,q'\}\to\Bbb Z_+$
such that $\sigma_n=\alpha$, $\sigma_{n+1}=\beta$, $\sigma_n'=\alpha'$ and  $\sigma_{n+1}'=\beta'$ for each $n\in D$.
Thus we now obtain that 
$(C_{n},F_{n})\sim(r,\sigma)$, 
$(C_{n+1},F_{n+1})\sim(q,\alpha)$, $(C_{n+2},F_{n+2})\sim(p,\beta)$,
$(C_{n}',F_{n}')\sim(r,\sigma')$, 
   $(C_{n+1}',F_{n+1}')\sim(q,\alpha')$ and $(C_{n+2}',F_{n+2}')\sim(p,\beta')$ for each $n\in D$.
   Since $(C_n+C_{n+1}+C_{n+2}, F_{n+2})\sim(pqr,  \sigma\diamond(\alpha\diamond\beta))$
   and 
   $(C_n'+C_{n+1}'+C_{n+2}', F_{n+2}')\sim(pqr,  \sigma'\diamond(\alpha'\diamond\beta'))$,
 we pass   to a corresponding bounded telescoping to deduce from~\thetag{5-7} that
   $$
   \align
   \sigma\diamond(\alpha\diamond\beta)(i)&=\sigma'\diamond(\alpha'\diamond\beta')(i+j_0)\quad\text{or}\\
    \sigma'\diamond(\alpha'\diamond\beta')(i)&=\sigma\diamond(\alpha\diamond\beta)(i+j_0)
   \endalign
   $$
   for  some positive integer $j_0\le pqr/2$ and each $i=1,\dots,pqr-j_0$. 
 It follows from Lemma~5.6 that $\alpha\diamond\beta(1)=\cdots =\alpha\diamond\beta(l_0)$ 
 or $\alpha\diamond\beta(pq)=\cdots=\alpha\diamond\beta(pq-l_0)$ for some positive integer $l_0\le pq/2$.
 This implies, in turn, that $\alpha(1)=\cdots=\alpha(q-1)$
 is constant.
 Hence $C_{n+1}$ is an arithmetic sequence for each $n\in D$.
 It follows now from Theorem~4.5 (or, more precisely, from the claim in the beginning of the proof of Theorem~4.5) that $T$ is rigid, a contradiction.

   (ii) is proved in a similar way.
   We only note that  \thetag{5-6}, the ergodicity condition in the statement of (ii) and Lemma~5.2 imply that
   $\lambda=\mu\times\mu'$, i.e. that $T$ and $S$ are disjoint whenever
   $s(i+i_0)-s(i_0)-s'(i+i_0')+s'(i_0')\ne 0$
  for some $i$ such that $-\min(i_0,i_0')\le i\le r-1-\max(i_0,i_0')$.
  Thus, if  $T$ and $S$ are not disjoint then we obtain as in (i) that~\thetag{5-7} holds.
As we have already shown  in the proof of (i), the condition  \thetag{5-7} (because it holds for every bounded telescoping of the original $(C,F)$-parameters of $T$ and $T'$) implies that $T$ is rigid, a contradiction.
   \qed
\enddemo

We consider an application of Corollary~5.5 to the inverse problem.

\proclaim{Corollary 5.7} Let $T$ be a $(C,F)$-action of $\Bbb Z$ with bounded parameters.
If $T$ is not rigid then
 $T$ and $T^{-1}$ are 
 isomorphic if and only if  $C_n^*=C_n$ eventually\footnote{For the definition of $C_n^*$ (and $F_n^*$ used in the proof of Corollary~5.7) see Lemma~2.8.}.
 If, moreover, $T$ is totally ergodic and $T$ is not isomorphic to $T^{-1}$ then $T$ and $T^{-1}$ are disjoint.
 \endproclaim

To prove the corollary we need  a measure theoretical analogue of Lemma~2.1.

\proclaim{Lemma 5.8} Let $T=(T_g)_{g\in G}$ and $T'=(T_g')_{g\in G}$ be two ergodic measure preserving free
 $G$-actions on $\sigma$-finite standard measure spaces  $(X,\mu)$ and $(X',\mu')$ respectively.
Let $A$ be a  Borel subset in $X$, let $A'$ be a Borel  subset in $X'$ such that $\mu(A)=\mu'(A)>0$ and let $\theta:A\to A'$  be
a measure preserving Borel bijection.
If, given $x\in A$ and $g\in G$,   we have that $T_gx\in A$ if and only if $T'_g\theta x\in A'$ and, moreover,  $\theta T_g x=T'_g\theta x$ 
then $T$ and $T'$ are (measure theoretically) conjugate.
\endproclaim
\demo{Proof}
Since $T$ is ergodic, there is a countable partition of $X\setminus A$ into subsets $A_g$, $g\in G\setminus\{1\}$, such that $T_{g}^{-1}A_g\subset A$ for each $g\in G$, $g\ne 1$.
We now set 
$$
\widetilde \theta x=
\cases
\theta x, &\text{if }x\in A\\
T_g'\theta T_g^{-1}x, &\text{if }x\in A_g.
\endcases
$$
It is routine to verify that $\widetilde\theta$ is Borel, one-to-one and onto (mod 0).
Of course, $\theta$ is measure preserving.
It is obvious  that $\widetilde\theta$ intertwines $T$ with $T'$.
\qed
\enddemo

\demo{Proof of Corollary 5.7} We note that $T^{-1}$ is a $(C,F)$-action associated with the sequence $(C_n^*,F_{n-1}^*)_{n\ge 1}$.
It is easy to verify that the sequence $(C_n^*,F_{n-1})_{n\ge 1}$  satisfies the conditions (I)--(III) and \thetag{1-4}.
It is bounded.
The $(C,F)$-action  of $\Bbb Z$ associated with $(C_n^*,F_{n-1})_{n\ge 1}$ is isomorphic to $T^{-1}$ by Lemma~5.8.
It remains to apply Corollary~5.5.
\qed
\enddemo

Corollary 5.7 yields immediately that the Chacon  map $T$ with 3 cuts and the Chacon map $R$ with 2 cuts are (measure theoretically) not isomorphic to their inverses.
Moreover, in view of Example~4.10,  $T$ is disjoint with  $T^{-1}$  and $R$ is disjoint with $R^{-1}$.

\subhead 5.3. Light mixing and bounded rank-one constructions
\endsubhead
Our purpose here is to show that the class of all rank-one transformations with bounded parameters is strictly bigger than the subclass of adapted rank-one transformations with bounded parameters.
For that we will use the property of light mixing.
Let $S=(S_n)_{n\in\Bbb Z}$ be a measure preserving $\Bbb Z$-action on a standard probability space $(Y,\nu)$. 
We recall that $S$ is called {\it lightly mixing} (see, for instance, \cite{Si}) if for all subsets $A,B\subset Y$ of positive measure,
$$
\liminf_{n\to\infty}\nu(S_nA\cap B)>0.
$$
By \cite{FrKi}, the Chacon map with 2 cuts is lightly mixing.
We now prove the following theorem.

\proclaim{Theorem 5.9} Let  $T=(T_n)_{n\in\Bbb Z}$ be  an adapted  rank-one action of $\Bbb Z$ and let
$T_1\sim (r_n,\sigma_n)_{n\ge 1}$ for some sequence of integers $r_n$ and maps $\sigma_n:\{1,\dots,r_n\}\to\Bbb Z_+$.
If there is $K>0$ such that $\sup_{n\ge 1}\max_{1\le i\le r_n}\sigma_n(i)\le K$ then
there are a sequence $n_m\to +\infty$ and  a polynomial $P(Z)=\nu_0+\nu_1 Z+\cdots+\nu_K Z^K$ with real coefficients such that $\min_{0\le i\le K}\nu_i\ge 0$, 
 $\sum_{0\le i\le K}\nu_i=1$ and
$$
T_{-h_{n_m}}\to P(T_1)\tag5-8
$$ 
in the weak operator topology\footnote{We consider $T_i$ in \thetag{5-8} as the unitary Koopman  operators in $L^2(X,\mu)$ generated by the transformation $T_i$, $i\in\Bbb Z$.}.
It follows that $T$ is not lightly mixing.
 \endproclaim

\demo{Proof}
Represent $T$ as a $(C,F)$-action associated with a  sequence $(C_{n}, F_{n-1})_{n\ge 1}$
such that
$F_n=\{0,1,\dots,h_n-1\}$ and 
 $(C_{n+1},F_{n+1})\sim(r_n,\sigma_n)$ for each $n\ge 0$.
 As above, we will denote by $s_n$ the integral of $\sigma_n$.
Since $T$ is adapted, i.e. $\sigma_n(r_n)=0$ for each $n$, it follows that
$$
\sup_{m,n>0}\max_{1\le i\le r_n\cdots r_{n+m}}(\sigma_n\diamond\cdots\diamond\sigma_{n+m})(i)\le K.
$$
Therefore, passing, if necessary, to a telescoping we may assume without loss of generality that $r_n\to\infty$ while the upper bound on the number of consecutive spacers does not increase.
Let $A$ and $B$ be cylinders in $X$.
Then $A=[A_n]_{n}$ and $B=[B_n]_n$ for some subsets $A_n$ and $B_n$ in $\Bbb Z$ such that $A_n\cup B_n\subset \{0,\dots,h_n-1-K\}$ for each sufficiently large $n$.
Then
$$
\align
\mu(T_{-h_n}A&\cap B)=\mu(T_{-h_n}[A_n+C_{n+1}]_{n+1}\cap B)\\
&=
\mu\left(\bigsqcup_{1\le i< r_n}[\sigma_n(i)+A_n+s_n(i-1)+(i-1)h_n]_{n+1}\cap [B_n]_n\right)\pm \frac 1{r_n}\\
&=
\mu\left(\bigsqcup_{1\le i< r_n}[(\sigma_n(i)+A_n)\cap B_n+s_n(i-1)+(i-1)h_n]_{n+1}\right)\pm \frac 1{r_n}\\
&=
\frac 1{r_n}\mu\left(\bigsqcup_{1\le i< r_n}[(\sigma_n(i)+A_n)\cap B_n]_n\right)\pm \frac 1{r_n}\\
&=
\frac 1{r_n}\mu\left(\bigsqcup_{1\le i\le r_n}[\sigma_n(i)+A_n]_n\cap [B_n]_n\right)\pm \frac 2{r_n}\\
&=
\frac 1{r_n}\mu\left(\bigsqcup_{1\le i\le r_n}T_{\sigma_n(i)}A\cap B\right)\pm \frac 2{r_n}.
\endalign
$$
Denote by $\vartheta_n$ the image of the equidistribution on $\{1,\dots,r_n\}$ under $\sigma_n$.
Then $\vartheta_n$ in a distribution on $\{0,1,\dots,K\}$ and 
$$
\mu(T_{-h_n}A\cap B)=\sum_{i=1}^{K}\vartheta_n(i)\mu(T_iA\cap B)\pm\frac 2{r_n}.
\tag5-9
$$
Let $\nu$ stand for a limit point of the sequence $(\vartheta_n)_{n=1}^\infty$.
Passing to the limit in \thetag{5-9} along a corresponding subsequence we obtain \thetag{5-8}.

To prove the second claim of the proposition, select subsets $A$ and $B$ of positive measure  in such a way that
$B\cap\bigcup_{i=0}^K T_iA=\emptyset$.
Then, in view of \thetag{5-9},  we have that $\mu(T^{-h_{n_m}}A\cap B)\to 0$ as $
m\to\infty$.
Hence $T$ is not lightly mixing.
\qed
\enddemo

\subhead{5.4. Quadchotomy theorem}\endsubhead
We conclude this section with a quadchotomy theorem on  ``structure'' of rank-one transformations with bounded parameters (Theorem~M).

\proclaim{Theorem 5.10} Let  $T$ be a  $(C,F)$-action of $\Bbb Z$ associated with a bounded sequence $(C_n,F_{n-1})_{n\ge 1}$ and let $
F_n=\{0,1,\dots,h_n-1\}$
  for some  integers $h_n\ge 1$ for each  $n\ge 0$.
Let $(C_{n+1},F_{n+1})\sim(r_n,\sigma_n)$ for each $n\ge 0$ and 
$K:=\sup_{n\ge 1}\sup_{1\le j\le r_n}\sigma_n(j)$.
Then one of the following four properties holds:
\roster
\item"\rom(i)" $T$ has MSJ (in particular, $T$ is weakly mixing and $C(T)=\{T_n\mid n\in\Bbb Z\})$.
 \item"\rom(ii)"
 $T$ is non-rigid, the group $\Lambda_T\subset\Bbb T$ of eigenvalues of $T$ is nontrivial but finite and
 the order of each $\lambda\in\Lambda_T$ does not exceed  $K$.
 For each  $\rho\in J_2^e(T)$ which is neither a graph-joining nor $\mu\times\mu$, there is $\lambda\in\Lambda_T\setminus\{1\}$
 and $n>0$ such that $\lambda^n=1$ and
 $\frac 1n\sum_{i=0}^{n-1}\rho\circ(I\times T_i)=\mu\times\mu$.
\item"\rom(iii)" $T$ is rigid,
the group $\Lambda_T\subset\Bbb T$ of eigenvalues of $T$ is finite and
 the order of each $\lambda\in\Lambda_T$ does not exceed  $K$.
\item"\rom(iv)" $T$ is an odometer of bounded type.
\endroster
  
\endproclaim

\demo{Proof}
Consider two cases: $T$ is rigid and $T$ is not rigid.
If $T$ is rigid then we apply Proposition~4.6 and obtain either (iii) or (iv).
If $T$ is not rigid but $T$ is totally ergodic, we apply Theorem~5.3 and obtain (i).
It remains to show that if $T$ is non-rigid but not totally ergodic then (iii) holds.
Since $T$ is non-rigid,
it follows from  Theorem~4.5 that  there is $M>0$ such that for each $n>0$, the sum $C_n+\cdots+C_{n+M}$ is not an arithmetic sequence.
Let $\lambda\in\Lambda_T$  and  let $\xi:X\to\Bbb C$ be a corresponding non-zero measurable eigenfunction for $T$
such  that $|\xi|=1$ almost everywhere on $X$.
Fix $\epsilon>0$.
Then we can  select $n>0$, an element $f\in F_{n-1}$ and a complex number $z\in\Bbb T\subset\Bbb C$ such that the set  $A:=\{x\in X\mid |\xi(x)-z|<\epsilon\}$ is of positive measure  and  the cylinder $[f]_{n-1}$ is $(1-\epsilon/r^M)$-{\it full of $A$}.
Let $m\ge 0$ be the smallest  integer such that the set $C:=C_n+\cdots+C_{n+m}$ is not an arithmetic sequence.
Of course, $m\le M$.
It follows that for each $c\in C$, the cylinder $[f+c]_{n+m}$ is $(1-\epsilon)$-full of $A$.
We note that $C=\{ih_{n-1}+s(i)\mid 0\le i< \# C\}$
for some map $s:\{0,1,\dots,\# C-1\}\to\Bbb Z_+$  and $0\le s(i+1)-s(i)\le \sigma_{n+m-1}(i+1)\le K$ for each $i=0,1,\dots,\# C-2$.
   Since for each $i\in\{0,\dots,\# C-2\}$, we have that  $T_{h_{n-1}+s(i+1)-s(i)}[f+ih_{n-1}+s(i)]_{n+m}= [f+(i+1)h_{n-1}+s(i+1)]_{n+m}$, there exists
 $x_i\in A $ such that $T_{h_{n-1}+s(i+1)-s(i)}x_i\in  A$.
 This yields that 
 $$
\max_{0\le i\le \# C-2} |\lambda^{h_{n-1}+s(i+1)-s(i)}-1|\le 2\epsilon.
$$
In turn, this inequality implies that $\lambda^k=1$ for some $k\in\{1,\dots, K\}$.
Indeed, otherwise we would obtain that $s(i+1)-s(i)=s(j+1)-s(j)$ for all $0\le i,j\le\# C-1$
and hence $C$ is an arithmetic sequence, a contradiction.
Thus the first statement of (iii) is proved.
If $\rho\in J_2^e(T)$ is not a graph-joining, it follows from the proof of Theorem~5.3
that $\rho(I\times T_n)=\rho$ for some $n>0$.
If $T_n$ is ergodic then $\rho=\mu\times\mu$ by Lemma~5.2.
Otherwise, there is $\lambda\in \Lambda_T\setminus\{1\}$ such that $\lambda^n=1$.
Then $\omega:=\frac 1n\sum_{i=0}^{n-1}\rho\circ(I\times T_i)\in J_2(T)$ and $\omega(I\times T_1)=\omega$.
Hence $\omega=\mu\times\mu$ by Lemma~5.2.
\qed
\enddemo

\head Appendix. Symbolic representation for  rank-one transformations
\endhead

 Let $T$ be a $(C,F)$-action of $\Bbb Z$ associated with a sequence $(C_n,F_{n-1})_{n\ge 1}$ satisfying (I)--(III) and such that $F_n=\{0,1,\dots,h_n-1\}$ for all $n$ and some sequence $(h_n)_{n=1}^\infty$.
The following definition is equivalent to \cite{AdFePe, Definition~5.1}.

\definition{Definition A.1} 
For $k\in\Bbb Z_+$, a  rank-one $(C,F)$-action $T=(T_i)_{i\in\Bbb Z}$
is 
{\it essentially $k$-expansive} if the partition $\Cal P_k$ of  $X$, given by $\Cal P_k:=\{[0]_k,X\setminus[0]_k\}$, generates the entire Borel $\sigma$-algebra under the action of $T$.
\enddefinition

For  $k\in\Bbb Z_+$ and $x\in X$, we set 
$$
\Cal P_k(x):=\cases 0,&\text{if }x\in[0]_k\\ 1,&\text{otherwise}.\endcases
$$
 It is easy to see that the map
$$
\phi^{(k)}:X\ni x\mapsto \phi^{(k)}(x)=(\Cal P_k(T_ix))_{i\in \Bbb Z}\in\{0,1\}^\Bbb Z
$$
is continuous. 
It intertwines $T$ with the shiftwise $\Bbb Z$-action $S$ on $\{0,1\}^\Bbb Z$.
Then  $T$ is essentially $k$-expansive  if and only if $\phi^{(k)}$ is an isomorphism of $(X,\mu,T)$ onto $(\{0,1\}^\Bbb Z,\mu\circ(\phi^{(k)})^{-1},S)$.
Thus $\phi^{(k)}$ provides a {\it constructive symbolic model} of $T$ (see a discussion about symbolic  models in \cite{AdFePe} and \cite{Fe2}).
Of course, if $T$ is $0$-expansive then $T$ is $k$-expansive for each $k>0$.
By \cite{Da6, Theorem~2.8, Remark~2.10}, each rank-one transformation
is (measure theoretically) isomorphic to a rank-one transformation $Q\sim (r_n',\sigma_n')_{n=1}^\infty$ with $\sigma_n'(r_n')>0$ for infinitely many $n$.
 If $Q$ is,  in addition, probability preserving  then---as was explained in \cite{AdFePe}---$Q$ is  essentially $0$-expansive in view of  \cite{Ka, Appendix}.
Below (see Theorem~A.3) we  slightly modify  the proof of \cite{Da6, Theorem~2.8} to obtain a finer condition on $(r_n',\sigma_n')_{n=1}^\infty$ which yields immediately  the $0$-expansiveness of $Q$.

The following proposition is a slight modification of \cite{AdFePe, Proposition~5.2}.
We give a  full proof of it to make the present paper  more self-contained.

\proclaim{Proposition A.2} Let $T_1\sim(r_n,\sigma_n)_{n=1}^\infty$ such that \thetag{1-5} and \thetag{1-6} hold.
If $\sigma_n(r_n)>\max\{\sigma_n(1),\dots,\sigma_n(r_n-1)\}$  for each $n>0$ then
$T$ is essentially $0$-expansive.
\endproclaim
\demo{Proof}
Let $\goth B_0$ stand for the smallest $T$-invariant sub-$\sigma$-algebra containing $[0]_0$.
Of course, if $[0]_{n}\in \goth B_0$ for each $n\ge 0$ then $\goth B_0$ coincides with the entire Borel $\sigma$-algebra.
We prove these inclusions by induction.
The inclusion $[0]_0\in\goth B_0$ is by definition.
Suppose that $[0]_n\in\goth B_0$ for some $n\ge 0$.
It follows from the condition of the proposition that 
$$
[0]_{n+1}=\{x\in [0_n]\mid T_{h_n+\sigma_n(i)}x\not\in[0]_n\text{ for each }i=1,\dots,r_n-1\},
$$
and hence $[0]_{n+1}=[0]_n\cap\bigcap_{i=1}^{r_n-1}T_{-h_n-\sigma_n(i)}(X\setminus[0]_n)\in\goth B_0$.
\qed
\enddemo
We now state and prove the main result of Appendix.

\proclaim{Theorem A.3} Every rank-one finite measure preserving $\Bbb Z$-action $T$ is (measure theoretically) isomorphic to a $(C,F)$-action which is  essentially $0$-expansive. \endproclaim

\demo{Proof} Let $T\sim(r_n,\sigma_n)_{n=1}^\infty$ such that \thetag{1-5} and \thetag{1-6} hold.
Passing, if necessary, to a telescoping, we may assume without loss of generality that
$\sum_{n=1}^\infty\frac1{r_n}<\infty$.
For each $n>0$, there is a unique $i_n\in\{1,\dots,r_n-1\}$
such that
$$
(i_n-1)h_n+\sum_{j=r_n-i_n}^{r_n}\sigma_n(j)\le\max_{1\le j<r_n}\sigma_n(j)<i_nh_n+\sum_{j=r_n-i_n+1}^{r_n}\sigma_n(j). \tag{A-1}
$$
Since $T$ is finite measure preserving,
$\sum_{n=1}^\infty\frac{\sum_{i=1}^{r_n}\sigma_n(i)}{r_1\cdots r_n}<\infty$ and  $\frac{h_n}{r_1\cdots r_{n-1}}\to \mu(X)<\infty$ as $n\to\infty$.
Therefore using the left-hand side inequality in~\thetag{A-1} we obtain that
$\sum_{n=1}^\infty\frac{(i_n-1)h_n}{r_1\cdots r_{n}}<\infty$ and hence 
$$
\sum_{n=1}^\infty\frac{i_n}{r_n}<\infty.\tag A-2
$$
Enumerate the elements of $C_n$ in the increasing order:
$C_n=\{c_n(0),\dots, c_n(r_n-1)\}$ with $0=c_n(0)<c_n(1)<\cdots<c_n(r_n-1)$.
We now set 
$$
C_n^*:=\{c_n(0),\dots,c_n(r_n-i_n-1)\}\quad\text{and}\quad F_{n-1}^*:=F_{n-1}
$$ 
for each $n>0$.
Then the  sequence $(C_n^*,F_{n-1}^*)_{n\ge 1}$ satisfies (I)--(III).
Denote by $T^*=(T^*_i)_{i\in\Bbb Z}$ the associated $(C,F)$-action.
Let $X^*$ stand for the corresponding  $(C,F)$-space  of $T^*$.
It is a routine to verify that $T^*\sim(r_n-i_n,\sigma_n^*)_{n=1}^\infty$, where
$$
\sigma_n^*(i):=\cases\sigma_n(i),&\text{if }1\le i<r_n-i_n\\ i_nh_n+\sum_{j=r_n-i_n+1}^{r_n}\sigma_n(j),&\text{if }j=r_n-i_n.\endcases\tag A-3
$$
It follows from the right-hand side inequality in \thetag{A-1}, \thetag{A-3} and Proposition~A.2 that $T$ is essentially $0$-expansive.
It remains to show that $T^*$ is  isomorphic to $T$.
We do this in the same way as in the proof of \cite{Da6, Theorem~2.8}.

Given $n\ge 0$, since the restriction of $\mu$ to $X_n:=F_n\times C_{n+1}\times C_{n+2}\times\cdots$ is proportional to the infinite product of the equidistributions on $F_n$ and $C_j$ for all $j>n$ and \thetag{A-2} holds, it follows from  the Borel-Cantelli lemma that  the subset
$$
Y_n:=\{x=(f_n,c_{n+1},c_{n+2},\dots)\in X_n\mid c_j\in C_j^*\text{ eventually in }j\}
$$
is of full measure in $X_n$.
We now define a Borel map $\phi_n:Y_n\to X^*$ by setting
$$
\phi_n(x):=(f_n+c_{n+1}+\cdots+c_m,c_{m+1},c_{m+2},\dots)\in F_m^*\times C_{m+1}^*\times
C_{m+2}^*\times\cdots\subset X^*,
$$
where $m$ is chosen in such a way that  $c_j\in C_j^*$ for all $j>m$.
Then we have that $Y_n\subset Y_{n+1}\subset\cdots$.
 Hence the subset $Y:=\bigcup_{n>0}Y_n$ is of full measure in $X$.
Moreover, $\phi_{n+1}\restriction Y_n=\phi_n$ for each $n$.
Hence we obtain a well-defined map $\phi:Y\to X^*$ by setting  $\phi\restriction Y_n=\phi_n$ for each $n>0$.
It is straightforward to verify that $\phi$ is a (measure theoretical) isomorphism of $X$ and $X^*$.
Moreover, $\phi$ intertwines $T$ with $T^*$, as desired.
\qed 
\enddemo

We now illustrate the explicit constructive procedure described in  the proof of~Theorem~A.3 on the example of  2-adic odometer.

\example{Example A.4} Let $T=(T_i)_{i\in\Bbb Z}$ be the 2-adic odometer,
i.e. $T_1\sim(r_n,\sigma_n)_{n=1}^\infty$, where $r_n=2$ and $\sigma_n\equiv 0$ for all $n>0$.
Equivalently, $T$ is the $(C,F)$-action associated with a sequence $(C_n,F_{n-1})_{n\ge 1}$ such that $F_n:=\{0,1,\dots, 2^n-1\}$ and $C_{n+1}=\{0,2^n\}$ for each $n\ge 0$.
Of course, $T$ is not essentially $0$-expansive because $\phi^{(0)}\equiv 0$.
Arguing as in the proof of Theorem~A.3 and passing to a telescoping we first obtain that $T$ is isomorphic to  a rank-one action $T'$ such that $T'\sim(2^n,\sigma_n')$, where $\sigma_n'\equiv 0$ for all $n$.
Then on every step of the cutting-and-stacking construction of $T'$ we ``replace'' the last (i.e. $2^n$-th) subtower\footnote{The hight of this subtower is $2^{n(n+1)/2}$.} with spacers.
In such a way we obtain a new rank-one 
$\Bbb Z$-action $T^*\sim(2^n-1,\sigma_n^*)_{n=1}^\infty$ such that  
$$
\sigma^*_n(i)=\cases 0,&\text{if }i=1,\dots,2^n-2\\ 2^{n(n+1)/2},&\text{if }i=2^n-1\endcases
$$
for each $n$.
By Proposition~A.2, $T^*$ is essentially $0$-expansive.
It is called a {\it constructive odometer} in \cite{AdFePe}.
By the proof of Theorem~A.3, $T^*$ is isomorphic to $T$.
\endexample

\enddocument